\documentclass[11pt,a4paper]{article}
\usepackage{graphicx}
\usepackage{amsmath}
\usepackage{amssymb}
\usepackage{amsthm}
\usepackage{geometry}
\usepackage{fancyhdr}
\usepackage{color} 
\usepackage{overpic}
\usepackage{mathabx}

\usepackage[all]{xy}
\newcommand{\al}{\alpha}

\newcommand{\be}{\beta}
\newcommand{\ga}{\gamma}
\newcommand{\Ga}{\Gamma}
\newcommand{\shm}{\underline{\rm SHM}}

\newcommand{\ra}{\rightarrow}

\newcommand{\xra}{\xrightarrow}

\newcommand{\rgl}{\rangle}
\newcommand{\lgl}{\langle}

\newcommand{\bpf}{\begin{proof}}
\newcommand{\epf}{\end{proof}}
\newcommand{\bthm}{\begin{thm}}
\newcommand{\ethm}{\end{thm}}
\newcommand{\bprop}{\begin{prop}}
\newcommand{\eprop}{\end{prop}}
\newcommand{\bcor}{\begin{cor}}
\newcommand{\ecor}{\end{cor}}
\newcommand{\blem}{\begin{lem}}
\newcommand{\elem}{\end{lem}}
\newcommand{\bdefn}{\begin{defn}}
\newcommand{\edefn}{\end{defn}}
\newcommand{\bexmp}{\begin{exmp}}
\newcommand{\eexmp}{\end{exmp}}
\newcommand{\brem}{\begin{rem}}
\newcommand{\erem}{\end{rem}}

\newcommand{\bdia}{\begin{displaymath}\xymatrix}
\newcommand{\edia}{\end{displaymath}}
\newcommand{\beq}{\begin{equation*}\begin{aligned}}
\newcommand{\eeq}{\end{aligned}\end{equation*}}

\newcommand{\intg}{\mathbb{Z}}
\newcommand{\real}{\mathbb{R}}

\newtheorem{thm}{\textbf {Theorem}}[section]
\newtheorem{cor}[thm]{\textbf{Corollary}}
\newtheorem{prop}[thm]{\textbf{Proposition}}
\newtheorem{lem}[thm]{\textbf{Lemma}}
\newtheorem{conj}[thm]{Conjecture}

\newtheorem{quest}[thm]{Question}

\theoremstyle{definition}
\newtheorem{defn}[thm]{\textbf{Definition}}

\newtheorem{exmp}[thm]{Example}

\theoremstyle{remark}
\newtheorem{rem}[thm]{Remark}

\newcommand{\khm}{\underline{\rm KHM}^-}
\newcommand{\khi}{\underline{\rm KHI}^-}
\newcommand{\shi}{\underline{\rm SHI}}

\title{Knot Floer homologies in monopole and instanton theory via sutures}

\author{Zhenkun Li}
\date{}

\begin{document}
\bibliographystyle{plain}

\maketitle
\abstract{In this paper we construct possible candidates for the minus versions of monopole and instanton knot Floer homologies. For a null-homologous knot $K\subset Y$ and a base point $p\in K$, we can associate the minus versions, $\khm(Y,K,p)$ and $\khi(Y,K,p)$, to the triple $(Y,K,p)$. We prove that a Seifert surface of $K$ induces a $\intg$-grading, and there is an $U$-map on the minus versions, which is of degree $-1$. We also prove other basic properties of them. If $K\subset Y$ is not null-homologous but represents a torsion class, then we can also construct the corresponding minus versions for $(Y,K,p)$. We also proved a surgery-type formula relating the minus versions of a knot $K$ with those of the dual knot, when performing a Dehn surgery of large enough slope along $K$. The techniques developed in this paper can also be applied to compute the sutured monopole and instanton Floer homologies of any sutured solid tori.}

\tableofcontents

\section{Introduction}
\subsection{Statement of results}
Floer homologies have become very important tools in the study of $3$-manifolds, since the first construction by Floer in \cite{floer1988instanton}. Among them, two major branches are the monopole Floer homology, which was introduced by Kronheimer and Mrowka \cite{kronheimer2007monopoles} and the Heegaard Floer homology, which was introduced by Ozsv\'ath and Szab\'o \cite{ozsvath2004holomorphic} or Rasmussen \cite{rasmussen2003floer}. For a closed oriented $3$-manifold $Y$, there are four flavors of homologies associated to $Y$ in each of the two theories, and they are isomorphic by work of Kutluhan, Lee and Taubes in \cite{kutluhan2010hf} and in subsequent papers. If there is a knot $K$ inside a $3$-manifold $Y$, then there are corresponding four flavors of homologies of the pair $(Y,K)$ in the Heegaard Floer theory. See Ozsva\'ath and Szab\'o \cite{ozsvath2004holomorphicknot}. However, some corresponding constructions in the monopole and instanton theory are missing. The only monopole or (non-singular) instanton Floer homology for knots in $3$-manifolds is a version based on sutured manifolds, which was introduced by Kronheimer and Mrowka in \cite{kronheimer2010knots} and was refined by Baldwin and Sivek in \cite{baldwin2015naturality}. The monopole version is proved to be isomorphic to the hat version of the knot Floer homology in Heegaard Floer theory, which was due to Baldwin and Sivek \cite{baldwin2016equivalence} or Lekili \cite{lekili2013heegaard}. In this paper, we construct Floer homologies associated to a based oriented null-homologous knot, which are candidates for the monopole and the instanton correspondences of the minus version of the knot Floer homology in Heegaard Floer theory.

\bthm\label{thm_1_1}
Suppose $Y$ is a closed connected oriented $3$-manifold and $K\subset Y$ is an oriented null-homologous knot. Suppose further that $S$ is a Seifert surface of $K$, and $p\in K$ is a base point. Then, we can associate the triple $(Y,K,p)$ a module $\khm(Y,K,p)$ over the mod $2$ Novikov Ring $\mathcal{R}$. It is well defined up to multiplication by a unit in $\mathcal{R}$. The Seifert surface $S$ induces a $\intg$ grading on $\khm(Y,K,p)$, which we denote by $\khm(Y,K,P,S,i)$. Moreover, the following properties hold.

(1) For $i>g=g(S)$, $\khm(Y,K,p,S,i)=0$.

(2) There is a map
$$U:\khm(Y,K,p)\ra\khm(Y,K,p)$$
that is of degree $-1$.

(3) There exists an $N_0\in\intg$ such that if $i<N_0$, then
$$U:\khm(Y,K,p,S,i)\cong \khm(Y,K,p,S,i-1).$$

(4) There exists an exact triangle
\begin{equation*}
\xymatrix{
{\khm}(Y,K,p)\ar[rr]^{U}&&{\khm}(Y,K,p)\ar[dl]^{\psi}\\
&\underline{\rm KHM}(Y,K,p)\ar[lu]^{\psi'}&
}	
\end{equation*}

(5) If $Y=S^3$ and $S$ realizes the genus of the knot, then we have
$$\khm(Y,K,p,S,i)\neq 0$$
for $i=g(S)$.
\ethm

The same construction can also be carried out in instanton theory. 
\bthm
Under the same settings as in Theorem \ref{thm_1_1}, we can construct $\khi(Y,K,p)$, using instanton Floer homology, so that all the properties (1)-(5) in the that theorem hold in the instanton settings.
\ethm

It is worth mentioning here that Kutluhan \cite{kutluhan2013seiberg} constructed another minus version of knot monopole Floer homology in a different way. He used the holonomy filtration for the construction.

\subsection{Outline of the proof}
In the current subsection, we only discuss in the monopole settings, and the constructions in the instanton settings are similar. The construction of $\khm(Y,K,p)$ is based on sutured monopole Floer homology. A sutured manifold $(M,\ga)$ is a compact oriented $3$ manifold with a closed oriented $1$-submanifold $\ga$ on $\partial{M}$, which we call the suture. The suture $\ga$ divides $\partial{M}$ into two parts, according to the orientations of $\ga$ and the $3$-manifold, which we call $R_-(\ga)$ and $R_+(\ga)$, respectively. Sutured manifolds were first introduced by Gabai in  \cite{gabai1983foliations}. Kronheimer and Mrowka then carried out the construction of the monopole and instanton Floer homologies on balanced sutured manifolds in \cite{kronheimer2010knots}. 

A sutured manifold $(M,\ga)$ is called balanced if $M$ and $R(\ga)$ both have no closed components and $\chi(R_-(\ga))=\chi(R_+(\ga))$. To define the sutured monopole Floer homology for such a pair $(M,\ga)$, Kronheimer and Mrowka constructed a closed $3$-manifold $Y$, together with a distinguishing surface $R$, out of $(M,\ga)$. The pair $(Y,R)$ is called a closure of $(M,\ga)$. Sometimes we simply call $Y$ a closure. The genus of the closure refers to the genus of the surface $R$. To construct a closure, one needs to find a compact connected oriented surface $T$, whose boundary is diffeomorphic to $\ga$, and then glue $[-1,1]\times T$ to $M$, with $[-1,1]\times \partial{T}$ identified with an annular neighborhood of $\ga\subset \partial{M}$. The surface $T$ is called an auxiliary surface. The new $3$-manifold after the gluing is called a pre-closure, and it has two boundary components, $R_+$ and $R_-$, of the same genus. Then, we can pick a diffeomorphism $h$ from $R_+$ to $R_-$ to glue the two boundary components together to obtain a closure $(Y,R)$. We call $h$ a gluing diffeomorphism.

To study the naturality of sutured monopole Floer homology, Baldwin and Sivek \cite{baldwin2014invariants} constructed  canonical maps between two different closures of a same balanced sutured manifold $(M,\ga)$. Their construction is only well-defined up to multiplication by a unit, so the closures and canonical maps form a projective transitive system and result in a canonical module $\shm(M,\ga)$, whose elements are well defined only up to a unit.

The construction of the (canonical) module $\khm(Y,K,p)$ was inspired by Etnyre, Vela-Vick and Zarev in \cite{etnyre2017sutured}, where they use a sequence of balanced sutured manifolds $(Y(K),\Ga_n)$ and the gluing maps in sutured (Heegaard) Floer theory, which was introduced by Honda, Kazez, and Mati\'c \cite{honda2008contact}, to construct a direct system. They proved that the direct limit is isomorphic to the classical minus version of knot Floer homology in Heegaard Floer theory. Here, $Y(K)=Y\backslash {\rm int}(N(K))$ is the knot complement, and $\Gamma_n$ consists of two curves on $\partial{Y(K)}\cong T^2$, which are of class $\pm(1,-n)$ under the framing induced by some Seifert surface. In this paper, we construct the same direct system in sutured monopole Floer theory. In particular, there is a commutative diagram
\begin{equation}\label{eq_commutative_diagram}
\xymatrix{
\shm(-Y(K),-\Ga_n)\ar[rr]^{\psi_{-,n+1}^n}\ar[dd]^{\psi_{+,n+1}^n}&&\shm(-Y(K),-\Ga_{n+1})\ar[dd]^{\psi_{+,n+2}^{n+1}}\\
&&\\
\shm(-Y(K),-\Ga_{n+1})\ar[rr]^{\psi_{-,n+2}^{n+1}}&&\shm(-Y(K),-\Ga_{n+2})
}	
\end{equation}
Here, the balanced sutured manifolds are the same as described above, and the maps come from gluing maps in sutured monopole monopole Floer theory, which were constructed by the author in \cite{li2018gluing}.

The commutativity of (\ref{eq_commutative_diagram}) is guaranteed by the functoriality of the gluing map. The crucial difference from the work of Etnyre, Vela-Vick and Zarev in \cite{etnyre2017sutured} is that, because of the involvement of closures, the construction of the grading in the monopole and the instanton settings is a delicate issue. We construct a grading in the direct limit in two steps.

The first step is to construct a grading on each $\shm(Y(K),\Ga_n)$, for all $n$, using the Seifert surface $S$. To construct such a grading, we work with a more general case, where $(M,\ga)$ is an arbitrary balanced sutured manifold, $S$ is a properly embedded surface whose boundary is connected, and $\partial{S}$ intersects $\ga$ transversely at $2n$ points. 

For the case $n=1$, the construction has already been carried out by Baldwin and Sivek in \cite{baldwin2018khovanov}. When $n=1$, we can pick a properly embedded arc $\al\subset T$, where $T$ is an auxiliary surface for $(M,\ga)$. When gluing $[-1,1]\times T$ to $M$, we require that the end points of $\al$ are glued to the two intersection points in $\partial{S}\cap\ga$, and, hence, $[-1,1]\times\al$ is glued to $S$ along $[-1,1]\times\partial{\al}$. Then, $S$ becomes a surface $\widetilde{S}$ properly embedded in the pre-closure $\widetilde{M}$. Note $\widetilde{M}$ has two boundary components $R_{+}$ and $R_-$, and the two boundary components of $\widetilde{S}$ are contained in different boundary components of $\widetilde{M}$. Thus, we can pick a gluing diffeomorphism $h:R_+\ra R_-$ which also identifies the two boundary components of $\widetilde{S}$. Hence, $\widetilde{S}$ becomes a closed surface $\bar{S}$ inside the closure $Y$ of $(M,\ga)$. The grading can be defined by looking at the pairing of the first Chern classes of the spin${}^c$ structures on $Y$ with the fundamental class of $\bar{S}$. This idea was first introduced by Kronheimer and Mrowka in \cite{kronheimer2010knots}, and, in \cite{baldwin2018khovanov}, Baldwin and Sivek proved that, when $n=1$, the definition of the grading is independent of all choices made in the construction and is well defined in $\shm(M,\ga)$.

For a general $n$, the basic idea to construct a grading is the same. However, there are more choices involved, and, thus, many new issues arise. For example, for a general $n$, we need to pick $n$ arcs $\al_1,...,\al_n$ instead of just one, and we need to specify which arc connects which pair of intersection points in $\partial{S}\cap\ga$. Thus, this leads to a new combinatorial problem which did not occur in Baldwin and Sivek \cite{baldwin2018khovanov}. We deal with this combinatorial problem in Subsection \ref{subsec_choice_of_type_I}. To conclude the proof, we also need a new interpretation of Baldwin and Sivek's canonical maps between different closures. We use simply the Floer excision introduced by Kronheimer and Mrowka in \cite{kronheimer2010knots} to construct an equivalent canonical map, which was originally introduced by Baldwin and Sivek in \cite{baldwin2015naturality}. This is covered in Subsection \ref{subsec_choices_of_type_III}.

When constructing the grading based on a surface $S$, we need the extra assumption that $n$ is odd. Recall that $|S\cap\ga|=2n$. If $n$ is even, then we need to perturb $S$ to create a new pair of intersection points and, thus, increase $n$ by $1$. There are two different ways of perturbations, which we call positive and negative stabilizations, and denote them by $S^+$ and $S^-$, respectively. Based on $S^+$ and $S^-$, we can construct two different gradings on $\shm(Y(K),\Ga_n)$. The relation between the two gradings will be the key to the second step of constructing a grading on the direct limit. Also, using the grading shifting property betweem $S^{+}$ and $S^-$, we can compute the sutured monopole Floer homology of a solid torus with any valid suture.

\bprop
Suppose $V$ is a solid torus and $\ga$ is a suture on $\partial{V}$ with $2n$ components and slope $\frac{p}{q}$, then
$$\shm(-V,-\ga)\cong\mathcal{R}^{(2^{n-1}\cdot |p|)}.$$
\eprop

Similarly, in instanton theory, we have the following.
\bprop
Suppose $V$ is a solid torus and $\ga$ is a suture on $\partial{V}$ with $2n$ components and slope $\frac{p}{q}$, then
$$\shi(-V,-\ga)\cong\mathbb{C}^{(2^{n-1}\cdot |p|)}.$$
\eprop

The second step of constructing a grading on the direct limit is to prove that maps in the commutative diagram (\ref{eq_commutative_diagram}) shift the grading in a desired way. To be more explicit, $\psi_{-,n+1}^n$ must be of degree $0$, while $\psi_{+,n+1}^n$ must be of degree $-1$. The construction of the maps $\psi_{\pm,n+1}^n$ relies on the by-pass attachments in the monopole and instanton settings, which are realized by contact handle attachments, as introduced by Baldwin and Sivek in \cite{baldwin2016contact,baldwin2016instanton}. 

It is a basic observation that the region we attach contact handles is disjoint from the Seifert surface $S$, hence if we look at the grading associated to the 'correct' surfaces, then $\psi_{-,n+1}^n$ and $\psi_{+,n+1}^n$ will both preserve the grading. However, the 'correct' surfaces involves both positive and negative stabilizations, while, to define a canonical grading on $\shm(Y(K),\Ga_n)$, we only use negative stabilizations. Hence, the problem is reduced to understanding the grading shifting between $S^{+}$ and $S^-$. 

To understand this grading shifting property, we first need a better understanding of the construction of the closures, the construction of canonical maps, and how spin${}^c$ structures on different closures are related by canonical maps. In particular, we prove the following result. 

\bprop\label{prop_1_2}
Suppose $(Y(K),\Ga_n)$ is the balanced sutured manifold described as above, and $Y_n$ is a closure of $(Y(K),\Ga_n)$. Suppose $\mathfrak{s}_1$ and $\mathfrak{s}_2$ are two spin${}^c$ structures on $Y_n$, so that they both support the sutured monopole Floer homology of $(Y(K),\Ga_n)$. Then, in terms of Poinc\'are duals of first Chern classes of the spin${}^c$ structures, the difference between $\mathfrak{s}_1$ and $\mathfrak{s}_2$ lies in $H_1(Y(K))$. More precisely, there is a $1$-cycle $x$ in $Y(K)$ so that
$$P.D.(c_1(\mathfrak{s}_1)-c_1(\mathfrak{s}_2))=[x]\in H_1(Y).$$
\eprop

Proposition \ref{prop_1_2} will be the basis for understanding the grading shifting property between the gradings associated to $S^+$ and $S^-$, which are the positive and negative stabilizations of $S$. We deal with the grading shifting property in Section \ref{sec_grading_shift}. We present the construction of the minus version in Subsection \ref{subsec_construction} and prove some basic properties of it in Subsection \ref{subsec_properties}. Most of the basic properties have been stated in theorem \ref{thm_1_1}. Besides them, we also prove that the direct system in the construction of the minus version stabilizes.

\bprop
For a fixed $i\in\intg$, there exists $N_1\in\intg$, such that for $n>N_1$, we have an isomorphism:
$$\psi_{-,n+1}^n:\shm(-Y(K),-\Ga_n,i)\cong \shm(-Y(K),-\Ga_{n+1},i).$$
\eprop

The techniques used in computing the sutured Floer homology of a solid torus can also be applied to knot complements. As a result, we obtain the following.

\bprop
Suppose $K\subset Y$ is a knot and $S\subset Y$ is a Seifert surface of $K$. Suppose $Y_{\phi}$ is the manifold obtain from $Y$ by doing a Dehn surgery along $K$ with slope $-\frac{p}{q}$ with $p,q>0$. We also have the dual knot $K_{\phi}\subset Y$. Then for any fixed $i$, there exists $N\in\real$, such that if the surgery slope $-\frac{p}{q}<N$, then we have
$$\khm(-Y_{\phi},K_{\phi},S,i)\cong\khm(-Y,K,S,i).$$

Moreover, a similar result in instanton theory also holds.
\eprop

{\bf Acknowledgements.} This material is based upon work supported by the National Science Foundation under Grant No. 1808794. The author would like to thank his advisor Tom Mrowka for his invaluable helps. The author would like to thank John Baldwin, Mariano Echeverria, Jianfeng Lin, Langte Ma, and Donghao Wang, Yi Xie for helpful conversations.

\section{Preliminaries}
\subsection{Balanced sutured manifolds and monopole Floer homology}\label{subsec_sutured_monopoles_and_natruality}
\bdefn\label{defn_balanced_sutured_manifold}
A {\it balanced sutured manifold} is a pair $(M,\ga)$ of a compact oriented $3$-manifold $M$ and a closed oriented $1$-submanifold $\ga\subset\partial M$. On $\partial{M}$, let $A(\ga)=\ga\times[-1,1]$ be an annular neighborhood of $\ga$, and let 
$$R(\ga)=\partial{M}\backslash{\rm int}(A(\ga)).$$
They satisfy the following requirements.

(1) Both $M$ and $R(\ga)$ have no closed components.

(2) If we orient $\partial{R}(\ga)=\partial{A}(\ga)=\ga\times\{\pm1\}$ in the same way as $\ga$, then the orientation on $\partial{R}(\ga)$ must induce a unique orientation on $R(\ga)$. This orientation is called the {\it canonical orientation} on $R(\ga)$. Use $R_+(\ga)$ to denote the part of $R(\ga)$ whose canonical orientation coincides with the boundary orientation of $\partial{M}$, and $R_-(\ga)$ the rest. 

(3). $\chi(R_+(\ga))=\chi(R_-(\ga)).$
\edefn

To define sutured monopole Floer homology, we need to construct a closed $3$-manifold out of a balanced sutured manifold $(M,\ga)$. Let $T$ be a connected oriented surface so that the following holds.

(1) There is an orientation reversing diffeomorphism
$$f:\partial{T}\ra \ga.$$

(2) $T$ has genus at least $2$.

After choosing such a $T$, we can use $f$ to glue a thickened $T$ to $M$:
$$\widetilde{M}=M\mathop{\cup}_{f}[-1,1]\times T.$$
The manifold $\widetilde{M}$ has two boundary components:
$$\partial\widetilde{M}=R_+\cup R_-,$$
where
$$R_{\pm}=R_{\pm}(\ga)\mathop{\cup}_{f} \{\pm 1\}\times T.$$
Let $h:R_+\ra R_-$ be an orientation preserving diffeomorphism, then we can form a closed $3$-manifold as
$$Y=\widetilde{M}\mathop{\cup}_{id\cup h}[-1,1]\times R_{+},$$
where $h:\{1\}\times R_{+}\ra R_-\subset \partial{\widetilde{M}}$ is the map just defined and $id:\{-1\}\times R_+\ra R_+\subset \partial\widetilde{M}$ is the identity on $R_+$. Let $R=\{0\}\times R_+\subset Y$, and we make the following definition.

\bdefn\label{defn_closure}
The manifold $\widetilde{M}$ is called a {\it pre-closure} of $(M,\ga)$. The pair $(Y,R)$ is called a {\it closure} of $(M,\ga)$. The choices $T,f,c$, and $h$ are called the auxiliary data. In particular, the surface $T$ is called an {\it auxiliary surface} and $h$ is a {\it gluing diffeomorphism}.
\edefn

\brem
Throughout this paper, we   require that $T$ is connected and has large enough genus. However, in general, the choice of auxiliary surface has more freedoms. See Kronheimer and Mrowka \cite{kronheimer2010knots}.
\erem

To construct local coefficients, we also need to choose a non-separating simple closed curve $\eta\subset R$. The base ring we use in the present paper is the mod 2 Novikov ring $\mathcal{R}$. For a detailed definition, readers are referred to \cite{baldwin2016contact}.

\bdefn\label{defn_set_of_spin_c_structures}
Suppose $Y$ is a closed connected oriented $3$-manifold and $R$ is a closed oriented surface inside $Y$, so that each component of $R$ has genus at least $2$. If $R$ is connected, we define the set of {\it top} spin${}^c$ structures as follows:
$$\mathfrak{S}(Y|R)=\{\rm{spin}^c~\rm{structure}~\mathfrak{s}~\rm{on}~Y|c_1(\mathfrak{s})[R]=2g(R)-2.\}$$

If $R$ is disconnected and let $R_1,...,R_n$ be its components, then we define
$$\mathfrak{S}(Y|R)=\mathop{\bigcap}_{i=1}^n\mathfrak{S}(Y|R_i).$$

For later references, we also define the set of {\it supporting} spin${}^c$ structures as follows: 
$$\mathfrak{S}^*(Y|R)=\{\mathfrak{s}\in\mathfrak{S}(Y|R)|\widecheck{HM}_{\bullet}(Y,\mathfrak{s};\Ga_{\eta})\neq0\}.$$
Here, $\widecheck{HM}_{\bullet}(Y,\mathfrak{s};\Ga_{\eta})$ is the to-version of monopole Floer homology with local coefficients associated to the pair $(Y,\mathfrak{s})$. more details, readers are referred to \cite{kronheimer2007monopoles}.
\edefn

\bdefn
The {\it sutured monopole Floer homology} of $(M,\ga)$ is defined to be
$$SHM(M,\ga)=HM(Y|R;\Gamma_{\eta}),$$
where
$$HM(Y|R;\Gamma_{\eta})=\bigoplus_{\mathfrak{s}\in\mathfrak{S}(Y|R)}\widecheck{HM}_{\bullet}(Y,\mathfrak{s};\Gamma_{\eta})$$
\edefn

The following lemmas from Kronheimer and Mrowka \cite{kronheimer2010knots} will be useful.

\blem\label{lem_surface_bundle_over_S_1}
Suppose $Y$ is a surface bundle over $S^1$ whose fibres are closed connected oriented surfaces of genus at least $2$. Let $R$ be a fibre and $\eta\subset R$ be a non-separating simple closed curve. Then, there is a unique spin${}^c$ structure $\mathfrak{s}$ on $Y$ so that the following is true.

(1) We have $c_1(\mathfrak{s})[R]=2g(R)-2$.

(2) We have $\widecheck{HM}_{\bullet}(Y,\mathfrak{s};\Ga_{\eta})\neq0$.

Moreover, for this spin${}^c$ structure $\mathfrak{s}$, we have
$$\widecheck{HM}_{\bullet}(Y,\mathfrak{s};\Ga_{\eta})\cong\mathcal{R},$$
where $\mathcal{R}$ is the base ring for local coefficients.
\elem

\blem\label{lem_adjunction_inequality}
Suppose $Y$ is a closed oriented $3$-manifold and $R\subset Y$ is an embedded closed connected oriented surface of genus at least $1$. Suppose further that $\mathfrak{s}$ is a spin${}^c$ structure such that
$$|c_1(\mathfrak{s})[R]|>2g(R)-2,$$
then we have
$$\widecheck{HM}_{\bullet}(Y,\mathfrak{s};\Ga_{\eta})=0$$
for any choice of local coefficients.
\elem

Floer excisions were introduced to the context of sutured monopole Floer homology by Kronheimer and Mrowka in \cite{kronheimer2010knots}. In the rest of the current subsection, we summarize the results that we need in later sections.

For $i=1,2$, suppose $Y_i$ is a closed connected oreinted $3$-manifold and $R_i\subset Y_i$ is an embedded closed connected oriented homologically essential surface of genus at least $2$. Let $\eta_i\subset R_i$ be a non-separating simple closed curve. When cutting $Y_i$ open along $R_i$, we get
$$\widetilde{Y}_i=Y_i\backslash{\rm int}(N(R_i)),$$
where $N(R_i)$ is a product neighborhood of $R_i\subset Y_i$. The manifold $\widetilde{Y}_i$ has two boundary components
$$\partial\widetilde{Y}_i=R_{i,+}\cup R_{i,-}.$$
We orient $R_{i,\pm}$ in the same way as $R_{i}$. There are parallel copies of $\eta_i$, which we call $\eta_{i,\pm}$, on the surfaces $R_{i,\pm}$. Pick an orientation preserving diffeomorphism 
$$h:R_{1}\ra R_{2}$$
so that $h(\eta_1)=\eta_2$. We can use $h$ to glue $R_{1,+}$ to $R_{2,-}$ and $R_{1,-}$ to $R_{2,+}$. Then, $\widetilde{Y}_1$ and $\widetilde{Y}_2$ are glued together to become a connected $3$-manifold which we call $Y$. Let $R\subset Y$ be the disjoint union of the surfaces $R_{1,+}$ and $R_{2,+}$ in $Y$. Let $\eta\subset R$ be the disjoint union of curves $\eta_{1,+}$ and $\eta_{2,+}$.

\begin{figure}[h]
\centering
\begin{overpic}[width=5.0in]{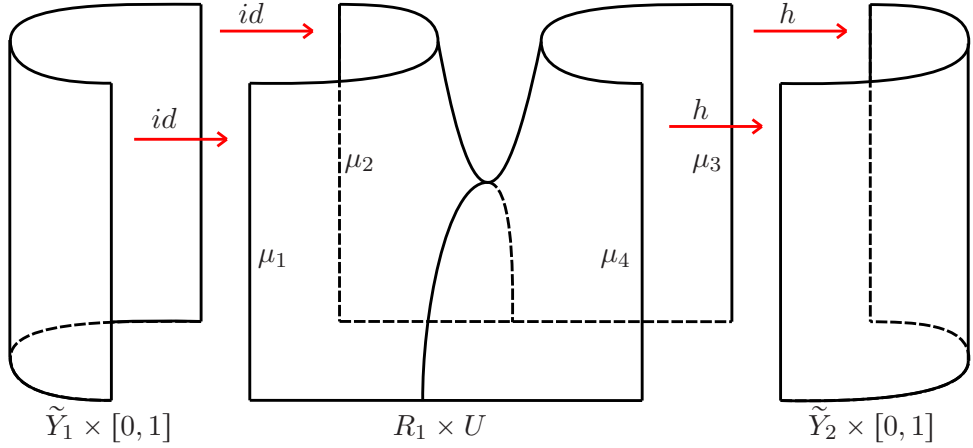}
	\put(4,-3){$\widetilde{Y}_1\times[0,1]$}
	\put(40,-3){$R_{1}\times U$}
	\put(83,-3){$\widetilde{Y}_2\times[0,1]$}
	\put(15,29){$id$}
	\put(24,40){$id$}
	\put(71,30){$h$}
	\put(80,40){$h$}
	\put(35,25){$\mu_2$}
	\put(26,15){$\mu_1$}
	\put(71,25){$\mu_3$}
	\put(61.5,15){$\mu_4$}
\end{overpic}
\vspace{0.05in}
\caption{Gluing three parts together to get $W$. The middle part is $R_{1}\times U$, while the $R_{1,+}$ directions shrink to a point in the figure.}\label{fig_excision_cobordism}
\end{figure}

There is a $4$-dimensional cobordism $W$ from $Y_1\sqcup Y_2$ to $Y$, which is constructed as follows: Let $U$ be the surface as depicted in Figure \ref{fig_excision_cobordism}. It has four vertical arcs as part of the boundary, and we can assume that each of them is identified with $[0,1]$. Now we can use the identity map and the map $h$ to glue three pieces $\widetilde{Y}_1$, $\widetilde{Y}_2$ and $U\times R_1$ together, to obtain the desired cobordism $W$. This cobordism $W$ then induces a map as in \cite{kronheimer2010knots}
\begin{equation}\label{eq_floer_excision_higer_genus}
HM(W): HM(Y_1\sqcup Y_2|R_1\cup R_2;\Gamma_{\eta_1\cup\eta_2})\ra HM(Y|R;\Ga_{\eta}).
\end{equation}

We can also cut and re-glue along tori. For $i=1,2$, let $Y_i$ be as above. Let $T_i\subset Y_i$ be a torus and $R_i\subset Y_i$ be a closed connected oriented surface so that $R_i$ intersects $T_i$ transversely along a circle $c_i$. Suppose $\eta_i\subset R_i$ is a simple closed curve so that $\eta_i$ intersects $c_i$ transversely at a point $p_i$. Let 
$$h:T_1\ra T_2$$
be an orientation preserving diffeomorphism so that $h(c_1)=c_2$ and $h(p_1)=p_2$. As above, we can cut $Y_i$ open along $T_i$ and re-glue using $h$ to obtain a connected $3$-manifold $Y$. There is a distinguishing surface $R$, obtained by cutting $R_i$ open along $c_i$ and re-glue using $h$. The curves $\eta_1$ and $\eta_2$ are also cut and re-glued to give rise to a simple closed curve $\eta\subset R\subset Y$. As above, there is a cobordism map
\begin{equation}\label{eq_floer_excision_torus}
HM(W): HM(Y_1\sqcup Y_2|R_1\cup R_2;\Gamma_{\eta_1\cup\eta_2})\ra HM(Y|R;\Ga_{\eta}).
\end{equation}

\bthm[Kronheimer and Mrowka \cite{kronheimer2010knots}]\label{thm_floer_excision} The maps (\ref{eq_floer_excision_higer_genus}) and (\ref{eq_floer_excision_torus}) are both isomorphisms.
\ethm

\subsection{The naturality of sutured monopole Floer homology}\label{subsec_naturality}

In \cite{baldwin2014invariants}, Baldwin and Sivek constructed a canonical map between two different closures of the same balanced sutured manifold. To do this, they also refined the definition of closures.

\bdefn\label{defn_marked_closure}
A {\it marked closure} $\mathcal{D}=(Y,R,r,m,\eta)$ of a balanced sutured manifold $(M,\ga)$ consists of the following.

(1) A closed connected oriented $3$-manifold $Y$.

(2) A closed connected oriented surface $R$ of genus at least two.

(3) An orientation preserving embedding
$$r:R\times[-1,1]\hookrightarrow Y.$$

(4) An orientation preserving embedding
$$m:M\hookrightarrow Y\backslash{\rm int}({\rm im}(r)).$$

(5) A non-separating simple closed curve $\eta\subset R$.

They satisfy the following requirements.

(a) The embedding $m$ extends to a diffeomorphism
$$M\mathop{\cup}_{f}T\times[-1,1]\ra Y\backslash{\rm int}({\rm im}(r)),$$
for some auxiliary data $(T,f)$.

(b) The embedding $m$ restricts to an orientation preserving embedding
$$R_{+}(\ga)\hookrightarrow r(R\times\{-1\}).$$

The {\it genus} of the marked closure $\mathcal{D}$ is referred to the genus of the surface $R$. We define
$$SHM(\mathcal{D})=\bigoplus_{\mathfrak{s}\in\mathfrak{S}(Y|r(R\times\{0\}))}\widecheck{HM}_{\bullet}(Y,\mathfrak{s};\Gamma_{r(\eta\times\{0\})}).$$
\edefn

\bthm[Baldwin and Sivek \cite{baldwin2015naturality}]\label{thm_canonical_maps}
Suppose $(M,\ga)$ is a balanced sutured manifold, then for any two marked closures $\mathcal{D}_1$ and $\mathcal{D}_2$ of $(M,\ga)$, there is a canonical map $\Phi_{\mathcal{D}_1,\mathcal{D}_2}$, well defined up to a unit, from $SHM(\mathcal{D}_1)$ to $SHM(\mathcal{D}_2)$. The canonical maps satisfy following properties.

(1) If $\mathcal{D}_1=\mathcal{D}_2$, then
$$\Phi_{\mathcal{D}_1,\mathcal{D}_2}\doteq id.$$
Here $\doteq$ means equal up multiplication by a unit.

(2) Suppose there is a third marked closure $\mathcal{D}_3$ for $(M,\ga)$, then we have
$$\Phi_{\mathcal{D}_1,\mathcal{D}_3}\doteq \Phi_{\mathcal{D}_2,\mathcal{D}_3}\circ \Phi_{\mathcal{D}_1,\mathcal{D}_2}.$$
\ethm

Hence, for a balanced sutured manifold $(M,\ga)$, marked closures $\mathcal{D}$ and canonical maps $\Phi$ fits into a projective transitive system, which is defined in \cite{baldwin2015naturality}. The projective system determines a canonical module, which we   denote by
$$\shm(M,\ga).$$
We can then talk about elements (up to multiplication by a unit) in that canonical module.

There is an extra ambiguity when dealing with knots in $3$-manifolds. Let $K\subset Y$ be a knot. The extra ambiguity comes from the choices of tubular neighborhoods of $K\subset Y$ to remove to obtain a knot complement. Fix a point $p\in K$. Suppose
$$\varphi:S^1\times D^2\hookrightarrow Y$$
is an embedding, where $D^2$ is the unit sphere in the complex plane, and $S^1=\partial{D}^2$. We require that
$$\varphi(S^1\times\{0\})=K,~{\rm and}~\varphi(\{1\}\times\{0\})=p.$$

Let $Y_{\varphi}=Y\backslash{\rm int}({\rm im}(\phi))$, and let $\ga_{\varphi}=\varphi(\{\pm1\}\times\partial{D}^2)$, with opposite orientations on two components. For each fixed $\varphi$, we have a well defined canonical module $\shm(Y(\varphi),\ga_{\varphi})$, and we want also relate different choices of $\varphi$.

Suppose $\varphi'$ is another embedding $S^1\times D^2\hookrightarrow Y$, satisfying the same conditions as $\varphi$. Pick a tubular neighborhood $N$ of $K\subset Y$ such that ${\rm im}(\varphi),{\rm im}(\varphi')\subset N$. Also, pick an ambient isotopy
$$f_t:Y\ra Y,~t\in[0,1]$$
such that the following is true.

(1) For any $t\in[0,1]$, $f_t(p)=p.$

(2) For any $t\in[0,1]$, $f_t$ restricts to identity outside $N\subset Y$.

(3) We have $f_1({\rm im}(\varphi))={\rm im}(\varphi')$.

(4) We have $f_1(\varphi(\{\pm1\}\times\partial{D}^2))=\varphi'(\{\pm1\}\times D^2)$.

It is clear that $f_1:(Y_{\varphi},\ga_{\varphi})\ra (Y_{\varphi'},\ga_{\varphi'})$ is a diffeomorphism between balanced sutured manifolds. Hence, we can define
$$\Psi_{\varphi,\varphi'}=\shm(f_1):\shm(Y_{\varphi},\ga_{\varphi})\ra \shm(Y_{\varphi'},\ga_{\varphi'}).$$

\bthm(Baldwin and Sivek \cite{baldwin2015naturality})
The map $\Psi_{\varphi,\varphi'}$ is well defined, i.e., is independent of the choices of the tubular neighborhood $N$ and the ambient isotopy $f_t$. Also, it has the following properties.

(1) We have $\Psi_{\varphi,\varphi}=id.$

(2) If there is a third embedding $\varphi''$, then
$$\Psi_{\varphi,\varphi''}=\Psi_{\varphi',\varphi''}\circ \Psi_{\varphi,\varphi'}.$$
\ethm

Thus, we know that $\{\shm(Y_{\varphi},\ga_{\varphi})\}$ and $\{\Psi_{\varphi,\varphi'}\}$ form a transitive system of projective transitive systems. Thus, they lead to a larger projective transitive system, and, hence, the monopole knot Floer homology $\underline{\rm KHM}(Y,K,p)$ is well defined (as a projective transitive system).

\subsection{Contact structures and contact elements}\label{subsec_contact_elements_and_contact_structures}
In this subsection we summarize the results related to contact geometry which we will use in later sections.
\bdefn\label{defn_contact_sutured_manifold}
A {\it contact sutured manifold} $(M,\ga,\xi)$ is a triple where $(M,\ga)$ is a balanced sutured manifold, and $\xi$ is a contact structure on $(M,\ga)$ so that $\partial{M}$ is convex and $\ga$ is the dividing set. The contact structure $\xi$ is said to be {\it compatible} with the balanced sutured manifold $(M,\ga)$.
\edefn

\bthm(Baldwin and Sivek \cite{baldwin2016contact})
Suppose $(M,\ga,\xi)$ is a contact sutured manifold, then we can associate an element
$$\phi_{\xi}\in\shm(-M,-\ga)$$
to it. This element is called the {\bf contact element}.
\ethm

\bdefn
Suppose $(M',\ga')$ is a balanced sutured manifold. A {\it sutured submanifold} $(M,\ga)$ of $(M',\ga')$ is another balanced sutured manifold so that $M\subset{\rm int}(M')$. 
\edefn

The gluing maps in sutured monopole Floer homology were define by the author in \cite{li2018gluing}, and it is crucial in the construction of the direct system in Section \ref{seciton_direct_system}.

\bthm\label{thm_gluing_map}
Suppose $(M,\ga)$ is a sutured submanifold of $(M',\ga')$ and suppose $Z=M'\backslash{\rm int}(M)$. Suppose $\xi$ is a contact structure on $Z$ so that $(Z,\ga\cup\ga',\xi)$ is a contact sutured manifold. Then, there is a map
$$\Phi_{\xi}:\shm(-M,-\ga)\ra\shm(-M',-\ga'),$$
so that the following is true.

(1) If $(M',\ga')$ is a sutured submanifold of $(M'',\ga'')$ and there is a contact structure $\xi'$ on $M''\backslash{\rm int}(M')$, making it a contact sutured manifold, then we have
$$\Phi_{\xi'}\circ\Phi_{\xi}=\Phi_{\xi\cup\xi'}:\shm(-M,-\ga)\ra\shm(-M'',-\ga'').$$

(2) Suppose $(M',\ga',\xi')$ is a contact sutured manifold and $\xi'|_Z=\xi,$ then we have 
$$\Phi_{\xi}(\phi_{\xi'|_M})=\phi_{\xi'}.$$
\ethm

\begin{figure}[h]
\centering
\begin{overpic}[width=4.5in]{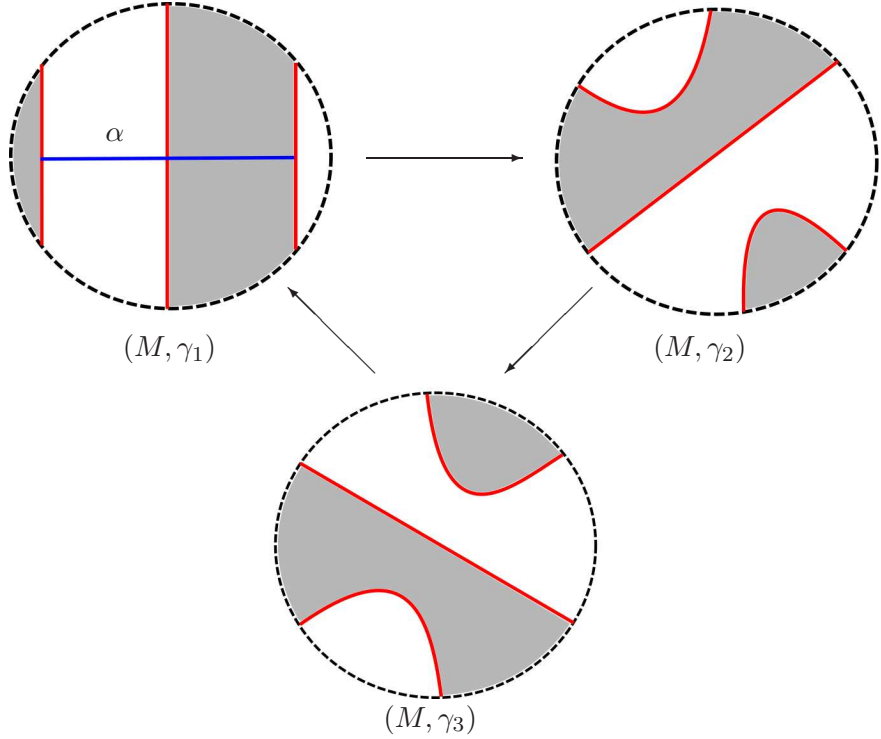}
	\put(43,-3){$(M,\ga_3)$}
	\put(13,40){$(M,\ga_1)$}
	\put(74,40){$(M,\ga_2)$}
	\put(11,65){$\al$}
	\put(41,63){\vector(1,0){18}}
	\put(67,48){\vector(-1,-1){10}}
	\put(42,38){\vector(-1,1){10}}
\end{overpic}
\vspace{0.05in}
\caption{The by-pass exact triangle.}\label{fig_by_pass}
\end{figure}

Suppose we have three balanced sutured manifold $(M,\ga_1)$, $(M,\ga_2)$, and $(M,\ga_3)$ so that the underlining $3$-manifold is the same, but the sutures are different. Suppose further that $\ga_1$, $\ga_2$, and $\ga_3$ are the same outside a disk $D\subset \partial{M}$, and, within the disk $D$, they are depicted as in Figure \ref{fig_by_pass}. We say that $(M,\ga_2)$ is obtained from $(M,\ga_1)$ by a by-pass attachment along the arc $\al$. Similarly, $(M,\ga_3)$ is obtained from a by-pass attachment from $(M,\ga_2)$ and $(M,\ga_1)$ from $(M,\ga_3)$. Then, we have the following theorem.

\bthm[Baldwin and Sivek \cite{baldwin2016contact}]\label{thm_by_pass_attachment}
There is an exact triangle relating the sutured monopole Floer homologies of the three balanced sutured manifolds:

\begin{equation}\label{eq_by_pass_general}
\xymatrix{
\shm(-M,-\ga_1)\ar[rr]^{\psi_{12}}&&\shm(-M,-\ga_2)\ar[dl]^{\psi_{23}}\\
&\shm(-M,-\ga_3)\ar[ul]^{\psi_{31}}&
}
\end{equation}
\ethm

In contact geometry, a by-pass is a half disk, which carries some particular contact structure, attached along a Legendrian arc to a convex surface. For more details, see Honda \cite{honda2000classification}. We can describe the maps in (\ref{eq_by_pass_general}) as follows: We explain the construction of the map $\psi_{12}$, and the other two are the same. Let $Z=\partial{M}\times[0,1]$, and we can pick the suture $\ga_1$ on $\partial{M}\times\{0\}$ as well as the suture $\ga_2$ on $\partial{M}\times\{1\}$. Then, there is a special contact structure $\xi_{12}$ on $Z$ that corresponds to the by-pass attachment and makes $(Z,\ga_1\cup\ga_2)$ a contact sutured manifold. Hence, we can attach $Z$ to $M$ by the identification $\partial{M}\times\{0\}=\partial{M}\subset M$. The result $(M\cup Z,\ga_2)$ is diffeomorphic to $(M,\ga_2)$ and we have
$$\psi_{12}=\Phi_{\xi_{12}}.$$
Here, $\Phi_{\xi_{12}}$ is the gluing map associated to $\xi_{12}$ as in Theorem \ref{thm_gluing_map}.

In Section \ref{seciton_direct_system}, we will use the by-passes on knot complements to construct the direct system. Let $K\subset Y$ be an oriented knot. Let $\lambda$ and $\mu$ be the longitude and meridian according to some framing of the knot. Let $\Ga_n$ be a suture on $\partial{Y(K)}$ which consists of two curves of class $\pm(\lambda-n\mu)$, and $\Ga_{\infty}$ consists of two meridians. In this case, $\partial{Y(K)}$ is a torus, and we have the following theorem due to Honda \cite{honda2000classification}.

\bthm\label{thm_contact_structures_on_T_2_cross_I}
There are two tight and minimal-twisting contact structures on $T^2\times[0,1]$ so that, for $i=1,2$, $T^2\times\{i\}$ is convex with dividing set being $\Ga_{n+i}$. These two contact structures correspond to two different by-pass attachments on $(Y(K),\Ga_n)$.
\ethm

\begin{figure}[h]
\centering
\begin{overpic}[width=5in]{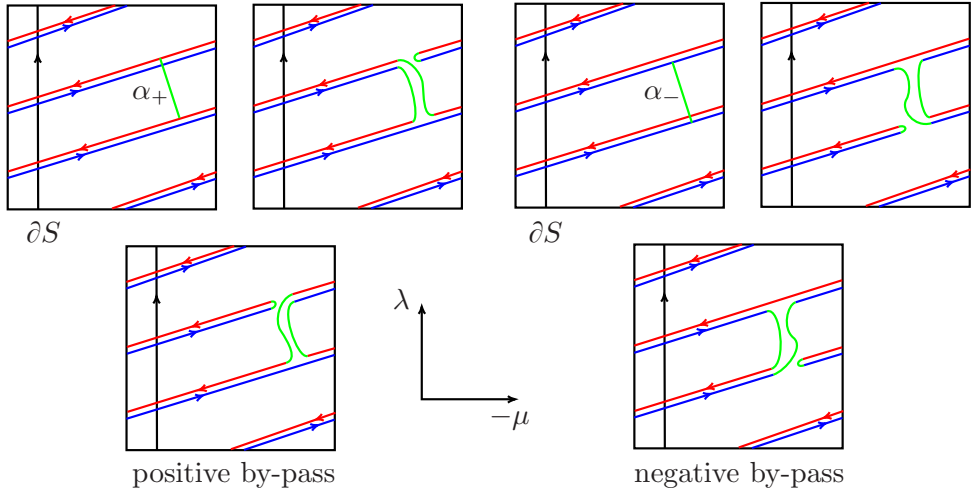}
	\put(13,-3){positive by-pass}
	\put(65,-3){negative by-pass}
	\put(13,37){$\al_{+}$}
	\put(66,37){$\al_{-}$}
	\put(40,15){$\lambda$}
	\put(2,22){$\partial{S}$}
	\put(54,22){$\partial{S}$}
	\put(50,3){$-\mu$}
\end{overpic}
\vspace{0.05in}
\caption{The positive and negative by-pass attachments for $(Y(K),\Ga_3))$. The squares represent the toroidal boundary of $Y(K)$. Note the contact structures $\xi_{\pm,2}$ correspond to the by-passes from the bottom one to the top left one in each by-pass triangle.}\label{fig_by_pass_3}
\end{figure}

\bdefn
We denote the two contact structures in Theorem \ref{thm_contact_structures_on_T_2_cross_I} by $\xi_{+,n}$ and $\xi_{-,n}$, respectively. The corresponding two by-passes are called {\it positive} and {\it negative}, respectively. The two by-passes can be distinguished by Figure \ref{fig_by_pass_3}.
\edefn

There are by-pass exact triangles associated to the positive and negative by-passes:

\begin{equation}\label{eq_by_pass_on_knot_complement}
\xymatrix{
\underline{\rm SHM}(-Y(K),-\Ga_{n+1})\ar[rr]^{\psi_{\pm,\infty}^{n+1}}&&\shm(-Y(K),-\Ga_{\infty})\ar[dl]^{\psi_{\pm,n}^{\infty}}
\\
&\shm(-Y(K),-\Ga_{n})\ar[lu]^{\psi_{\pm,n}^{n+1}}&
}
\end{equation}

Note we have $\psi_{\pm,n}^{n+1}=\Phi_{\xi_{\pm,n}}$. We also need the following facts.

\bprop[Honda \cite{honda2000classification}]\label{prop_composition_of_contact_structures}
On $T^2\times[0,2]$, the two contact structures $\xi_{-,n}\cup \xi_{+,n+1}$ and $\xi_{+,n}\cup\xi_{-,n+1}$ are the same.
\eprop

\bcor\label{cor_commutative_diagram}
We have a commutative diagram
\begin{equation*}
\xymatrix{
\shm(Y(K),\Ga_n)\ar[rr]^{\psi_{-,n+1}^n}\ar[dd]^{\psi_{+,n+1}^n}&&\shm(Y(K),\Ga_{n+1})\ar[dd]^{\psi_{+,n+2}^{n+1}}\\
&&\\
\shm(Y(K),\Ga_{n+1})\ar[rr]^{\psi_{-,n+2}^{n+1}}&&\shm(Y(K),\Ga_{n+2})
}
\end{equation*}
\ecor

\bpf
The corollary follows from proposition \ref{prop_composition_of_contact_structures} and theorem \ref{thm_gluing_map}.
\epf

There is a second way to interpret the maps $\psi_{\pm}$ associated to by-pass attachments by Ozbagci. In \cite{ozbagci2011contact}, he proved that a by-pass attachment can be realized by attaching a contact $1$-handle followed by a contact $2$-handle. In sutured monopole Floer theory, we have maps associated to the contact handle attachments, due to Baldwin and Sivek \cite{baldwin2016contact}. So, we can compose those contact handle attaching maps to define $\psi_{\pm}$. This is the original way that Baldwin and Sivek constructed the by-pass maps (when they defined by-pass maps, there was no construction of gluing maps) and proved the existence of the exact triangle. The two interpretations are the same because of the functoriality of the gluing maps. We will use this second point of view in the proof of Proposition \ref{prop_general_degree_shifting_for_knot_complement}.

\section{An Alexander grading}\label{sec_alexander_grading}
\subsection{The construction}
\bdefn\label{defn_stabilization_grading}
Suppose $(M,\ga)$ is a balanced sutured manifold, and $S$ is a properly embedded oriented surface. A {\it stabilization} of $S$ is an isotopy of $S$ to a surface $S'$, so that the isotopy creates a new pair of intersection points:
$$\partial S'\cap\ga=(\partial{S}\cap\ga)\cup \{p_+,p_-\}.$$
We require that there are arcs $\al\subset \partial{S'}$ and $\be\subset \ga$, oriented in the same way as $\partial{S'}$ and $\ga$, respectively, such that the following is true.

(1) We have $\partial{\al}=\partial{\be}=\{p_+,p_-\}$.

(2) The curves $\al$ and $\be$ cobound a disk $D$ so that ${\rm int}(D)\cap (\ga\cup \partial{S}')=\emptyset$.

\begin{figure}[h]
\centering
\begin{overpic}[width=4.0in]{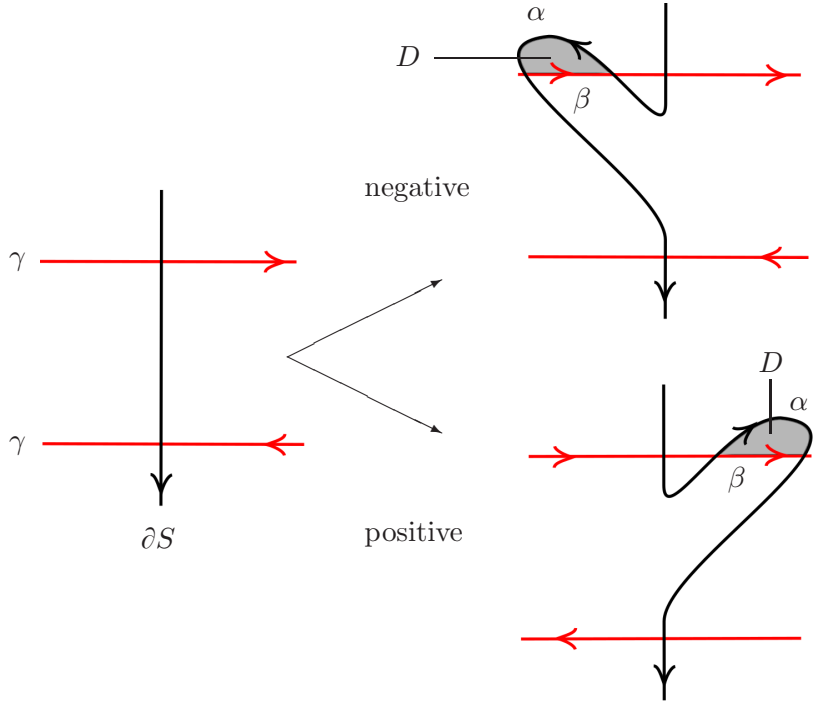}
	\put(13,20){$\partial{S}$}
	\put(-4,57){$\ga$}
	\put(-4,33){$\ga$}
	\put(32,45){\vector(2,1){20}}
	\put(32,45){\vector(2,-1){20}}
	\put(42,21){positive}
	\put(42,66){negative}
	\put(46,83){$D$}
	\put(51,84){\line(1,0){15}}
	\put(63,89){$\al$}
	\put(69,77.5){$\be$}
	\put(93,43){$D$}
	\put(94.5,42){\line(0,-1){7}}
	\put(97,38){$\al$}
	\put(89,28){$\be$}
\end{overpic}
\vspace{0.05in}
\caption{The positive and negative stabilizations of $S$.}\label{fig_pm_stabilization_of_surfaces}
\end{figure}

The stabilization is called {\it negative} if $D$ can be oriented so that $\partial{D}=\al\cup \be$ as oriented curves. it is called {\it positive} if $\partial{D}=(-\al)\cup\be$. See Figure \ref{fig_pm_stabilization_of_surfaces}.

We denote by $S^{\pm k}$ the result of performing $k$ many positive or negative stabilizations of $S$.
\edefn

The following lemma is straightforward.
\blem\label{lem_positive_negative_stabilization_and_decomposition}
Suppose $(M,\ga)$ is a balanced sutured manifold, and $S$ is a properly embedded oriented surface. Suppose $S^+$ and $S^-$ are the results of doing a positive and negative stabilization on $S$, respectively. Then, we have the following.

(1) If we decompose $(-M,-\ga)$ along $S$ or $S^-$, then the resulting two balanced sutured manifolds are diffeomorphic.

(2) If we decompose $(-M,-\ga)$ along $S^+$, then the resulting balanced sutured manifold $(M',\ga')$ is not taut, as $R_{\pm}(\ga')$ would both become compressible.
\elem

Suppose $(M,\ga)$ is a balanced sutured manifold, and $S$ is a properly embedded oriented surface. Suppose further that $S$ has precisely one boundary component and $\partial{S}$ intersects $\ga$ at $2n$ points. Since $\ga$ is parallel to the boundary of $R_+(\ga)$, it is null-homologous, so the algebraic intersection number of $\partial{S}$ with $\ga$ on $\partial{M}$ must be zero. We also assume that $n=2k+1$ is odd, as this can be achieved by a stabilization of $S$ if needed. Suppose the intersection points are $p_1,...,p_{2n}$, and they are indexed according to the orientation of $\partial S$.

Now pick a connected auxiliary surface $T$ for $(M,\ga)$, which is of large enough genus. Let $f:\partial{T}\ra\ga$ be an orientation reversing diffeomorphism and let $p_i'=f^{-1}(p_i)$. Suppose $\al_1,...,\al_n$ are pair-wise disjoint simple arcs on $T$, so that the following is true.

(1) The classes $[\al_1],...,[\al_n]$ are linearly independent in $H_1(T,\partial{T})$. 

(2) We have that $\partial{\al}_1=\{p_1',p_2'\}$, and, for all $1\leq i\leq k$, we have
$$\partial{\al}_{2i}=\{p'_{4i-1},p'_{4i+2}\},~{\rm and}~\partial{\al}_{2i+1}=\{p'_{4i},p'_{4i+1}\}.$$
Let
$$\widetilde{M}=M\mathop{\cup}_{id\times f}[-1,1]\times T,~{\rm and}~\widetilde{S}=S\mathop{\cup}_{id\times f}(\bigcup_{i=1}^n[-1,1]\times \al_i).$$
We know that
$$\partial{\widetilde{M}}=R_+\cup R_-,~{\rm and}~\partial{\widetilde{S}}\cap R_{\pm}=\bigcup_{i=1}^{k+1}C_{i,\pm}.$$
Here we require that for $i=1,...,k+1$, 
$$\al_{2i-1}\times\{\pm1\}\subset C_{i,\pm}.$$
Pick an orientation preserving diffeomorphism $h:R_+\ra R_-$ so that for $i=1,...,k+1$,
$$h(C_{i,+})=C_{i,-}.$$
Then, we can use $h$ and $\widetilde{M}$ to obtian a closure $(Y,R)$ of $(M,\ga)$. The boundary components of the surface $\widetilde{S}$ are glued with each other under $h$, so $\widetilde{S}$ becomes a closed surface $\bar{S} \subset Y$. From the construction, we know that
$$\chi(\bar{S})=\chi(S)-n.$$
We pick a non-separating simple closed curve $\eta\subset R$, so that $\eta$ is disjoint from $\bar{S}\cap R$ and represents a class which is linearly independent from the classes represented by the components of $\bar{S}\cap R$ in $H_1(R)$.

\bdefn\label{defn_grading}
We say that the surface $\bar{S}\subset Y$ is {\it associated to} the surface $S\subset M$. We can use $\bar{S}$ to define a grading on ${ SHM}(M,\ga)$ as follows.
$${ SHM}(M,\ga,S,i)=\bigoplus_{\substack{\mathfrak{s}\in\mathfrak{S}(Y|R)\\ c_1(\mathfrak{s})[\bar{S}]=2i}}\widecheck{HM}_{\bullet}(Y,\mathfrak{s};\Gamma_{\eta}).$$
We say that this grading is {\it associated to} the surface $S\subset M$. When using the language of marked closures, the closure $(Y,R)$ corresponds to a marked closure $\mathcal{D}=(Y,R,m,r,\eta)$, and we write the grading as
$${\rm SHM}(\mathcal{D},S,i).$$

The grading on ${\rm SHM}(\mathcal{D})$ also induces a grading on $\shm(M,\ga)$, as stated in Theorem \ref{thm_well_definedness_of_the_grading}. We also say it is {\it associated to} $S$ and write
$$\shm(M,\ga,S,i).$$ 
\edefn

\bthm\label{thm_well_definedness_of_the_grading}
When $\partial S$ is connected, the grading on $\shm(M,\ga)$ associated to $S$ is well-defined. That is, it is independent of all the choices made in the construction.
\ethm

\bpf
There are four types of choices we made in the construction of the grading:

I. The point $p_1$ on $\partial{S}\cap \ga$.

II. The choice of the arcs $\al_1,...,\al_n$ on $T$.

III. The choice of the gluing diffeomorphism $h$.

IV. The genus of the closure.

The proof of Theorem makes up the rest of the current section. In particular, the results are stated in corollary \ref{cor_independence_of_type_I}, corollary \ref{cor_independent_of_choice_of_type_II}, proposition \ref{prop_well_definedness_of_grading_34}, and lemma \ref{lem_independence_of_choice_of_type_V}.
\epf

In \cite{baldwin2018khovanov}, Baldwin and Sivek have already dealt with the choices of type II, III and IV. Among them, the idea for type IV can be adapted to the setting of the current paper verbatim, so we do not bother to write down the proof again.

\blem[Baldwin and Sivek \cite{baldwin2018khovanov}]\label{lem_independence_of_choice_of_type_V}
The definition of the grading on $\shm(M,\ga)$ associated to the surface $S\subset M$ is independent of choices of type IV.
\elem

To deal with the choices of type II, we have the following lemma.

\blem\label{lem_base_change_on_surface}
Suppose $T$ is a compact connected oriented surface-with-boundary and is of large enough genus. Suppose further that $\{\al_1,...,\al_n\}$ is a set of properly embedded simple arcs on $T$ so that the following is true.

(1) The arcs $\al_1,...,\al_n$ are pair-wise disjoint.

(2) The arcs represent linearly independent classes $[\al_1],...,[\al_n]$ in $H_1(T,\partial{T})$.

Suppose $\{\al_1',...,\al_n'\}$ is another set of properly embedded simple arcs so that the following is true.

(3) For $i=1,...,n$, we have $\partial{\al}_i=\partial{\al}_i'$.

(4) The set of arcs $\{\al_1',...,\al_n'\}$ also satisfies the above conditions (1) and (2).

Then, there is an orientation preserving diffeomorphism $h:T\ra T$ so that $h$ fixes the boundary of $T$, and, for $i=1,...,n$, we have
$$h(\al_i)=\al_i'.$$
\elem
\bpf
Suppose $N$ is a product neighborhood of 
$$\al_1\cup...\cup\al_n\subset T.$$
Let $\widetilde{T}=T\backslash{\rm int}(N).$ The boundary $\partial{\widetilde{T}}$ consists of the following:
$$\partial{\widetilde{T}}=(\partial{T}\cap\widetilde{T})\cup(\bigcup_{i=1}^n\al_{i,+}\cup\al_{i,-}).$$
Here, $\al_{i,\pm}$ are parallel copies of $\al_i$, being part of the boundary of the product neighborhood $N$. From condition (2), we know that $\widetilde{T}$ is connected. Also, by construction,
$$\chi(\widetilde{T})=\chi(T)+n.$$ 

Similarly, we can pick $N'$ to be a product neighborhood of 
$$\al_1'\cup...\cup\al_n'\subset T,$$
and take
$$\widetilde{T}'=T\backslash{\rm int}(N'),~{\rm and}~\partial{\widetilde{T}'}=(\partial{T}\cap\widetilde{T}')\cup(\bigcup_{i=1}^n\al'_{i,+}\cup\al'_{i,-}).$$

By condition (3), we can assume that $N\cap\partial{T}=N'\cap\partial{T}$, so there is an orientation preserving diffeomorphism
$$f:\partial\widetilde{T}\ra\partial\widetilde{T}'$$
so that
$$f|_{\partial{T}\cap\widetilde{T}}=id,~{\rm and}~f(\al_{i,\pm})=\al'_{i,\pm}$$
for all $i=1,...,n$. Since we have
$$\chi(\widetilde{T}')=\chi(T)+n=\chi(\widetilde{T}),$$
the diffeomorphism $f$ extends to a diffeomorphism
$$g:\widetilde{T}\ra \widetilde{T}'.$$
Thus, we can glue $\widetilde{T}$ and $\widetilde{T}'$ along $\al_{i,\pm}$ and $\al'_{i,\pm}$, and $g$ is glued to become a diffeomorphism 
$$h:T\ra T$$
that is the desired one.
\epf

As discussed in \cite{baldwin2018khovanov}, Lemma \ref{lem_base_change_on_surface} gives rise to the following corollary.

\bcor\label{cor_independent_of_choice_of_type_II}
The grading on $\shm(M,\ga)$ associated to the surface $S\subset M$ is independent of choices of type II.
\ecor

We deal with the choices of type III in Subsection \ref{subsec_choices_of_type_III} and the choices of type I in Subsection \ref{subsec_choice_of_type_I}.
 
\subsection{A reformulation of Canonical maps}\label{subsec_choices_of_type_III}
In this subsection, we give an alternative description of the canonical maps $\Phi_{\mathcal{D},\mathcal{D}'}$, which was originally constructed by Baldwin and Sivek in \cite{baldwin2015naturality} for two different marked closures of the same genus. For our convenience, we only study the a special case as described in the following paragraph.

Suppose $(M,\ga)$ is a balanced sutured manifold and $T$ is a connected auxiliary surface. Let
$$\widetilde{M}=M\cup [-1,1]\times T,~\partial{M}=R_+\cup R_-.$$
Suppose $h_1$ and $h_2$ are two different gluing diffeomorphisms, and there are corresponding marked closures $\mathcal{D}_1=(Y_1,R_+,r_1,m,\eta)$ and $\mathcal{D}_2=(Y_2,R_+,r_2,m,\eta)$, respectively. Here, we choose the same non-separating simple closed curve $\eta$ on $R_+$ to support local coefficients. 

Let $h= h_1^{-1}\circ h_2$, and $Y^h$ be the mapping torus of $h$, i.e., the manifold obtained from $R_+\times[-1,1]$ by identifying $R_{+}\times\{1\}$ with $R_{+}\times\{-1\}$ via $h$. Then, we can obtain $Y_2$ from $Y_1$ and $Y^h$ as follows. Cut $Y_1$ open along $R_+\times\{0\}$ and cut $Y^h$ along $R_+\times\{0\}$. We can re-glue them via the identity map on $R_+$ to get a connected manifold. This resulting manifold is precisely $Y_2$. As in Theorem \ref{thm_floer_excision}, there is a cobordism $W$ from $Y_1\sqcup Y^h$ to $Y_2$, and $W$ induces an isomorphism: 
$$HM(W):HM(Y_1\sqcup Y^h|R_+\cup R_+)\ra HM(Y_2|R_+).$$
Note, from Lemma \ref{lem_surface_bundle_over_S_1}, we know that
$$HM(Y^h|R_+)\cong \mathcal{R}.$$
Let $a$ be a generator of $HM(Y^h|R_+)$ and let $\iota$ be the map
$$\iota:HM(Y_1|R_+)\ra HM(Y_1|R_+)\otimes HM(Y^h|R_+)\cong HM(Y_1\sqcup Y^h|R_+\cup R_+)$$
defined by
$$\iota(x)=x\otimes a.$$
We have the following proposition.

\bprop\label{prop_new_definition_for_canonical_map}
The canonical map $\Phi_{\mathcal{D}_1,\mathcal{D}_2}$ can be re-interpreted as
$$\Phi_{\mathcal{D}_1,\mathcal{D}_2}\doteq HM(W)\circ\iota.$$
\eprop

Before proving the proposition, we first use it to prove the fact that the definition of the grading is independent of the choices of type III. Suppose $(M,\ga)$ is a balanced sutured manifold and $S\subset M$ is a properly embedded surface with precisely one boundary component, so that $\partial{S}$ intersects $\ga$ at $2n$ points for some odd $n=2k+1$. Suppose further that, in the construction of the grading induced by $S$, the choices of type I, II, IV are fixed. This means that there is a connected auxiliary surface $T$ for $(M,\ga)$ and n arcs $\al_1,...,\al_n$ so that the following holds

(1) We have
$$\partial(\al_1\cup...\cup\al_n)=\partial{S}\cap \ga.$$

(2) Let
$$\partial(M\cup [-1,1]\times T)=R_+\cup R_-,~{\rm and}~\widetilde{S}=S\bigcup_{i=1^n}([-1,1]\times \al_i),$$
then we have
$$\partial{\widetilde{S}}\cap R_{\pm}=C_{1,\pm},...,C_{k+1,\pm}.$$

Suppose there are two gluing diffeomorphisms $h_1$ and $h_2$ so that, for $i=1,2$
$$h_i(C_{1,+}\cup...\cup C_{k+1,+})=C_{1,-}\cup...\cup C_{k+1,-}.$$
Suppose further that there are marked closures $\mathcal{D}_1=(Y_1,R_+,m,r_1,\eta)$ and $\mathcal{D}_2=(Y_2,R_+,m,r_2,\eta)$ corresponding to $h_1$ and $h_2$, respectively. Here, we choose the same non-separating simple closed curve $\eta\subset R_+$ to construct local coefficients. We have the following proposition.
\bprop\label{prop_well_definedness_of_grading_34}
For any $i\in\intg$, we have
$$\Phi_{\mathcal{D}_1,\mathcal{D}_2}:{\rm SHM}(\mathcal{D}_1,S,i)\xra{\cong}{\rm SHM}(\mathcal{D}_2,S,i).$$
As a result, the definition of the grading on $\shm(M,\ga)$ is independent of the choices of type III.
\eprop
\bpf
Let $h= h^{-1}_1\circ h_2$, and form $Y^h$ as in Proposition \ref{prop_new_definition_for_canonical_map}. From Lemma \ref{lem_surface_bundle_over_S_1}, there is a unique spin${}^c$ structure $\mathfrak{s}_0$ so that
$$HM(Y^h|R_+)= \widecheck{HM}_{\bullet}(Y^h,\mathfrak{s}_0;\Ga_{\eta})\cong \mathcal{R}.$$

There are tori inside $Y^h$: The cylinders $C_{i,+}\times[-1,1]\subset R_+\times[-1,1]$ are glued via $h$ to become a union of tori $T$. Lemma \ref{lem_adjunction_inequality} tells us that
$$c_1(\mathfrak{s}_0)[T]=0.$$

Let $\bar{S}_1\subset Y_1$ and $\bar{S}_2\subset Y_2$ be the surfaces induced by $S\subset M$ as in the construction of the grading. We know that there is a $3$-dimensional cobordism from $S_1\sqcup T$ to $S_2$ inside the the cobordism $W$. The construction of this ($3$-dimensional) cobordism is similar to that of the Floer excisions. If $\mathfrak{s}$ is a spin${}^c$ structure on $W$, which contributes non-trivially to the cobordism map $HM(W)$, then $\mathfrak{s}$ must restrict to $\mathfrak{s}_0$ on $Y^h$. Hence, we know that
$$c_1({\mathfrak{s}})([\bar{S}_2])=c_1({\mathfrak{s}})([\bar{S}_1]+[T])=c_1({\mathfrak{s}})([\bar{S}_1])+c_1({\mathfrak{s}_0})([T])=c_1({\mathfrak{s}})([\bar{S}_1]).$$
Thus, $HM(W)$ preserves the grading and so does $\Phi^g_{\mathcal{D}_1,\mathcal{D}_2}$, by Proposition \ref{prop_new_definition_for_canonical_map}.
\epf

Now we proceed to prove proposition \ref{prop_new_definition_for_canonical_map}. There are a few preparations we need.

\blem\label{lem_functoriality_of_new_canonical_maps}
Under the settings of Proposition \ref{prop_new_definition_for_canonical_map}, suppose we have a third gluing diffeomorphism $h_3$, $h'=h_2^{-1}\circ h_3$, and $h''=h\circ h'=h^{-1}_1\circ h_3$. Construct $W'$, $W''$, $\iota'$, and $\iota''$ just in the same way as we construct $W$ and $\iota$. Then, we have an identity:
\begin{equation}\label{eq_functoriality_of_new_canonical_maps}
HM(W'')\circ \iota''\doteq HM(W')\circ \iota'\circ HM(W)\circ \iota.
\end{equation}
\elem

\bpf
Let $Y_{h'}$ and $Y_{h''}$ be the mapping tori of $h'$ and $h''$, respectively.
Since $h''=h\circ h'$, there is an excision cobordism from $Y_{h}\sqcup Y_{h''}$ to $Y_{h''}$ just as we construct $W$, $W'$, and $W''$. Call this cobordism $-W_e^{\vee}$, and let $W_e$ be the cobordism from $Y_{h''}$ to $Y_{h}\sqcup Y_{h'}$, obtained by putting $-W_e^{\vee}$ up side down and then reversing the orientation. By Theorem \ref{thm_floer_excision} and Lemma \ref{lem_surface_bundle_over_S_1}, it is straightforward to see that
$$HM(W\cup W'\cup W_e)\circ\iota_3\doteq HM(W')\circ \iota'\circ HM(W)\circ \iota.$$
Hence, to prove (\ref{eq_functoriality_of_new_canonical_maps}), it is enough to show that
\begin{equation}\label{eq_identifying_cobordism_maps}
HM(W\cup W'\cup W_e)\doteq HM(W'').	
\end{equation}

\begin{figure}[h]
\centering
\begin{overpic}[width=4.0in]{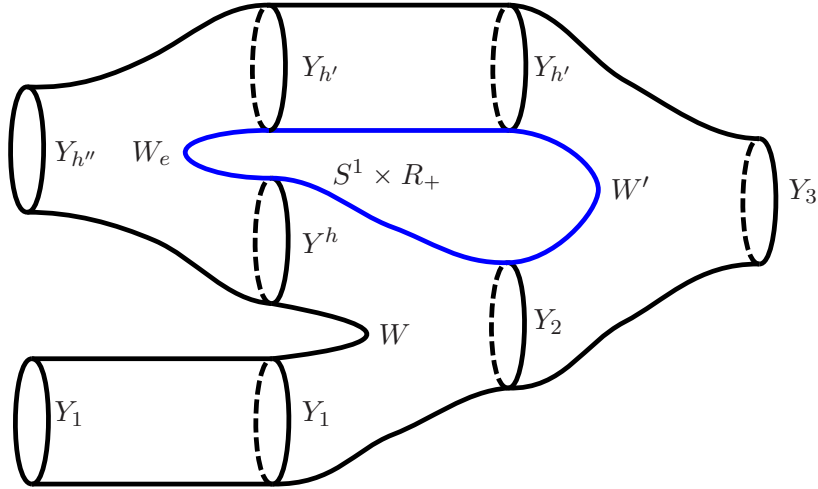}
	\put(42,40){$S^1\times R_+$}
	\put(48,19){$W$}
	\put(38,9){$Y_1$}
	\put(38,31){$Y^h$}
	\put(38,54){$Y_{h'}$}
	\put(68,21){$Y_2$}
	\put(68,54){$Y_{h'}$}
	\put(78,38){$W'$}
	\put(101,38){$Y_3$}
	\put(6,43){$Y_{h''}$}
	\put(16,43){$W_e$}
	\put(6,9){$Y_1$}
\end{overpic}
\vspace{0.05in}
\caption{The union $W\cup W'\cup W_e$. The (blue) curve in the middle represents the $3$-manifold $S^1\times R_+$ to cut along.}\label{fig_compo_of_exision_cob}
\end{figure}

However, we can cut $W'\cup W'\cup W_e$ open along the $3$-manifold $ S^1\times R_{+}$, as depicted in Figure \ref{fig_compo_of_exision_cob} and glue back two copies of $D^2\times R_+$. The resulting $4$-manifold is exactly $W''$. Hence, from Proposition 2.5 in \cite{kronheimer2010knots}, (\ref{eq_identifying_cobordism_maps}) holds true and we conclude the proof of lemma \ref{lem_functoriality_of_new_canonical_maps}.
\epf

\bcor\label{cor_identity_canonical_maps}
If $h_1=h_2$, then we have
$$HM(W)\circ\iota\doteq id.$$
\ecor

\bpf
From Theorem \ref{thm_floer_excision}, we know that
$$HM(W)\circ\iota$$
is an isomorphism. From Lemma \ref{lem_functoriality_of_new_canonical_maps}, we know that
$$HM(W)\circ\iota\circ HM(W)\circ\iota\doteq HM(W)\circ\iota.$$
Hence, the corollary follows.
\epf

\bpf[Proof of Proposition \ref{prop_new_definition_for_canonical_map}.]
Suppose $h$ is decomposed into Dehn twists:
$$h\sim D^{e_1}_{a_1}\circ...\circ D^{e_n}_{a_n},$$
as in Baldwin and Sivek \cite{baldwin2015naturality}. From Theorem \ref{thm_canonical_maps} and Lemma \ref{lem_functoriality_of_new_canonical_maps}, it is suffice to deal with the case when $n=1$, i.e., there is only one Dehn twist involved.

When $e_1=1$, the Dehn twist is positive. In this case, the canonical map $\Phi^g_{\mathcal{D}_1,\mathcal{D}_2}$ is constructed using the cobordism $W$, as in the hypothesis of Proposition \ref{prop_new_definition_for_canonical_map}, with the boundary component $Y^h$ capped off by the total space of a relative minimal Lefschetz fibration, see Lemma 4.9 in Baldwin and Sivek \cite{baldwin2015naturality}. Since such a Lefschetz fibration has relative monopole invariant being a unit in $\mathcal{R}$, as in Proposition B1 in \cite{baldwin2015naturality}, we conclude that
$$\Phi^g_{\mathcal{D}_1,\mathcal{D}_2}\doteq HM(W)\circ\iota.$$

When $e_1=-1$, the Dehn twist is negative. We can instead look at the canonical map $\Phi^g_{\mathcal{D}_2,\mathcal{D}_1}$. It corresponds to $h^{-1}$ and is constructed using a positive Dehn twist. Suppose we construct $W'$ and $\iota'$ out of $h^{-1}$, just as we construct $W$ and $\iota$ out of $h$. Then, from the previous case we know that
$$\Phi^g_{\mathcal{D}_2,\mathcal{D}_1}\doteq HM(W')\circ\iota'.$$
Then, the identity
$$\Phi^g_{\mathcal{D}_1,\mathcal{D}_2}\doteq HM(W)\circ\iota.$$
follows from Theorem \ref{thm_canonical_maps}, Lemma \ref{lem_functoriality_of_new_canonical_maps} and Corollary \ref{cor_identity_canonical_maps}.
\epf

\subsection{Pairing of the intersection points}\label{subsec_choice_of_type_I}
In this subsection, we deal with type I choices, i.e., the choice of $p_1$ among all intersection points in $S\cap\ga$. 

Let us first pick an arbitrary intersection point in $\partial{S}\cap\ga$ as $p_1$. We need to relax the requirement in the construction of the grading that $\partial{\al}_i$ are chosen to be a special pair of points in $S\cap \ga$. To record the data of the end points of $\al_{i}$, we make the following definition.

\bdefn\label{defn_pairing}
Suppose we have a collection of $n$ pair of numbers
$$\mathcal{P}=\{(i_1,j_1),...,(i_n,j_n)\}$$
so that 
$$\{i_1,j_1,...,i_n,j_n\}=\{1,2,...,2n\},$$
and, for all $l=1,...,n$, we have
$$i_l\notequiv j_l~({\rm mod}~2).$$
Then, we call such a collection $\mathcal{P}$ a {\it pairing} of {\it size} $n$. 
\edefn

Suppose $(M,\ga)$ is a balanced sutured manifold and $S\subset M$ is a properly embedded oriented surface. Suppose further that $S$ has a connected boundary, and it intersects $\ga$ at $2n=4k+2$ points. Those points are labeled by $p_1,...,p_{4k+2}$, according to the orientation of $\partial{S}$, with an arbitrary chosen starting point $p_1$. Continuing, suppose $\mathcal{P}=\{(i_l,j_l)\}_{l=1}^n$ is a pairing of size $n$, $T$ is an auxiliary surface of $M$, and $\al_1,...,\al_n$ are pair-wise disjoint simple arcs so that the following is true. 

(1) The arcs $\al_1$,..., $\al_n$ represent linearly independent classes in $H_1(T,\partial{T})$.

(2) For $l=1,...,n$, we have
$$\partial{\al_l}=\{p_{i_l},p_{j_l}\}.$$

Then, as in Definition \ref{defn_grading}, we can construct
$$\widetilde{M}=M\cup T\times[-1,1],~\widetilde{S}_{\mathcal{P}}=S\cup (\bigcup_{l=1}^n\al_l\times[-1,1]).$$
We have
$$\partial{\widetilde{M}}=R_+\cup R_-,\partial\widetilde{S}_{\mathcal{P}}\cap R_{\pm}=C_{1,\pm}\cup C_{s_{\pm},\pm}.$$
In general, the numbers of intersection circles, $s_{+}$ and $s_{-}$, are not necessarily equal to each other, so we make the following definition.

\bdefn\label{defn_balanced_pairing}
A pairing $\mathcal{P}$ is called balanced if $s_-=s_+$.
\edefn

\bexmp\label{exmp_pairings}
Here are some examples of the pairings. Assume $n=2k+1$ is odd.

(1) The simplest possible pairing
$$\mathcal{P}=\{(1,2),(3,4),...,(4k+1,4k+2)\}$$
has $s_-=1$ and $s_+=n$, or $s_-=n$ and $s_+=1$, depending on the choice of the starting point $p_1$, so it is not a balanced paring for $n>1$.

(2) In Definition \ref{defn_grading}, we have a paring arising from the construction of the grading:
$$\mathcal{P}^g=\{(1,2),(3,6),(4,5),...,(4k-1,4k+2),(4k,4k+1)\}.$$
This is an example of a balanced pairing, with $s_+=s_-=k+1$. 

(3) There is a very special balanced pairing with $s_+=s_-=1$:
$$\mathcal{P}^s=\{(1,2k+2),(2,2k+3),...,(2k+1,4k+2)\}.$$
\eexmp

If $(M,\ga)$, $S$, and $p_1$ are chosen as above, and we are equipped with a balanced pairing $\mathcal{P}$, then we can repeat the construction in Definition \ref{defn_grading} and define a grading on $\shm(M,\ga)$. By Corollary \ref{cor_independent_of_choice_of_type_II}, Proposition \ref{prop_well_definedness_of_grading_34}, and Lemma \ref{lem_independence_of_choice_of_type_V}, the grading depends only on the choice of $p_1$ and $\mathcal{P}$. Since $S$ and $p_1$ are fixed throughout this subsection, we omit them from the notation and write, in a moment, the grading as
$$\shm(M,\ga,\mathcal{P},i).$$

There is an operation we can perform on balanced pairings. Suppose $\mathcal{P}$ is a balanced pairing and we pick two indices $l_1$ and $l_2$ so that the following two conditions hold.

(i) The two arcs $\{1\}\times \al_{l_1}$ and $\{1\}\times \al_{l_2}$ are not contained in the same boundary component of ${\widetilde{S}_{\mathcal{P}}}.$ 

(ii) The two arcs $\{-1\}\times \al_{l_1}$ and $\{-1\}\times \al_{l_2}$ are not contained in the same boundary component of $\partial{\widetilde{S}}.$ 

Then, we can perform the following operation on $\mathcal{P}$: Suppose, in the two pairs $(i_{l_1},j_{l_1})$ and $(i_{l_2},j_{l_2})$, $i_{l_1}$ and $i_{l_2}$ are odd (and the two other numbers must be even), then we can obtain a new pairing $\mathcal{P'}$ out of $\mathcal{P}$ by removing the two pairs $(i_{l_1},j_{l_1})$ and $(i_{l_2},j_{l_2})$ from $\mathcal{P}$ and add two new pairings $(i_{l_1},j_{l_2})$ and $(i_{l_2},j_{l_1})$. 

\bdefn
We call the above operation the {\it cut and glue} on parings. Two pairings are called {\it equivalent} if one is obtained from the other by a cut and glue operation.
\edefn






\bexmp\label{exmp_cut_and_glue}
If $n=3$, $\mathcal{P}=\{(1,2),(3,6),(5,4)\}$, $l_1=1$, and $l_2=3$ ($l_1=1$ and $l_2=2$ do not meet the requirements of performing a cut and glue operation), then the resulting pairing $\mathcal{P}'$ is
$$\mathcal{P}'=\{(1,4),(3,6),(2,5)\},$$ 
and it is balanced.

It is obvious that the equivalence is an equivalent relation. Also, the result of a cut and glue operation on a balanced pairing is still a balanced one.
\eexmp

\begin{figure}
\centering
\begin{overpic}[width=5in]{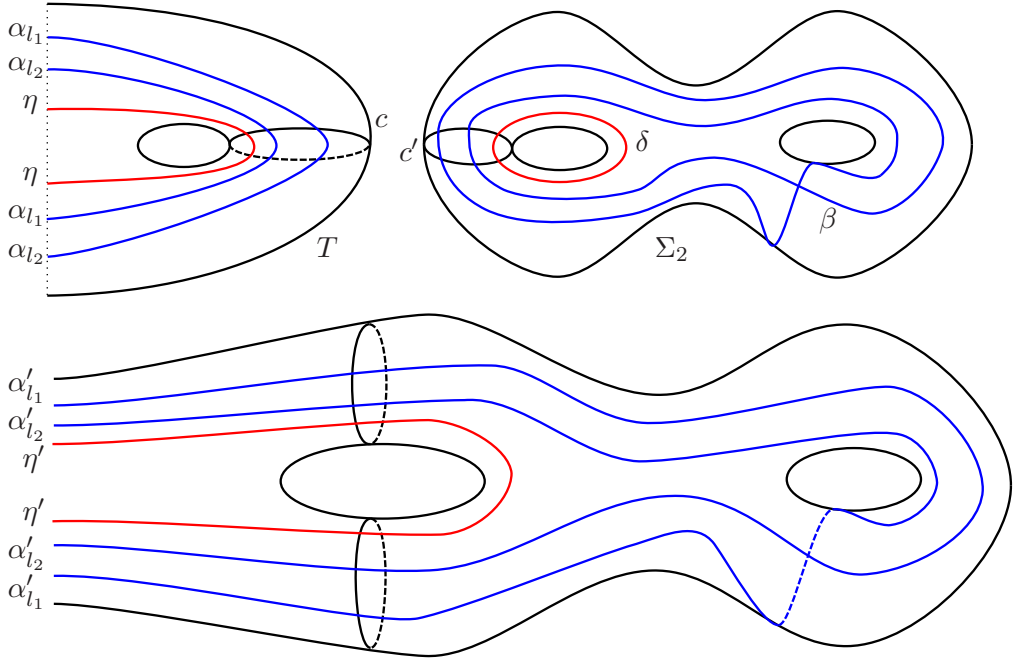}
	\put(28,42){$T$}
	\put(63,42){$\Sigma_2$}
	\put(-2.5,57.5){$\eta$}
	\put(-4,65){$\al_{l_1}$}
	\put(-4,61.5){$\al_{l_2}$}
	\put(-2.5,50){$\eta$}
	\put(-4,46){$\al_{l_1}$}
	\put(-4,42){$\al_{l_2}$}
	\put(34,55.5){$c$}
	\put(36.5,52){$c'$}
	\put(61,53){$\delta$}
	\put(80,45){$\be$}
	\put(-2.5,20){$\eta'$}
	\put(-4,28){$\al_{l_1}'$}
	\put(-4,24){$\al_{l_2}'$}
	\put(-2.5,14.5){$\eta'$}
	\put(-4,6.5){$\al_{l_1}'$}
	\put(-4,10.5){$\al_{l_2}'$}
\end{overpic}
\vspace{0.05in}
\caption{The auxiliary surface $T$ and the surface $\Sigma_2$}\label{fig_cut_and_glue_surfaces_along_tori}
\end{figure}

\blem\label{lem_equivalent_pairings_lead_to_same_grading}
Suppose a cut and glue operation on a balanced pairing $\mathcal{P}$ associated to the two indices $l_1$ and $l_2$ gives rise to a new balanced pairing $\mathcal{P}'$, then, for all $i\in\intg$, we have
$$\shm(M,\ga,\mathcal{P},i)=\shm(M,\ga,\mathcal{P}',i).$$
\elem
\bpf
At this point, we have shown that the choices of type II, III, and IV do not make difference on the definition of the grading. So, once $\mathcal{P}$ is chosen, we can freely choose other auxiliary data to construct the grading. Let $T$ and $\al_1,...,\al_n$ be chosen, and the pre-closure $\widetilde{M}$ as well as the properly embedded surface $\widetilde{S}_{\mathcal{P}}$ have been constructed. We can assume that they are chosen so that there is a curve $c$ intersecting both $\al_{l_1}$ and $\al_{l_2}$ transversely at one point. See Figure \ref{fig_cut_and_glue_surfaces_along_tori}. The requirements (i) and (ii) make sure that $\{\pm1\}\times \al_{l_1}$ and $\{\pm1\}\times \al_{l_2}$ lie in four different boundary components of $\widetilde{S}_{\mathcal{P}}$. So, there is an orientation preserving diffeomorphism $h:R_+\ra R_-$, where $\partial{\widetilde{M}}=R_+\cup R_-$, so that
$$h(\partial{\widetilde{S}}\cap R_+)=\partial{\widetilde{S}}\cap R_-, ~ h(c\times\{1\})=c\times\{-1\},$$
$$h(\al_{l_1}\times\{1\})=\al_{l_1}\times\{-1\},~{\rm and}~h(\al_{l_2}\times\{1\})=\al_{l_2}\times\{-1\}.$$

Let
$$Y=\widetilde{M}\mathop{\cup}_{id\cup h} [-1,1]\times R_+,~{\rm and}~R=\{0\}\times R$$
be a closure of $(M,\ga)$. The surface $\widetilde{S}_{\mathcal{P}}$ becomes a closed surface $\bar{S}_{\mathcal{P}}\subset Y$. We can also choose a simple closed curve $\eta$ on $R=\{0\}\times R_+$ so that $\eta$ is disjoint from $\widetilde{S}_{\mathcal{P}}$ and $\eta$ intersects $c\times\{0\}$ transversely at one point. Hence, we obtain a marked closure $\mathcal{D}=(Y,R,m,r,\eta)$, where $m$ and $r$ are both inclusions.

By definition, we have
$$SHM(\mathcal{D},\mathcal{P},i)=\bigoplus_{\substack{\mathfrak{s}\in\mathfrak{S}(Y|R)\\ c_1(\mathfrak{s})[\bar{S}_{\mathcal{P}}]=2i}}\widecheck{HM}_{\bullet}(Y,\mathfrak{s};\Gamma_{\eta}).$$

Let $\Sigma_2$ be a closed connected oriented surface of genus $2$. Let $c'$, $\delta$ and $\be$ be three simple closed curves on $\Sigma_2$, as depicted in Figure \ref{fig_cut_and_glue_surfaces_along_tori}.

Let $Y_{\Sigma}$ be the $3$-manifold $S^1\times \Sigma_2$. There is a torus $\Sigma=S^1\times c\subset Y$ and a torus $\Sigma'=S^1\times c'\subset S^1\times\Sigma_2$. We can choose an orientation preserving diffeomorphism $h':\Sigma\ra \Sigma'$ so that, for all $t\in S^1$, we have $h'(\{t\}\times c)=\{t\}\times c'$ as well as
$$h'(\{t\}\times((\al_{l_1}\cap c)\cup(\al_{l_2}\cap c)))=\{t\}\times(\be\cap c').$$

We can use $\Sigma$, $\Sigma'$, and $h'$ to perform a Floer excision on $Y\sqcup Y_{\Sigma}$. The result is a $3$-manifold $Y'$, with a distinguishing surface $R'$, obtained from $R\sqcup \Sigma_2$ by cutting and re-gluing along the two curves $c$ and $c'$. The surface $\bar{S}_{\mathcal{P}}\subset Y$ also becomes a new closed surface $\bar{S}_{\mathcal{P}'}\subset Y'$, obtained from $\bar{S}\sqcup(S^1\times\be)$ by cutting and re-gluing along four curves $S^1\times (\al_{l_1}\cap c)$, $S^1\times (\al_{l_2}\cap c)$, and $S^1\times(\be\cap c')$ (there are two intersection points of $\be$ with $c'$). The curve $\eta$ together with $\delta\subset\Sigma_2$ gives rise to a simple closed curve $\eta'\subset R'$. See Figure \ref{fig_cut_and_glue_surfaces_along_tori}. Hence, we get a new marked closure $\mathcal{D}'=(Y',R',m',r',\eta')$. The Floer excision results in a cobordism $W$ from $Y\sqcup Y_{\Sigma}$ to $Y'$ and a map
$$HM(W):HM(Y\sqcup Y_{\Sigma}|R\cup \Sigma_2;\Ga_{\eta\cup\delta})\ra HM(Y'|R';\Ga_{\eta'}).$$

Let $a\in HM(Y_{\Sigma}|\Sigma_2;\Ga_\delta)\cong\mathcal{R}$ be a generator. Then, we can define
$$\iota: HM(Y|R;\Ga_{\eta})\ra HM(Y'|R';\Ga_{\eta'})$$
as $\iota(x)=x\otimes a$ and we know that
$$\Phi_{\mathcal{D},\mathcal{D}'}=HM(W)\circ \iota,$$
by the definition of Canonical maps in Baldwin and Sivek \cite{baldwin2015naturality}.

The surface $\bar{S}_{\mathcal{P}'}\subset Y'$ can also be obtained from the balanced pairing $\mathcal{P}'$, which is obtained by performing a cut and glue operation on $\mathcal{P}$ associated to the two indices $l_1$ and $l_2$. Just as we did in the proof of Proposition \ref{prop_well_definedness_of_grading_34}, we conclude that, for all $i$,
$$\Phi_{\mathcal{D},\mathcal{D}'}(SHM(\mathcal{D},\mathcal{P},i))=SHM(\mathcal{D'},\mathcal{P'},i).$$
This concludes the proof of Lemma \ref{lem_equivalent_pairings_lead_to_same_grading}.
\epf

\bdefn
Two balanced pairings $\mathcal{P},\mathcal{P}'$ are called {\it connected} if there is a sequence of balanced pairings
$$\mathcal{P}_0=\mathcal{P},\mathcal{P}_1,...,\mathcal{P}_n=\mathcal{P}',$$
so that, for all $i=0,1,...,n-1,$ $\mathcal{P}_i$ and $\mathcal{P}_{i+1}$ are equivalent.
\edefn

\blem\label{lem_connect_two_special_pairings}
For any odd $n$, the two special balanced pairings $\mathcal{P}^g$ and $\mathcal{P}^s$ in Example \ref{exmp_pairings} are connected to each other. 
\elem

\bpf
In Example \ref{exmp_cut_and_glue}, we have shown that
$$\{(1,2),(3,6),(4,5)\}~{\rm and}~\{(1,4),(2,5),(3,6)\}$$
are equivalent. In a similar way, we can also show that
$$\{(1,6),(2,4),(3,5)\}~{\rm and}~\{(1,4),(2,5),(3,6)\}$$
are equivalent. So, 
$$\{(1,2),(3,6),(4,5)\}~{\rm and}~\{(1,6),(2,4),(3,5)\}$$
are connected. The later one can be thought of being obtained from the former one by sliding the arc $\al_1$, which originally joined the points $p_1$ and $p_2$, over the two arcs $\al_2$ and $\al_3$.

If we ignore the pairs $(2,4)$ and $(3,5)$ and look at
$\{(1,6),(7,10),(8,9)\}$, then the above argument applies again and we can connect it to
$\{(1,10),(6,9),(7,8)\},$ and this can be thought of further sliding $\al_1$ over $\al_4$ and $\al_5$. We can repeat this step for many times.

{\bf Case 1.} If $n$ is of the form $4k+1$. In this case, we can slide $\al_1$ over to join $p_1$ with $p_{4k+2}$. Hence, $\mathcal{P}^g$ is connected to a new balanced pairing 
\beq
\mathcal{P}'=&\{(1,n+1=4k+2),(2,5),(3,4),...,(4k-2,4k+1),(4k-1,4k),\\
&(4k+3,4k+6),(4k+4,4k+5),...,(8k-1,8k+2),(8k,8k+1)\}.
\eeq
Then, we can perform cut and glue operations on pairs $(4l-2,4l+1)$ and $(4l-2+n,4l+1+n)$ as well as on pairs $(4l-1,4l)$ and $(4l-1+n, 4l+n)$, for all $1\leq l\leq k$. The result of these operations is nothing but the special balanced paring $\mathcal{P}^s$ introduced in Example \ref{exmp_pairings}. Hence, we are done.

{\bf Case 2.} If $n$ is of the form $4k+3$. In this case, we can slide $\al_1$ to join $p_1$ with $p_{4k+2}$, so the balanced pairing $\mathcal{P}^g$ is connected to 
\beq
\mathcal{P}'=&\{(1,4k+2),(2,5),(3,4),...,(4k-2,4k+1),(4k-1,4k),\\
&(4k+3,4k+6),(4k+4,4k+5),...,(8k+3,8k+6),(8k+4,8k+5)\}.
\eeq

Perform another cut and glue operation on pairs $(1,4k+2)$ and $(4k+4,4k+5)$, then we get a new balanced pairing
\beq
\mathcal{P}'=&\{(1,n+1=4k+4),(2,5),(3,4),...,(4k-2,4k+1),(4k-1,4k),\\
&(4k+2,4k+5),(4k+3,4k+6),...,(8k+3,8k+6),(8k+4,8k+5)\}.
\eeq
There is, then, an arc joining $p_{4k+2}$ and $p_{4k+5}$, and we can slide it over to join $p_{4k+5}$ and $p_2$. Similarly, there is an arc joining $p_{4k+3}$ with $p_{4k+6}$, and we can slide it over to join $p_{4k+3}$ with $p_{8k+6}$. Then, $\mathcal{P}^g$ is connected to a new balanced pairing
\beq
\mathcal{P}''=&\{(1,n+1=4k+4),(2,n+2=4k+5),(n=4k+3,2n=8k+6),\\
&(3,6),(4,5)...(4k-1,4k+2),(4k,4k+1)\\
&(4k+6,4k+9),(4k+7,4k+8),...,(8k+2,8k+5),(8k+3,8k+4)\}.
\eeq

Finally, we can perform cut and glue operations on pairs $(4l-1,4l+2)$ and $(4l-1+n,4l+2+n)$ as well as on $(4l,4l+1)$ and $(4l+n,4l+1+n)$, for all $1\leq l\leq k$, then the final result is $\mathcal{P}^s$, and we conclude the proof of Lemma \ref{lem_connect_two_special_pairings}.
\epf

\bcor\label{cor_independence_of_type_I}
The definition of the grading on $\shm(M,\ga)$ is independent of choices of type I.
\ecor
\bpf
It is straightforward to check that if we use the special balanced pairing $\mathcal{P}^s$, then the surface $\widetilde{S}_{\mathcal{P}^s}$ is the same for all possible choices of the starting point $p_1$. Hence the corollary follows from Lemma \ref{lem_equivalent_pairings_lead_to_same_grading} and Lemma \ref{lem_connect_two_special_pairings}.
\epf

\brem
We want to use $\mathcal{P}^{g}$ in the definition of grading because it is more convenient to use this construction to discuss about the positive and negative stabilizations (see Definition \ref{defn_stabilization_grading}), as we will see in Subsection \ref{sec_grading_shift}.
\erem

Though we only discussed some special pairings, we would like to make the following conjecture. Note the concept of balancedness, equivalence, connectedness defined above can be reached in a purely combinatorial way and is independent of all the topological input.

\begin{conj}
Any two balanced pairings of the same size $n$, where $n$ is odd, are connected.	
\end{conj}

\section{The grading shifting property}\label{sec_grading_shift}
\subsection{A naive version}
Suppose $(M,\ga)$ is a balanced sutured manifold and suppose $S$ is a properly embedded surface in $M$ with a connected boundary. In Definition \ref{defn_grading}, we constructed a grading on $\shm(M,\ga)$ associated to $S$, when $|\partial S\cap \ga|=2n$ with $n$ being odd. If $n$ is even, then we introduce, in Definition \ref{defn_stabilization_grading}, positive and negative stabilizations $S^{\pm}$ that both increase $n$ by $1$. It is a natural question to ask how the gradings associated to $S^+$ and $S^-$ are related to each other. The following proposition is a first answer to this question.

\bprop\label{prop_naive_version_of_degree_shifting}
Suppose $(M,\ga)$ is a balanced sutured manifold, $S\subset M$ is a properly embedded surface with a connected boundary, and that $\partial{S}$ intersects $\ga$ transversely at $2n$ points with $n=2k>0$ odd. Suppose further that the balanced sutured manifold obtained by decomposing $(-M,-\ga)$ along $S$ is taut. Let $S^+$ and $S^-$ are the positive and negative stabilizations of $S$, respectively. Suppose $S$ is of genus $g$ and let
$$g_c=g+k.$$
Then, we have
$$\shm(-M,-\ga,S^-,g_c)\subset \shm(-M,-\ga,S^+,g_c-1).$$
\eprop

We need the following lemma before proving Proposition \ref{prop_naive_version_of_degree_shifting}.
\blem[Kronheimer and Mrowka \cite{kronheimer2010knots}]\label{lem_decomposition_give_top_grading}
Suppose $(M,\ga)$ is a balanced sutured manifold and $S$ is properly embedded surface inside $M$ so that $\partial{S}$ is connected and $|\partial{S}\cap \ga|=2n$ with $n$ odd. Let
$$g_c=\frac{n-1}{2}+g(S),$$
then we know that
$$\shm(M,\ga,S,i)=0$$
for all $i>g_c$, and
$$\shm(M,\ga,S,g_c)\cong\shm(M',\ga'),$$
where $(M',\ga')$ is the balanced sutured manifold obtained from $(M,\ga)$ by decomposing along $S$.
\elem

\bpf
This is a reformulation of Proposition 6.9 in Kronheimer and Mrowka \cite{kronheimer2010knots}, using our definition of the gradings in Definition \ref{defn_grading}. The fact that $\shm(M,\ga,S,i)=0$ for all $i>g_c$ follows directly from the adjunction inequality in Lemma \ref{lem_adjunction_inequality}.
\epf

\bpf[Proof of proposition \ref{prop_naive_version_of_degree_shifting}.]
If we have two different negative stabilizations $S^-_1$ and $S^-_2$, then we know from Lemma \ref{lem_positive_negative_stabilization_and_decomposition} and Lemma \ref{lem_decomposition_give_top_grading} that
$$\shm(-M,-\ga,S^-_1,g_c)\cong\shm(-M',-\ga')\cong\shm(-M,-\ga,S^-_2,g_c),$$
where $(M',\ga')$ is obtained from $(-M,-\ga)$ by performing a sutured manifold decomposition along $S$. Hence, we can choose a special negative stabilization to deal with.

\begin{figure}[h]
\centering
\begin{overpic}[width=5.0in]{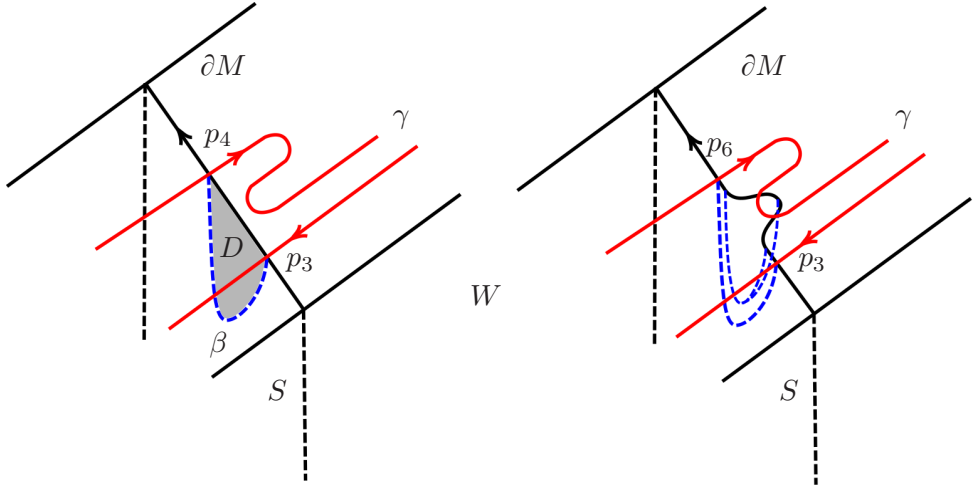}
	\put(40,38){$\ga$}
	\put(48,19){$W$}
	\put(80,9){$S$}
	\put(76,43){$\partial{M}$}
	\put(92,38){$\ga$}
	\put(20,43){$\partial{M}$}
	\put(27,9){$S$}
	\put(29,23){$p_3$}
	\put(20.5,36){$p_4$}
	\put(22,24){$D$}
	\put(21,14){$\be$}
	\put(82,23){$p_3$}
	\put(72.5,35){$p_6$}
\end{overpic}
\vspace{0.05in}
\caption{A negative stabilization of $S$. Positive stabilizations are similar.}\label{fig_negative_stabilization_surface}
\end{figure}

Suppose the intersection points of $\partial{S}\cap \ga$ are labeled as $p_1,...,p_{2n}$ according to the orientation of $\partial S$. When labeling the points, we need to pick a suitable $p_1$ so that the new pair of intersection points created by the positive or negative stabilization lie between $p_3$ and $p_4$. Let $\be'\subset \partial{S}$ be part of $\partial{S}$ so that $\partial{\be'}=\{p_3,p_4\}$ and $\be'$ contains no other intersection points $p_j$ for $j\neq 3,4$. Let $\be\subset S$ be a properly embedded arc so that $\partial{\be}=\{p_3,p_4\}$, $\be$ and $\be'$ co-bound a disk on $D$, and when performing positive and negative stabilizations, the isotopy on $S$ can be fixed outside the disk $D$. Now if we use the same starting point $p_1$ to label $\partial{S}^{\pm}\cap \ga$, then the new pair of intersection points are both $p_4$ and $p_5$ in the two cases. See Figure \ref{fig_negative_stabilization_surface}.

Suppose $T$ is an auxiliary surface for $(M,\ga)$ of large enough genus. When constructing the grading associated to $S^{\pm}$, we need to choose linearly independent arcs $\al_{1},\al_2,\al_3^{\pm},\al_4...,\al_{n+1}\subset T$ and the special pairing $\mathcal{P}^g$, which is defined in Example \ref{exmp_pairings}, to make it clear what are the end points of the arcs $\al_i$. Here, $\al_{3}^{\pm}$ correspond to the different surfaces $S^{\pm}$, while $T$ and all other arcs $\al_i$, for $i\neq3$, can be chosen to be the same for both $S^+$ and $S^-$. In the pre-closure $\widetilde{M}=M\cup [-1,1]\times T$, we have two surfaces $\widetilde{S}^{+}$ and $\widetilde{S}^{-}$. After picking suitable gluing diffeomorphisms $h^{\pm}$, we get two marked closures
$$\mathcal{D}^{+}=(Y^{+},R^{+},r^{+},m^{+},\eta^{+})~{\rm and}~\mathcal{D}^{-}=(Y^{-},R^{-},r^{-},m^{-},\eta^{-})$$
so that there are closed surfaces $\bar{S}^{+}$ and $\bar{S}^{+}$ inside $Y^{+}$ and $Y^-$, respectively, and the gradings associated to $S^+$ and $S^-$ are defined by looking at the pairings between the first Chern classes of the spin${}^c$ structures on $Y^+$ and $Y^-$ with the fundamental classes of $\bar{S}^{+}$ and $\bar{S}^{-}$, respectively. Note the genus of $\bar{S}^{+}$ and $\bar{S}^-$ are both $g_c+1=g+k+1$.

From Proposition \ref{prop_new_definition_for_canonical_map}, we know that the canonical map $\Phi_{-\mathcal{D}^-,-\mathcal{D}^+}$ can be interpreted in terms of a Floer excision cobordism $W$ from $-Y^-\sqcup-Y^h$, where $Y^h$ is the mapping torus of $h=(h^-)^{-1}\circ h^+$, to $-Y^+$. 

We can construct a special closed surface of genus $2$ as follows. Recall we have an arc $\be\subset S$, and since the isotopies for positive or negative stabilizations are supported in the interior of the disk $D$, $\be$ also lies in $\bar{S}^{\pm}$. Let $\delta=\be\cup(\al_2\times\{0\})\subset \bar{S}^{\pm}$ be a closed curve. Then, the curve $\delta$ cuts each of $\bar{S}^{\pm}$ into two parts. One part contains $S\backslash {\rm int}(D)$ and the other part is a connected oriented surface $T^{\pm}\subset \bar{S}^{\pm}$ of genus 1 and with boundary $\delta$. Inside $W$, we can define
$$\Sigma_2=T^-\cup[0,1]\times\delta\cup -T^+\subset W.$$
It is straightforward to see that, in $H_2(W)$,
$$[\bar{S}^-]=[\bar{S}^+]+[\Sigma_2].$$
Hence, by the adjunction inequality in dimension $4$, which is a $4$-dimensional analogue of Lemma \ref{lem_adjunction_inequality}, we have

\beq
\Phi_{-\mathcal{D}^-,-\mathcal{D}^+}(SHM(-\mathcal{D},S^-,g_c))\subset&~~~SHM(-\mathcal{D},S^+,g_c+1)\\
&\oplus SHM(-\mathcal{D},S^+,g_c)\\
&\oplus SHM(-\mathcal{D},S^+,g_c-1).
\eeq
The adjunction inequality also implies that $SHM(-\mathcal{D},S^+,g_c+1)=0$. If we decompose $(-M,-\ga)$ along $S^+$, and suppose $(M',\ga')$ is the resulting balanced sutured manifold, then, by Lemma \ref{lem_positive_negative_stabilization_and_decomposition}, $R_{\pm}(\ga')$ is compressible and so
$$SHM(-\mathcal{D},S^+,g_c)\cong SHM(-M',-\ga')=0.$$
The first isomorphism follows from Lemma \ref{lem_decomposition_give_top_grading} and the second equality follows again from the adjunction inequality in Lemma \ref{lem_adjunction_inequality}.

Hence, the only possibility left is
$$\Phi_{-\mathcal{D}^-,-\mathcal{D}^+}(SHM(-\mathcal{D},S^-,g_c))\subset SHM(-\mathcal{D},S^+,g_c-1)$$
and we we conclude the proof of Proposition \ref{prop_naive_version_of_degree_shifting}.
\epf

\subsection{Knot complements with two-component sutures}
In this section, we focus on the case when the balanced sutured manifold $(M,\ga)$ is the complement of a null-homologous knot, i.e., $M=X(K)=X\backslash{\rm im}(N(K))$, where $X$ is a closed connected oriented $3$-manifold and $K\subset X$ is a null-homologous knot. Also, we assume that $\ga$ has two components. Under these conditions, we show that the result of Proposition \ref{prop_naive_version_of_degree_shifting} holds for not only the top grading but also all gradings.

\bprop\label{prop_degree_shifting_formula_for_toroidal_boundary_with_two_sutures}
Suppose $(M=X(K),\ga)$ is the balanced sutured manifold as described in the above paragraph. Suppose further that $S$ is a Seifert surface of the knot $K$, viewed as a properly embedded surface in $M$, so that $|\partial{S}\cap \ga|=2n$. Then, for any $p,k,l\in \intg$ such that $n+p$ is odd, we have
$$\shm(-M,-\ga,S^{p},l)=\shm(-M,-\ga,S^{p+2k},l-k).$$
Note $S^p$ is defined as in Definition \ref{defn_stabilization_grading}, and, in particular, $S^0=S$.
\eprop

Before proving Proposition \ref{prop_degree_shifting_formula_for_toroidal_boundary_with_two_sutures}, we will first deal with the following related proposition.

\bprop\label{prop_difference_of_spin_c_structures}
Suppose $(Y,R)$ is a closure of $(-M,-\ga)$, and let $\mathfrak{s}_1,\mathfrak{s}_2\in\mathfrak{S}^*(Y|R)$ (see Definition \ref{defn_set_of_spin_c_structures}) be two supporting spin${}^c$ structures on $Y$. Then, there is a $1$-cycle $x$ inside $M$, so that
$$P.D.c_1(\mathfrak{s}_1)-P.D.c_1(\mathfrak{s}_2)=[x]\in H_1(Y).$$
Note the cycle is contained in $M$ but the identity is on the whole $Y$.
\eprop

We now describe the closures of $(-M,-\ga)$. Note if $(Y,R)$ is a closure of $(M,\ga)$, then $(-Y,-R)$ is a closure of $(-M,-\ga)$. So, in the following discussion, we only describe the closures of $(M,\ga)$, and, to get closures of $(-M,-\ga)$, one simply reverses the orientations.

\begin{figure}[h]
\centering
\begin{overpic}[width=4.5in]{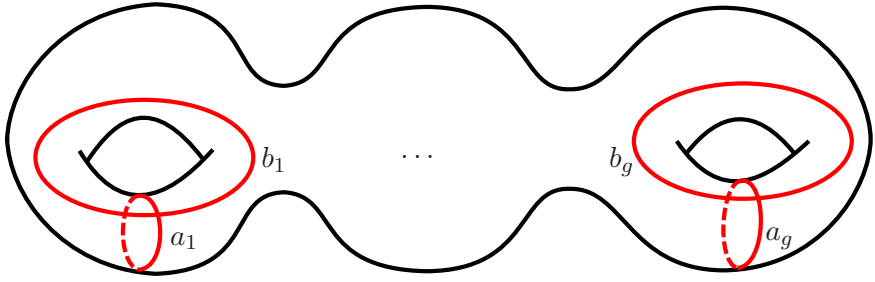}
	\put(19,4){$a_1$}
	\put(87.5,4.5){$a_g$}
	\put(29.5,13){$b_1$}
	\put(69.5,13){$b_g$}
	\put(45.5,14){$\dots$}
\end{overpic}
\vspace{0.05in}
\caption{The surface $\Sigma_g$.}\label{fig_Sigma_g}
\end{figure}

Let $\Sigma_g$ be a closed oriented connected surface of large enough genus $g$. Its first homology is generated by the classes $[a_1],[b_1],...,[a_g],[b_g]$, as depicted in Figure \ref{fig_Sigma_g}.

Let $T=\Sigma_g\backslash{\rm int}(N(a_1))$ be a surface obtained from $\Sigma_g$ by cutting $\Sigma_g$ open along $a_1$, then $T$ can be viewed as an auxiliary surface for $(M,\ga)$. Let
$$\widetilde{M}=M\cup [-1,1]\times T$$
be a pre-closure of $(M,\ga)$, and let
$$\partial{\widetilde{M}}=R_{+}\cup R_{-}.$$

If we choose a special gluing diffeomorphism $h^0:R_{+}\ra R_{-}$ so that $h_{T\times\{1\}}=id$, then we get a special marked closure 
$$\mathcal{D}^0=(Y^0,R, r^0,m^0,\eta).$$
Similar to the closures described in Section 5.1 in \cite{kronheimer2010knots}, the closure $(Y^0,R)$ can be achieved as follows: Let $\Sigma_g$ be the surface as in Figure \ref{fig_Sigma_g}, and let $Y_{\Sigma}=S^1\times\Sigma_g$. By abusing the notations, use $a_1$ to also denote the curve $\{1\}\times \al_1\subset Y_{\Sigma}$. Let $N(a_1)$ be a tubular neighborhood of $a_1\subset Y_{\Sigma}$. Note $a_1\subset \{1\}\times\Sigma_g$, so there is a framing on $\partial{N(a_1)}$ induced by $\{1\}\times\Sigma_g$. Let $\lambda_a$, $\mu_a$ be the longitude and meridian, respectively.

Then, we have
$$Y^0=M\mathop{\cup}_{\phi}(Y_{\Sigma}\backslash {\rm int}(N(a_1))).$$
Here,
$$\phi: \partial{N(a_1)}\ra \partial{M}$$
sends the two copies of $\lambda_a$ to the suture $\ga$. Note there are canonical ways to identify $R_{\pm}$ with $\Sigma_g$. So, in the marked closure $\mathcal{D}_0$, we have $R=\Sigma_g$.

\blem\label{lem_reference_closure_satisfies_proposition}
Proposition \ref{prop_difference_of_spin_c_structures} is true for $-Y^0$.
\elem

\bpf
From the Mayer-Vietoris sequece, we know that there is an exact sequence
$$H_1(T^2)\ra H_1(M)\oplus H_1(Y_{\Sigma}\backslash {\rm int}(N(a_1)))\ra H_1(Y^0)\ra0,$$
where $T^2=\partial{M}=\partial{(Y_{\Sigma}\backslash {\rm int}(N(a_1)))}.$
Hence, we conclude that
$$H_1(Y^0)=H_1(M)\oplus H_1(Y_{\Sigma}\backslash {\rm int}(N(a_1)))\slash\sim,$$
where $\sim$ is the relation induced by the gluing map $\phi:$
$$[\lambda_{a}]\sim \phi_*([\lambda_{a}]),~[\mu_{a}]\sim\phi_{*}([\mu_a]).$$
A direct calculation shows that
$$H_1(Y_{\Sigma}\backslash {\rm int}(N(a_1)))=\lgl [\mu_a],[a_1],[b_1],...,[a_g],[b_g],[s^0]\rgl,$$
where $s^0$ corresponds to the $S^1$ direction in $Y_{\Sigma}=\Sigma_g\times S^1$. 
Hence, we can write
\begin{equation}\label{eq_first_homology_of_the_reference_closure}H_1(Y^0)=H_1(M)\oplus\lgl [b_1],[a_2],[b_2],...,[a_g],[b_g],[s^0]\rgl.
\end{equation}
This is because $a_1$ and $\mu_a$ are absorbed into $H_1(M)$.

Suppose $\mathfrak{s}\in\mathfrak{S}^*(-Y^0|-\Sigma_g)$, then we can express $P.D.c_1(\mathfrak{s})$ in terms of the above basis. The coefficient of $[s]$ can be fixed by the evaluation
$$c_1(\mathfrak{s})[-\Sigma_g]=2g-2.$$
There are no $[b_1],[a_2],[b_2]...[a_g],[b_g]$ terms, since we can apply the adjunction inequality in Lemma \ref{lem_adjunction_inequality} to tori $a_1\times S^1,b_2\times S^1..., a_g\times S^1\subset Y^0$ to rule out those classes. The rest of the terms must then lie in $H_1(M)$. So, if we look at the difference (of the Poincar\'e dual of their first Chern class) of two supporting spin${}^c$ structures, it must lie in $M$.
\epf

Now we deal with general closures of $(-M,-\ga)$. As above, we have the pre-closure
$$\widetilde{M}=M\cup T\times[-1,1],$$
where $T=\Sigma_{g}\backslash N(a_1).$ Also, recall
$$\partial{\widetilde{M}}=R_+\cup R_-.$$
Note, as in the above discussion, there are canonical ways to identify $R_+$ and $R_-$ with $\Sigma_g$. We can pick any orientation preserving diffeomorphism $h:R_+\ra R_-$ to get a closure $(Y,\Sigma_g)$ of $(M,\ga)$, or a marked closure
$$\mathcal{D}=(Y,\Sigma_g,r,m,\eta).$$
In particular, the special marked closure $\mathcal{D}^0$ in Lemma \ref{lem_reference_closure_satisfies_proposition} corresponds to taking $h=h^0=id.$

Let $Y^h$ be the mapping torus of the diffeomorphism $h:\Sigma_g\ra \Sigma_g$, then we can reinterpret $Y$ as
$$Y=M\cup_{\phi}(Y^h\backslash {\rm int}(N(a_1))).$$
From Proposition \ref{prop_new_definition_for_canonical_map}, we know that the canonical map $\Phi_{\mathcal{D}_0,\mathcal{D}}$ can be obtained from a cobordism $W$ from $Y^0\sqcup Y^h$ to $Y$. The cobordism $W$ arises from the Floer excision as in Subsection \ref{subsec_naturality}. The computation of the first homologies of $Y$, $Y^h$ and $W_1$ are straightforward and we can describe them as follows
\begin{equation}\label{eq_first_homology_of_a_general_closure}
H_1(Y)=H_1(M)\oplus\lgl [\mu_a],[a_1],[b_1],...,[a_g],[b_g],[s]\rgl\slash\sim_{\phi,h}
\end{equation}
\begin{equation}\label{eq_first_homology_of_the_mapping_torus}
H_1(Y^h)=\lgl [a_1],[b_1]...,[a_g],[b_g],[s^{h}]\rgl\slash\sim_h
\end{equation}
\begin{equation}\label{eq_first_hmology_of_the_cobordism}
H_1(W)=H_1(M)\oplus\lgl [\mu_a],[a_1],[b_1],...,[a_g],[b_g],[s^0],[s^{h}]\rgl\slash\sim_{\phi,h}.
\end{equation}
Here, $s$ is a circle inside $Y$ which intersects $\Sigma_g$ once. We can isotope $h$ so that $h$ has a fixed point $p\in\Sigma_g$, then, inside $Y$, there is a circle $s=\{p\}\times S^1$. The class $s^h$ is similar.
The relations $\sim_{\phi,h}$ are
$$[a_1]\sim\phi_*([a_1]),~[\mu_a]\sim\phi_*([\mu_a]),~[a_i]\sim h([a_i]),~[b_i]\sim h([b_i]).$$
The relations $\sim_h$ are
$$[a_i]\sim h([a_i]),~[b_i]\sim h([b_i]).$$

\blem\label{lem_injective_inclusion_for_first_homology}
The inclusion $i:Y\hookrightarrow W$ induces an injective map
$$i_*:H_1(Y)\hookrightarrow H_1(W).$$
\elem

The following proposition is a built-in property of monopole Floer homology.
\blem\label{lem_spin_c_structure_shall_extend}
Suppose $(W,\nu)$ is an oriented cobordism between two pairs $(Y,\eta)$ and $(Y',\eta')$. Suppose further that $\mathfrak{s}$ is a spin${}^c$ structure on $Y$ and $\mathfrak{s}'$ is a spin${}^c$ structure on $Y'$ so that
$$\widecheck{HM}(W,\nu)(\widecheck{HM}_{\bullet}(Y,\mathfrak{s};\Gamma_{\eta}))\cap \widecheck{HM}_{\bullet}(Y',\mathfrak{s}';\Gamma_{\eta'})\neq \{0\},$$
then, we know that
$$i_*(P.D.c_1(\mathfrak{s}))=i_*'(P.D.c_1(\mathfrak{s'}))\in H_1(W).$$
Here, $i:Y\ra W$ and $i':Y'\ra W'$ are the inclusions.
\elem

Recall we defined $\mathfrak{S}^*(-Y^0|-\Sigma_g)$ to be the set of supporting spin${}^c$ structures as in Definition \ref{defn_set_of_spin_c_structures}. We can also define
$$\mathfrak{PDS}^*(-Y^0|-\Sigma_g)=\{P.D.c_1(\mathfrak{s})|\mathfrak{s}\in\mathfrak{S}^*(-Y^0|-\Sigma_g)\}.$$
We can define $\mathfrak{PDS}^*(-Y|-\Sigma_g)$ similarly. Then, we have the following lemma.
\blem\label{lem_the_map_rho}
Suppose we have the closures 
$$-\mathcal{D}_0=(-Y^0,-\Sigma_g,r,m,-\eta),~-\mathcal{D}=(-Y,-\Sigma_g,r,m,-\eta)$$ for $(-M,-\ga)$, the mapping torus $Y^h$ and the cobordism $W$ from $Y^0\sqcup Y^h$ to $Y$ defined as above. Suppose $\mathfrak{s}_h$ is the unique supporting spin${}^c$ structure on $-Y^h$ satisfying the statement of Lemma 2.6. Then, there exists a map
$$\rho:\mathfrak{PDS}^*(-Y^0|-\Sigma_g)\ra H_1(Y)$$ so that $\mathfrak{PDS}^*(-Y|-\Sigma_g)\subset {\rm im}(\rho)$ and $\rho$
satisfies the following property $(*)$: Suppose we have spin${}^c$ structures $\mathfrak{s}\in \mathfrak{S}^*(-Y^0|-\Sigma_g)$ and $\mathfrak{s}'\in \mathfrak{S}^*(-Y|-\Sigma_g)$, so that
$$\widecheck{HM}(-W)(\widecheck{HM}_{\bullet}(-Y^0,\mathfrak{s};\Gamma_{-\eta})\otimes \widecheck{HM}_{\bullet}(-Y^h,\mathfrak{s}_h;\Gamma_{-\eta}))\cap \widecheck{HM}_{\bullet}(-Y,\mathfrak{s}';\Gamma_{-\eta})\neq\emptyset,$$
then
$$P.D.c_1(\mathfrak{s}')=\rho(P.D.c_1(\mathfrak{s})).$$
\elem

\bpf
Suppose $\mathfrak{s}\in \mathfrak{S}^*(-Y^0|-\Sigma_g)$ is any supporting spin${}^c$ structure. We define the image $\rho(P.D.c_1(\mathfrak{s}))$ as follows. Pick any spin${}^c$ structure $\mathfrak{s}_W$ on $-W$ so that the following is true

(1) We have $\widecheck{HM}(-W,\mathfrak{s}_W,\nu)\neq 0$.

(2) We have $\mathfrak{s}_W|_{-Y^0}=\mathfrak{s}$.

Then, we define $\rho(P.D.c_1(\mathfrak{s}))=P.D.c_1(\mathfrak{s}_W|_{-Y})$. We now show that this map is well defined. Suppose we have another spin${}^c$ structure $\mathfrak{s}_W'$ on $-W$ so that condition (1) and (2) are both satisfied, then we need to show that
$$P.D.c_1(\mathfrak{s}_W|_{-Y})=P.D.c_1(\mathfrak{s}_W'|_{-Y}).$$
Let $i:Y\ra W$ be the inclusion. We know that there is an exact sequence
$$H_2(W,Y)\xra{\partial} H_1(Y)\xra{i_*} H_1(W).$$
By Lemma \ref{lem_injective_inclusion_for_first_homology} and the exactness, we know that ${\rm im}(\partial)={\rm ker}(i_*)=0$.
However, clearly we have
$$\partial(P.D.c_1(\mathfrak{s}_W)-P.D.c_1(\mathfrak{s}_W'))=P.D.c_1(\mathfrak{s}_W|_{-Y})-P.D.c_1(\mathfrak{s}_W'|_{-Y}),$$
and, thus, we conclude that
$$P.D.c_1(\mathfrak{s}_W|_{-Y})=P.D.c_1(\mathfrak{s}_W'|_{-Y}).$$

The property $(*)$ follows from the construction of $\rho$ and Lemma \ref{lem_spin_c_structure_shall_extend}. The fact that $\mathfrak{PDS}^*(-Y|-\Sigma_g)\subset {\rm im}(\rho)$ follows directly from the fact that $-W$ induces an isomorphism as in Theorem \ref{thm_floer_excision}.
\epf

\bpf[Proof of proposition \ref{prop_difference_of_spin_c_structures}.]
We need a more explicit description of the map $\rho$ in Lemma \ref{lem_the_map_rho}. Using the notations in that lemma, we have a supporting spin${}^c$ structure $\mathfrak{s}$ on $-Y^0$ and a (unique) supporting spin${}^c$ structure $\mathfrak{s}_h$ on $-Y^h$. By Lemma \ref{lem_reference_closure_satisfies_proposition}, we can write
$$P.D.c_1(\mathfrak{s})=[x]+(2-2g)[s^0],$$
where $[x]\in H_1(M)\subset H_1(Y^0)$ and $s^0$ is the class as in (\ref{eq_first_homology_of_the_reference_closure}). Also, we can write
$$P.D.c_1(\mathfrak{s}_h)=[y^h]+(2-2g)[s^h],$$
where $[y^h]$ is a linear combination of the classes $[a_1],...,[b_g]$ in $H_1(Y^h)$, which is described in (\ref{eq_first_homology_of_the_mapping_torus}).

Now we claim that
$$\rho(P.D.c_1(\mathfrak{s}))=[x]+[y^h]+(2-2g)[s]\in H_1(Y).$$
This is because the cycles $x\subset Y^0$ and $x\subset Y$ co-bound annuli $[0,1]\times x$ inside $W$, $y^h\subset Y^h$ and $y^h\subset Y$ co-bound annuli $[0,1]\times y^h$ inside $W$ and $s^0\subset Y^0$, and $s^h\subset Y^h$ and $s\subset Y$ co-bound a pair of pants in side $W$. Thus, inside $W$ we can find an explicit $2$-chain $c$ so that 
$$c\cap Y^0=\partial c\cap Y^0=P.D.c_1(\mathfrak{s}),~c\cap Y^h=\partial c\cap Y^h=P.D.c_1(\mathfrak{s}_h),$$
and
$$c\cap Y=\partial c\cap Y=[x]+[y^h]+(2-2g)[s].$$
Thus, as in the proof of Lemma \ref{lem_the_map_rho}, the injectivity of $i_*$ in Lemma \ref{lem_injective_inclusion_for_first_homology} implies that
\begin{equation}\label{eq_explicit_description_of_rho}
\rho(P.D.c_1(\mathfrak{s}))=[x]+[y^h]+(2-2g)[s].
\end{equation}
With this explicit formula, Proposition \ref{prop_difference_of_spin_c_structures} follows directly.
\epf

\bcor\label{cor_rho_is_bijective}
If the inclusion $j:M\ra Y$ induces an injective homomorphism
$$j_*:H_1(M)\ra H_1(Y),$$
then the map $\rho$ in lemma \ref{lem_the_map_rho} is in fact a bijection:
$$\rho:\mathfrak{PDS}^*(-Y^0|-\Sigma_g)\ra \mathfrak{PDS}^*(-Y|-\Sigma_g).$$
\ecor

\bpf
It is straightforward from (\ref{eq_first_hmology_of_the_cobordism}) to check that when $j_*$ is injective, the inclusion $j^0:M\ra W$ also induces an injective homomorphism
$$j^0_*:H_1(M)\ra H_1(W).$$
Then, the injectivity of $\rho$ follows directly from the description of $\rho$ in (\ref{eq_explicit_description_of_rho}), since $[y^h]$ and $(2-2g)[s]$ in that formula are fixed and the only variance is $[x]$ which is represented by a cycle in $M$. Once injectivity is established, the surjectivity follows immediately from the fact that $W$ induces an isomorphism.
\epf

\bpf[Proof of proposition \ref{prop_degree_shifting_formula_for_toroidal_boundary_with_two_sutures}]
Recall that we have a balanced sutured manifold $(M,\ga)$, where $M=X(K)$ is the complement of a null-homologous knot $K\subset X$, and $\ga$ has two components. Also, we have a Seifert surface $S$ of $K$ that can be viewed as a properly embedded surface in $M$. Let $|\partial{S}\cap \ga|=2n$. For any $p$ so that $n+p$ is odd, we can perform stabilizations, as introduced in Definition \ref{defn_stabilization_grading}, and apply the construction in Definition \ref{defn_grading} to obtain a grading
$$\shm(-M,-\ga,S^p,l).$$
As in Definition \ref{defn_grading}, we can construct a marked closure
$$\mathcal{D}_{p}=(Y_p,\Sigma_g,r_p,m_p,\eta)$$
so that $S^p\subset M$ extends to a closed surface $\bar{S}^p\subset Y_p$, which leads to the grading associated to $S^p$.

We claim that the inclusion $m_p:M\ra Y_p$ for any $p$ satisfies the condition in Corollary \ref{cor_rho_is_bijective}, that is,
$$(m_p)_*:H_1(M)\ra H_1(Y_p)$$
is injective. So, the corollary applies.

To prove this claim, first note that $M=X(K)$ is the knot complement of a null-homologous knot, so we can compute directly that
$$H_1(M)=H_1(X)\oplus\lgl [\mu_K]\rgl,$$
where $\mu_K$ is a meridian circle of $K$ inside $M=X(K)$. From the construction of $Y_p$, we know that
$$Y_p=M\mathop{\cup}_{\phi}(Y^{h_p}\backslash{\rm int}(N(a_1))),$$
where $h_p:\Sigma_g\ra\Sigma_g$ is an orientation preserving diffeomorphism, $Y^{h_{p}}$ is the mapping torus of $h_p$. Also, we can compute
$$H_1(Y_p)=H_1(M)\oplus\lgl [\mu_{a}],[a_1],...,[b_g],[s_p]\rgl\slash\sim_{\phi,h_p}$$
as in (\ref{eq_first_homology_of_a_general_closure}). Thus, the relations $\sim_{\phi,h_p}$ only possibly affect the class $[\mu_k]\in H_1(M)$ but nothing in $H_1(M)$. Hence, to show that $(m_p)_*$ is injective, it is enough to show that $(m_p)_*([\mu_k])$ is of infinite order. Yet this last thing is obvious, since, inside $Y_p$, $\mu_K$ intersects $\bar{S}_p$ transversely at one point.

Thus, we get a bijection
$$\rho_p:\mathfrak{PDS}^*(-Y^0|-\Sigma_g)\ra \mathfrak{PDS}^*(-Y_p|-\Sigma_g)$$
by Corollary \ref{cor_rho_is_bijective}. Here, $(Y^0,\Sigma_g)$ or $\mathcal{D}^0=(Y^0,\Sigma_g,r^0,m^0,\eta)$ is the special (marked) closure of $(M,\ga)$ as described in Lemma \ref{lem_reference_closure_satisfies_proposition}.

Similarly, we have a surface $S^{p+2k}\subset M$, a marked closure $\mathcal{D}_{p+2k}=(Y_{p+2k},\Sigma_g,r_{p+2k},m_{p+2k},\eta)$, an extension $\bar{S}_{p+2k}$ of $S$ inside $Y_{p+2k}$, and a bijection
$$\rho_{p+2k}:\mathfrak{PDS}^*(-Y^0|-\Sigma_g)\ra \mathfrak{PDS}^*(-Y_{p+2k}|-\Sigma_g).$$
Thus we can define
$$\rho^{p}_{p+2k}=\rho_{p+2k}\circ\rho_p^{-1}:\mathfrak{PDS}^*(-Y_p|-\Sigma_g)\ra \mathfrak{PDS}^*(-Y_{p+2k}|-\Sigma_g).$$

Also, from Proposition \ref{prop_new_definition_for_canonical_map}, Lemma \ref{lem_the_map_rho}, and the functoriality of the canonical maps, we know that $\rho$ has the following significant property: If $\mathfrak{s}\in \mathfrak{S}^*(-Y_p|-\Sigma_g)$ and $\mathfrak{s}'\in \mathfrak{S}^*(-Y_{p+2k}|-\Sigma_g)$ are supporting spin${}^c$ structures so that
$$\Phi_{-\mathcal{D}_p,-\mathcal{D}_{p+2k}}(\widecheck{HM}_{\bullet}(-Y_p,\mathfrak{s};\Gamma_{-\eta}))\cap \widecheck{HM}_{\bullet}(-Y_{p+2k},\mathfrak{s}';\Gamma_{-\eta})\neq \emptyset,$$
then we have
$$P.D.c_1(\mathfrak{s}')=\rho(P.D.c_1(\mathfrak{s})).$$
From the explicit description of $\rho$ in (\ref{eq_explicit_description_of_rho}), we know that 
$\mathfrak{s}_1,\mathfrak{s}_2\in \mathfrak{S}^*(-Y_i|-\Sigma_g)$ and $\mathfrak{s}'_1,\mathfrak{s}_2'\in \mathfrak{S}^*(-Y_{i+2k}|-\Sigma_g)$ are supporting spin${}^c$ structures so that
$$\Phi_{-\mathcal{D}_p,-\mathcal{D}_{p+2k}}(\widecheck{HM}_{\bullet}(-Y_p,\mathfrak{s}_1;\Gamma_{-\eta}))\cap \widecheck{HM}_{\bullet}(-Y_{p+2k},\mathfrak{s}'_1;\Gamma_{-\eta})\neq \emptyset,$$
and
$$\Phi_{-\mathcal{D}_p,-\mathcal{D}_{p+2k}}(\widecheck{HM}_{\bullet}(-Y_p,\mathfrak{s}_2;\Gamma_{-\eta}))\cap \widecheck{HM}_{\bullet}(-Y_{p+2k},\mathfrak{s}'_2;\Gamma_{-\eta})\neq \emptyset,$$
then there exists a $1$-cycle $x\subset M$ so that
\begin{equation}\label{eq_difference_of_spin_c_structures_are_the_same_1}
	P.D.c_1(\mathfrak{s}_1)-P.D.c_1(\mathfrak{s}_2)=[x]\in H_1(Y_p),
\end{equation}
and
\begin{equation}\label{eq_difference_of_spin_c_structures_are_the_same_2}
	P.D.c_1(\mathfrak{s}_1')-P.D.c_1(\mathfrak{s}_2')=[x]\in H_1(Y_{p+2k}).
\end{equation}

Recall, in Definition \ref{defn_grading}, the grading is obtained by the evaluation of the first Chern classes of the supporting spin${}^c$ structures and by Theorem \ref{thm_well_definedness_of_the_grading}, the grading is preserved by the canonical map. Hence, the above equalities (\ref{eq_difference_of_spin_c_structures_are_the_same_1}) and (\ref{eq_difference_of_spin_c_structures_are_the_same_2}) imply that there is a fixed integer $l_0$ so that, for any $l\in \intg$, we have
$$\shm(-M,-\ga,S^p,l)=\shm(-M,-\ga,S^{p+2k},l-l_0).$$
If we go through the construction of $\rho$, we know that $\rho$ is not only independent of $l\in\intg$, but also independent of the interior of $M$ and $S$ (and is only related to the data $\partial{S}$, $p$, $k$ and $\ga$.) Thus in order to figure out the value of $k$, we can only look at the basic case where $M$ is the complement of a trefoil inside $S^3$. The convenience is that, when decomposing $(M,\ga)$ along $S$ and $-S$, the resulting sutured manifolds are both taut. 

{\bf Case 1.}
If $p<0$ and $p+2k<0$. From Lemma \ref{lem_positive_negative_stabilization_and_decomposition} and Lemma \ref{lem_decomposition_give_top_grading} we know that the top non-vanishing degree of $\shm(-M,-\ga,S^p)$ is 
$$l=\frac{n-p-1}{2}+g(S),$$
while the first non-vanishing degree of $\shm(-M,-\ga,S^{p+2k})$
is
$$l'=\frac{n-p-2k-1}{2}+g(S).$$
However, from the above discussion we know that
$$l'=l-l_0$$
so $l_0=k$. 

{\bf Case 2.} If $p=-1$ and $k=1$ or $p=1$ and $k=-1$. Then $l_0=1=k$ from proposition \ref{prop_naive_version_of_degree_shifting}.

{\bf Case 3.} If $p>0$ and $p+2k>0$. Then we can look at the surface $-S\subset M$. Note positive stabilizations of $S$ are negative stabilizations of $-S$. Hence this is reduced to case 1 and we still have $l_0=k$.

{\bf Case 4.} If $p$ and $p+2k$ are of difference sign, and is not in case 2. We can apply case 1, 2, and 3 above and conclude that $l_0=k$.

So, in summary, we always have $l_0=k$ and we conclude the proof of Proposition \ref{prop_degree_shifting_formula_for_toroidal_boundary_with_two_sutures}.
\epf

\subsection{Floer homologies on a sutured solid torus}\label{subsec_sutured_solid_torus}
As a first application of the grading shifting property, we compute the sutured monopole Floer homology of any sutured solid tori. The same result in sutured Heegaard Floer theory can be found in Juh\'asz \cite{juhasz2010polytope}.

Suppose $V=S^1\times D^2$ is a solid torus. Let $\lambda$ denote a longitude $S^1\times\{t\}$ where $t\in\partial{D}^2$ and let $\mu$ denote a meridian $\{s\}\times \partial{D}^2$ where $s\in S^1$. Suppose further $\ga$ is a suture on $V$ so that $(V,\ga)$ is a balanced sutured manifold. Then, $\ga$ is parametrized by two quantities, $n$ and $s$, where $2n$ is the number of components of $\ga$ and $s$ is the slope of the suture. In this subsection, we write the suture $\ga$ as $\ga^n_{(q,-p)}$. We write the slope $s$ as $(q,-p)$, and this is to keep our notations consistent with the ones in Honda \cite{honda2000classification}. Note $(q,-p)$ means going around longitude $-p$ times and meridian $q$ times. We always assume that $p\geq0$.

\bprop\label{prop_solid_torus_with_two_sutures}
Suppose $(V,\ga^2_{(q,-p)})$ is defined as in the above paragraph. Then, we have
$$\shm(-V,-\ga^2_{(q,-p)})=\mathcal{R}^p.$$
\eprop

\bpf
If $p=|q|$, then $p=\pm q=1$, since they are co-prime. Then, $(V,\ga_{(1,\pm1)}^2)$ is diffeomorphic to a product sutured manifold $(A\times[-1,1],\partial{A}\times\{0\})$, where $A$ is an annulus. Thus, we know
$$\shm(-V,-\ga_{(1,-1)}^2)\cong\mathcal{R}.$$

From now on, we assume that $p>q>0$. If not we can achieve this assumption by applying diffeomorphisms of the solid torus $V$. We want to re-interpret the by-pass exact triangle as follows: We have a basic by-pass exact triangles

\begin{equation}\label{eq_initial_exact_triangle}
\xymatrix{
&\shm(-V,-\ga^2_{(1,-1)})\ar[rd]^{\psi_{-,2}}&\\
\shm(-V,-\ga^2_{(1,0)})\ar[ru]^{\psi_{-,1}}&&\shm(-V,-\ga^2_{(0,-1)})\ar[ll]^{\psi_{-,0}}\\
}
\end{equation}

Here $\psi_{-,0}=\psi_{-,0}^{\infty}$, $\psi_{-,1}=\psi_{-,1}^{0}$, and $\psi_{-,2}=\psi_{-,\infty}^{1}$, under the notations in (\ref{eq_by_pass_on_knot_complement}). 

Recall, from Subsection \ref{subsec_contact_elements_and_contact_structures}, that the map $\psi_{-,1}$ (as well as the other two) is interpreted as a gluing map: Suppose we have $(-V,-\ga^2_{(1,0)})$ and an identification $T^2=S^1\times \partial{D}^2$, then we can glue $[0,1]\times T^2$ to $V$ via the identification $id:\partial{V}=S^1\times \partial{D}^2\xra{=}\{0\}\times T^2$. Suppose $\{0\}\times T^2$ is equipped with the suture $\ga^2_{(1,0)}$, and $T^{2}\times\{1\}$ is equipped with the suture $\ga^2_{(1,-1)}$, then we can identify $(V,\ga^2_{(1,-1)})$ with $(V\cup [0,1]\times T^2,\ga^2_{(1,-1)})$. There exists a compatible contact structure $\xi_{-,1}$ on $([0,1]\times T^2,\ga^2_{(1,0)}\cup\ga^2_{(1,-1)})$ so that we have
$$\psi_{-,0}=\Phi_{\xi_{-,0}}:\shm(-V,-\ga^2_{(1,0)})\ra \shm(-V,-\ga^2_{(1,-1)}).$$

When dealing with other sutures, we can also glue $(T^2\times[0,1],\ga^2_{(1,0)}\cup\ga^2_{(1,-1)})$ to $V$, but along a diffeomorphism
$$g:\{0\}\times T^2\ra \partial V,$$
instead of the identity map. Such a map needs to be orientation preserving and, hence, is parametrized by an element in $SL_2(\intg)$. We can pick the map $g$ corresponding to the matrix
\begin{equation*}
A=\left(
\begin{array}{cc}
	q-q'&-q'\\
	p'-p&p'
\end{array}
\right)\in SL_2(\intg),
\end{equation*}
where $p'q-pq'=1$, $p'\leq p$, $q'\leq q$, $q''=p-p'$, and $p''=p-p'$. (Such $p',q',p'',q''$ are unique.)

Then, the suture $\ga^2_{(1,0)}$ on $T^2\times\{0\}$ is glued to $\ga^2_{(q,-p)}$ on $\partial{V}$ and the suture $\ga^2_{(1,-1)}$ on $T^2\times\{1\}$ now becomes the suture $\ga^2_{(q',-p')}$. As in Formula (\ref{eq_initial_exact_triangle}), they still fit into an exact triangle

\begin{equation}\label{eq_exact_triangle_for_two_sutures}
\xymatrix{
&\shm(-V,-\ga^2_{(q,-p)})\ar[rd]^{\psi_{-,2}}&\\
\shm(-V,-\ga^2_{(q'',-p'')})\ar[ru]^{\psi_{-,1}}&&\shm(-V,-\ga^2_{(q',-p')})\ar[ll]^{\psi_{-,0}}\\
}
\end{equation}

We claim that $\psi_{-,0}=0$. Let $D_{p}$ be a meridian disk of $V$ which intersects $\ga^2_{(q,-p)}$ at $2p$ points, then, from a similar argument as in Proposition \ref{prop_general_degree_shifting_for_knot_complement} (which we will prove later), we have
$$\psi_{-,0}(\shm(-V,-\ga^2_{(q',-p')},D_{p'}^{-(p-p')},i))\subset\shm(-V,-\ga^2_{(q'',-p'')},D_{p''}^{+(p-p'')},i)$$
for any $i\in\intg$.

We only deal with the case when $p'$ is odd and $p''$ is even. Other cases are similar. From the construction of the grading in Definition \ref{defn_grading}, we know that there is a suitable marked closure $\mathcal{D}_{p'}=(Y_{p'},R,r,m,\eta)$ and a closed surface $\bar{D}_{p'}\subset Y_{p'}$ so that the grading is defined via the evaluations of the first Chern classes of spin${}^c$ structures on the fundamental class of $\bar{D}_{p'}$. From the construction, we know that
$$\chi(\bar{D}_{p'})=\chi(D_{p'})-p'=1-p'.$$
Hence, the adjunction inequality in Lemma \ref{lem_adjunction_inequality} implies that
$$\shm(-V,-\ga^2_{(q',-p')},D_{p'},i)=0$$
if $i<\frac{1-p'}{2}$.
Then, from the grading shifting property in Proposition \ref{prop_degree_shifting_formula_for_toroidal_boundary_with_two_sutures}, we know that
$$\shm(-V,-\ga^2_{(q',-p')},D_{p'}^{-p''},i)=\shm(-V,-\ga^2_{(q',-p')},D_{p'},i+(\frac{p''}{2})).$$
Thus, we know
\begin{equation}\label{eq_vanishing_at_low_grading}
\shm(-V,-\ga^2_{(q',-p')},D_{p'}^{-p''},i)=0
\end{equation}
if $i<\frac{1-p'+p''}{2}$. Note, by definition, $p''=p-p'$.

The above argument for $D_{p'}$ applies to $D_{p''}^+$ as well. Note $p''$ is assumed to be even, so we need to perform a positive stabilization on $D_{p''}$ to construct the grading. The adjunction inequality in Lemma \ref{lem_adjunction_inequality} again implies that
\begin{equation}\label{eq_adjunction_inequality_D_p''_+}
\shm(-V,-\ga^2_{(q'',-p'')},D_{p''}^+,i)=0
\end{equation}
if $i>\frac{p''}{2}$. However, from Lemma \ref{lem_decomposition_give_top_grading}, we know that
$$\shm(-V,-\ga^2_{(q'',-p'')},D_{p''}^+,\frac{p''}{2})\cong\shm(M',\ga'),$$
where $(M',\ga')$ is the result of doing a sutured manifold decomposition on $(-V,-\ga^2_{(q'',-p'')})$ along the surface $D_{p''}^+$. From Lemma \ref{lem_positive_negative_stabilization_and_decomposition}, we know that
\begin{equation}\label{eq_top_grading_vanishing}
\shm(-V,-\ga^2_{(q'',-p'')},D_{p''}^+,\frac{p''}{2})\cong\shm(M',\ga')=0.
\end{equation}
The grading shifting property in Proposition \ref{prop_degree_shifting_formula_for_toroidal_boundary_with_two_sutures}, then, implies
$$\shm(-V,-\ga^2_{(q'',-p'')},D_{p''}^{+p'},i)=\shm(-V,-\ga^2_{(q'',-p'')},D_{p''}^+,i-\frac{p'-1}{2}).$$
The above equality, together with (\ref{eq_adjunction_inequality_D_p''_+}) and (\ref{eq_top_grading_vanishing}), implies that

$$\shm(-V,-\ga^2_{(q'',-p'')},D_{p''}^{+p'},i)=0$$
if $i \geq\frac{1-p'+p''}{2}$. Compare this with (\ref{eq_vanishing_at_low_grading}), we can see that $\psi_{-,0}=0.$

Once we conclude that $\psi_{-,0}=0$, we can compute
$\shm(-V,-\ga^2_{(q,-p)})$ by the induction, and Proposition \ref{prop_solid_torus_with_two_sutures} follows. 
\epf

\brem
As in Honda \cite{honda2000classification}, the two slopes $(q',-p')$ and $(q'',-p'')$ can be written out explicitly in terms of the continued fraction of $(q,-p)$. Note we have assumed $p>q$. Suppose
$$-\frac{p}{q}=r_1-\frac{1}{r_2-\frac{1}{r_3-...}},$$
where it is a finite continued fraction, and $r_j<-1$ for all $j$. We can write
\begin{equation}\label{eq_expression_in_terms_of_continued_fraction}
-\frac{p}{q}=[r_1,r_2,...,r_k].
\end{equation}
Under this notation, we have
$$-\frac{p'}{q'}=[r_1,r_2...,r_{k-1}],~-\frac{p''}{q''}=[r_1,r_2...,r_{k-1}+1],$$
and in the above notation, we identify $[r_1,...,r_{j-1},r_j,-1]$ with $[r_1,...,r_{j-1},r_j+1].$
\erem

\begin{figure}[h]
\centering
\begin{overpic}[width=5in]{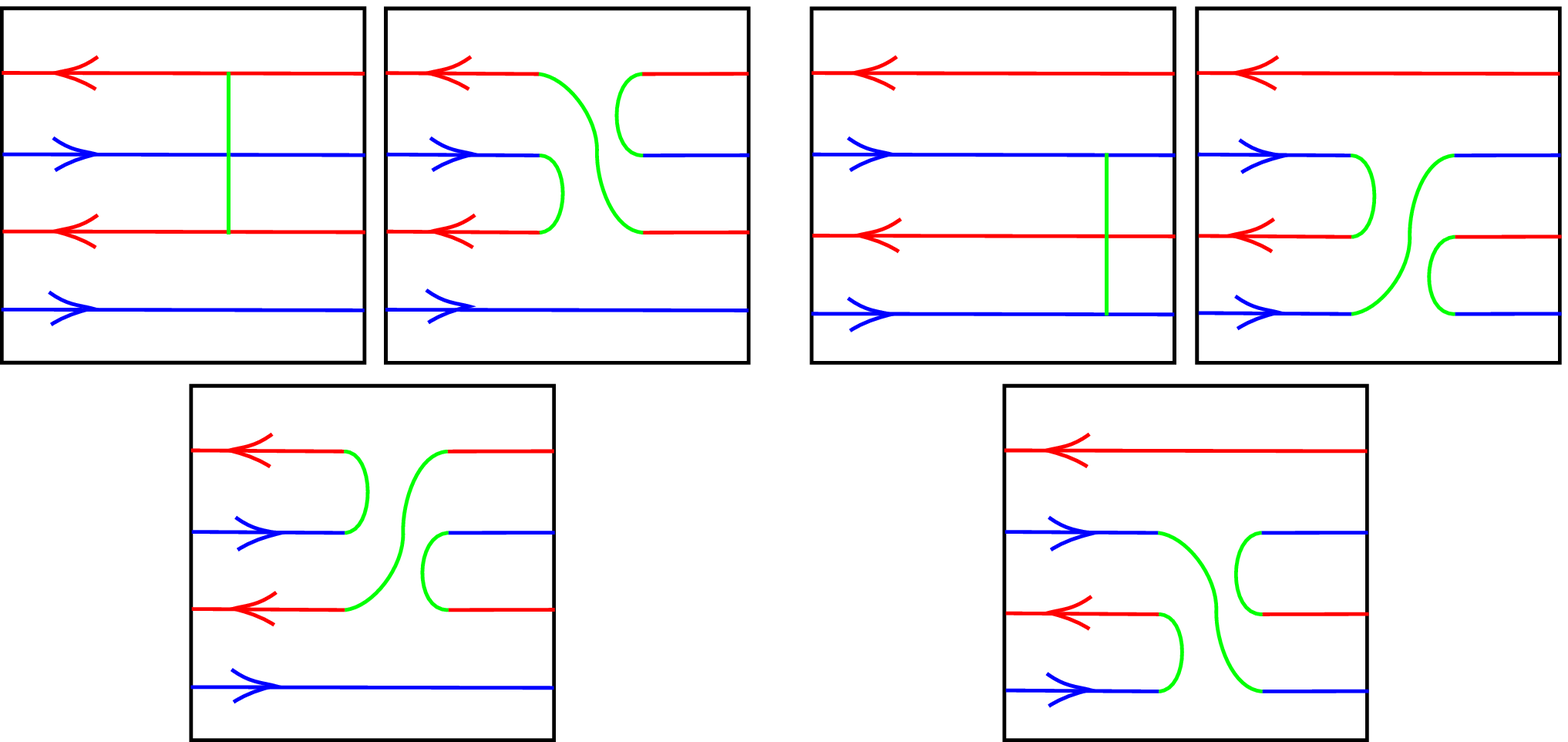}
	\put(15,-4){Positive by-passes}
	\put(65,-4){Negative by-passes}
\end{overpic}
\vspace{0.05in}
\caption{The positive and negative by-passes.}\label{fig_by_pass_5}
\end{figure}

Now we deal with the general sutures $\ga^{n}_{(q,-p)}$ for $n>1$. There are two types by-passes relating $(V,\ga^{2n+2}_{(q,-p)})$ and $(V,\ga^{2n}_{(q,-p)})$. We call them positive and negative by-passes according to Figure \ref{fig_by_pass_5}. They give rise to by-pass exact triangles:

\begin{equation}
\xymatrix{
&\shm(-V,-\ga^{2n+2}_{(q,-p)})\ar[rd]^{\psi_{\pm,n}^{n+1}}&\\
\shm(-V,-\ga^{2n}_{(q,-p)})\ar[ru]^{\psi_{\pm,n+1}^n}&&\shm(-V,-\ga^{2n}_{(q,-p)}\ar[ll]^{\psi_{\pm,n}^n})
}	
\end{equation}

\brem\label{rem_no_unique_by_passes}
Unlike the case of two sutures where there are exactly two different possibilities of by-passes, in the case where $\ga$ has more than two components, positive and negative by-passes are not unique. Here, we just pick two specific by-passes so that they are 'adjacent' to each other. This is crucial to the proof of Lemma \ref{lem_composition_to_be_identity}.
\erem

\blem\label{lem_composition_to_be_identity}
For any $n\in\intg$ and slope $(q,-p)$, we have
$$\psi_{-,n}^{n+1}\circ\psi_{+,n+1}^n=\psi_{+,n}^{n+1}\circ\psi_{-,n+1}^n=id:\shm(-V,-\ga^{2n}_{(q,-p)})\ra\shm(-V,-\ga^{2n}_{(q,-p)}).$$
\elem

\bpf
We will only prove that $\psi_{-,n}^{n+1}\circ\psi_{+,n+1}^n=id$. The other is the same.

From  \cite{baldwin2016contact} or \cite{ozbagci2011contact} we know that a by-pass attached along an arc $\al$ can be thought of as attaching a pair of contact $1$-handle and $2$-handle. The contact one handle is attached along the two end points $\partial{\al}$ while the contact two handle is attached along a Legendrian curve 
$$\be=\al\cup\al',$$
where $\al'$ is an arc on the contact $1$-handle intersecting the dividing set once.

\begin{figure}[h]
\centering
\begin{overpic}[width=5in]{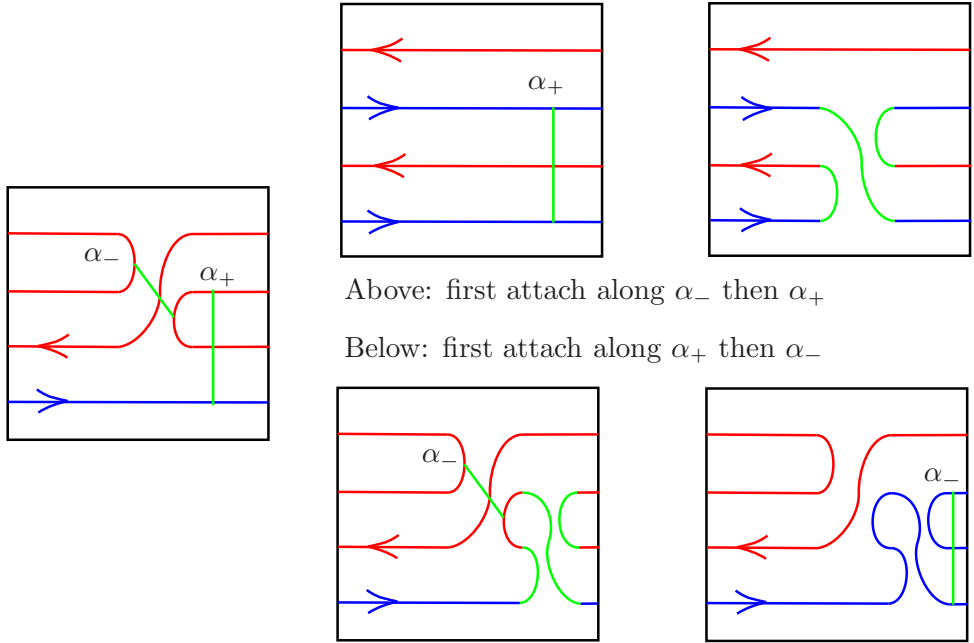}
	\put(35,36){Above: first attach along $\al_-$ then $\al_+$}
	\put(35,30){Below: first attach along $\al_+$ then $\al_-$}
	\put(20,38){$\al_+$}
	\put(8,40){$\al_-$}
	\put(54,58){$\al_+$}
	\put(43,19){$\al_-$}
	\put(95,17){$\al_-$}
\end{overpic}
\caption{Reversing the order of by-pass attachments. Bottom right picture: we can isotope $\al_-$ to this new position, where we can see directly that the by-pass is trivial.}\label{fig_by_pass_exchange}
\end{figure}

Now $\psi_{-,n}^{n+1}\circ\psi_{+,n+1}^n$ corresponds to first attaching a by-pass along $\al_+$ and then attaching another one along $\al_-$, as in Figure \ref{fig_by_pass_exchange}. However, in terms of contact handle attachments, the two pairs of handles are disjoint from each other, so we can reverse the order of attachments: Instead, we can first attach a by-pass along $\al_-$ and then along $\al_+$. If we attach a by-pass along $\al_-$ first, we can see from Figure \ref{fig_by_pass_exchange} that this is a trivial by-pass as discussed in Honda \cite{honda2002gluing}. In that paper, it is proved that a trivial by-pass does not change the contact structure. From theorem \ref{thm_gluing_map}, we conclude that a trivial by-pass induces the identity map. Then, the second by-pass attached along $\al_+$ is also trivial and, hence, again induces the identity map. Thus, we conclude that $\psi_{-,n}^{n+1}\circ\psi_{+,n+1}^n=id$.
\epf

\bcor
We know that
$$\shm(-V,-\ga_{(q,-p)}^{2n})\cong\mathcal{R}^{(2^{n-1}\cdot p)}.$$
\ecor
\bpf
From Lemma \ref{lem_composition_to_be_identity}, we know that $\psi_{\pm,n}^{n+1}$ is surjective while $\psi_{\pm,n+1}^n$ is injective. Hence, we can conclude the statement by using the by-pass exact triangles and the induction.
\epf

\bcor
We have
$$|\pi_0({\rm Tight}(V,\ga_{(q,-p)}^{2n}))|\geq 2^{n-1}\cdot|r_1+1|\cdot...\cdot|r_{k-1}+1|\cdot|r_k|.$$
\ecor

\bpf
First assume $n=1$. In \cite{honda2000classification}, Honda explained how to construct any compatible tight contact structures on a sutured solid torus: First we start with the standard tight contact structure on $(V,\ga_{(1,-1)}^2)$. Then, we can glue $k$ different layers $T^2\times[i-1,i]$, for $1\leq i\leq k$, to $V$, so that, on $T^2\times[i-1,i]$, $T^2\times\{i-1\}$ has the dividing set $\ga_{(1,-1)}^2$, while $T^2\times\{i\}$ has the dividing set $\ga_{(1,1-r_{i})}^2$. We glue $T^2\times\{0\}$ to $\partial{V}$ via identity, while glue $T^2\times\{i\}\subset T^2\times[i,i+1]$ to $T^2\times\{i\}\subset T^2\times[i-1,i]$ so that the dividing sets on these two surfaces are identified. 

Each layer $T^2\times[i-1,i]$ is further decomposed into the composition of $-1-r_i$ (or $-r_k$ for the last layer) many by-passes. There are two by-passes: One corresponds to the map $\psi_{-,1}$ in formula (\ref{eq_exact_triangle_for_two_sutures}), and the other corresponds to some $\psi_{+,1}$ in a similar by-pass exact triangle. Use the inductive step as introduced in \cite{honda2000classification}, which Honda used to construct tight contact structures on a sutured solid torus, we see that all the contact structures that Honda constructed have distinct contact elements. Hence, there are at least $|r_1+1|\cdot...\cdot|r_{k-1}+1|\cdot|r_k|$ many different contact structures. 

When $n$ is bigger than $1$, we proceed by induction. Suppose, for $n=l$, there are at least $m_l=2^{l-1}\cdot|r_1+1|\cdot...\cdot|r_{k-1}+1|\cdot|r_k|$ many different non-zero contact elements $\psi_{\xi_1},...,\psi_{\xi_{m_l}}\in\shm(-V,\ga_{(q,-p)}^{2l})$. From Lemma \ref{lem_composition_to_be_identity}, we know that $\psi_{+,l+1}^l$ and $\psi_{-,l+1}^l$ are both injective,
$$\psi_{\pm,l}^{l+1}\circ\psi_{\pm,l+1}^l=0,~{and}~\psi_{\mp,l}^{l+1}\circ\psi_{\pm,l+1}^l=id.$$
The first equality is due to the exactness of the by-pass triangle, and the second is again Lemma \ref{lem_composition_to_be_identity}. Hence, we know that, inside $\shm(-V,\ga_{(q,-p)}^{2l+2})$, there are at least $m_{l+1}=2^{l}\cdot|r_1+1|\cdot...\cdot|r_{k-1}+1|\cdot|r_k|$ many different contact elements 
$$\psi_{\pm,l+1}^{l}(\phi_{\xi_{1}}),...,\psi_{\pm,l+1}^{l}(\phi_{\xi_{m_l}}).$$
Hence, we are done.
\epf

\brem
When $n=1$, the above argument gives an alternative way to provide a tight lower bound of $|\pi_0(Tight(V,\ga^2_{(q,-p)}))|$, which is originally done by Honda \cite{honda2000classification}.

When $n>1$, as mentioned in Remark \ref{rem_no_unique_by_passes}, there are not just two by-passes, so this lower bound, a priori, need not to be tight. However, one could try to study the impact of all other by-pass attachments to see if we could improve the lower bound.
\erem

\brem
We can use a meridian disk of the solid torus to define a grading on $\shm(-V,-\ga_{(q,-p)}^{2n})$. The above method is also capable of computing the graded homology.
\erem

\section{The direct system and the direct limit}\label{seciton_direct_system}
\subsection{The construction}\label{subsec_construction}
Suppose $Y$ is a closed oriented $3$-manifold, and $K\subset Y$ is an oriented knot with a Seifert surface $S\subset Y$. Suppose further that $p\in K$ is a fixed base point and $\varphi:S^1\times D^2\hookrightarrow Y$ is an embedding as in Subsection \ref{subsec_naturality}, i.e., we require that
$$\varphi(S^1\times\{0\})=k,~{\rm and}~\varphi(\{1\}\times\{0\})=p.$$
Then, we have a $3$-manifold with boundary $Y_{\varphi}=Y\backslash {\rm int}({\rm im}(\varphi))$. The Seifert surface $S$ induces a framing on $\partial{Y_{\varphi}}$. We call the meridian $\mu_{\varphi}$ and the longitude $\lambda_{\varphi}$. Let $\Ga_{n,\varphi}$ be a collection of two disjoin parallel oppositely oriented simple closed curves on $\partial{Y_{\varphi}}$, each of class $\pm(\lambda_{\varphi}-n\mu_{\varphi})$. Then, we have a balanced sutured manifold $({Y_{\varphi}},\Ga_{n,\varphi})$.

Suppose $\varphi'$ is another embedding, then we also have $(({Y_{\varphi'}},\Ga_{n,\varphi'}))$. Suppose $f_t$ is the ambient isotopy defined as in Subsection \ref{subsec_naturality}, relating $\varphi$ and $\varphi'$. We have the following lemma.

\blem\label{lem_ambiguity_from_choices_of_knot_complements}
The diffeomorphism $f_1$ is a diffeomorphism from $({Y_{\varphi}},\Ga_{n,\varphi})$ to $({Y_{\varphi'}},\Ga_{n,\varphi'})$.
\elem

\bpf
It is enough to show that $f_1$ sends the framing $(\mu_{\varphi},\lambda_{\varphi})$ on $\partial{Y_{\varphi}}$ to the framing $(\mu_{\varphi'},\lambda_{\varphi'})$ on $\partial{Y_{\varphi'}}$.

By construction, $f_1$ sends $\mu_{\varphi}$ to $\mu_{\varphi'}$. $f_1$ must also preserve $\lambda_{\varphi}$, since $f_t$ is an isotopy, and $\lambda_{\varphi}$ can be characterized by the fact that it represents a generator of the map
$$i_*:H_1(\partial{Y_{\varphi}})\ra H_1(Y_{\varphi}),$$
where $i: \partial{Y_{\varphi}}\ra Y_{\varphi}$ is the inclusion.
\epf

\bcor\label{cor_well_definedness_of_gamma_n_suture}
There is a transitive system (of projective transitive systems) 
$$\{\shm({Y_{\varphi}},\Ga_{n,\varphi})\}~{\rm and}~\{\Psi_{\varphi,\varphi'}=\shm(f_1)\}.$$
 So, we obtain a canonical module $\shm(Y,K,p,n)$ associated to the quadruple $(Y,K,p,n)$. 
\ecor

Once Lemma \ref{lem_ambiguity_from_choices_of_knot_complements}, we can fix a knot complement to study with. Suppose $Y(K)=Y\backslash {\rm int}(N(K))$ be a knot complement and let $\lambda$ and $\mu$ be the longitude and meridian, respectively, with respect to the framing on $\partial{Y(K)}$ induced by the Seifert surface $S$. For any $n\in\intg_+$, use $\Ga_n$ to denote the suture on $\partial Y(K)$ consisting of a pair of simple closed curves of class $\pm(\lambda-n\mu)$, and use $\Ga_{\infty}$ to denote the suture on $\partial Y(K)$ consisting of a pair of meridians.

\bdefn[Kronheimer and Mrowka \cite{kronheimer2010knots}, or Baldwin and Sivek \cite{baldwin2015naturality}]\label{defn_original_knot_Floer_homology}
Define
$$\underline{\rm KHM}(Y,K,p)=\shm(Y(K),\Ga_{\infty}).$$
\edefn

\bdefn\label{defn_minus_version_of_knot_monopoles}
Define the {\it minus version of monopole knot Floer homology} of a based knot $K\subset -Y$, which is denoted by $\khm(-Y,K,p)$, to be the direct limit of the direct system
$$...\ra\shm(-Y(K),\Ga_n)\xra{\psi^n_{-,n+1}}\shm(-Y(K),\Ga_{n+1})\ra...$$
Here, the maps $\psi^n_{-,n+1}$ are defined in the exact triangle (\ref{eq_by_pass_on_knot_complement}).
By Corollary \ref{cor_commutative_diagram}, the maps $\{\psi_{+,n+1}^n\}_{n\in\intg_+}$ induce a map on $\khm$, which we call $U$:
$$U:\khm(-Y,K,p)\ra\khm(-Y,K,p).$$
\edefn

Next, we construct a grading on the direct limit $\khm(-Y,K,p)$. Suppose $S_n$ is the Seifert surface of $K$ so that $S_n$ intersects $\Ga_n$ at $2n$ points. Then, we have the following proposition.

\bprop\label{prop_general_degree_shifting_for_knot_complement}
Suppose $n$ is even, then, for any $i\in\intg$, we have
$$\psi^n_{\pm,n+1}(\shm(-Y(K),-\Ga_n,S^{\pm}_n,i))\subset\shm(-Y(K),-\Ga_{n+1},S_{n+1},i).$$

Suppose $n$ is odd, then we have for any $i\in\intg$
$$\psi^n_{\pm,n+1}(\shm(-Y(K),-\Ga_n,S^{\pm2}_{n},i))\subset\shm(-Y(K),-\Ga_{n+1},S^{\pm}_{n+1},i).$$
\eprop

\begin{figure}[h]
\centering
\begin{overpic}[width=4.5in]{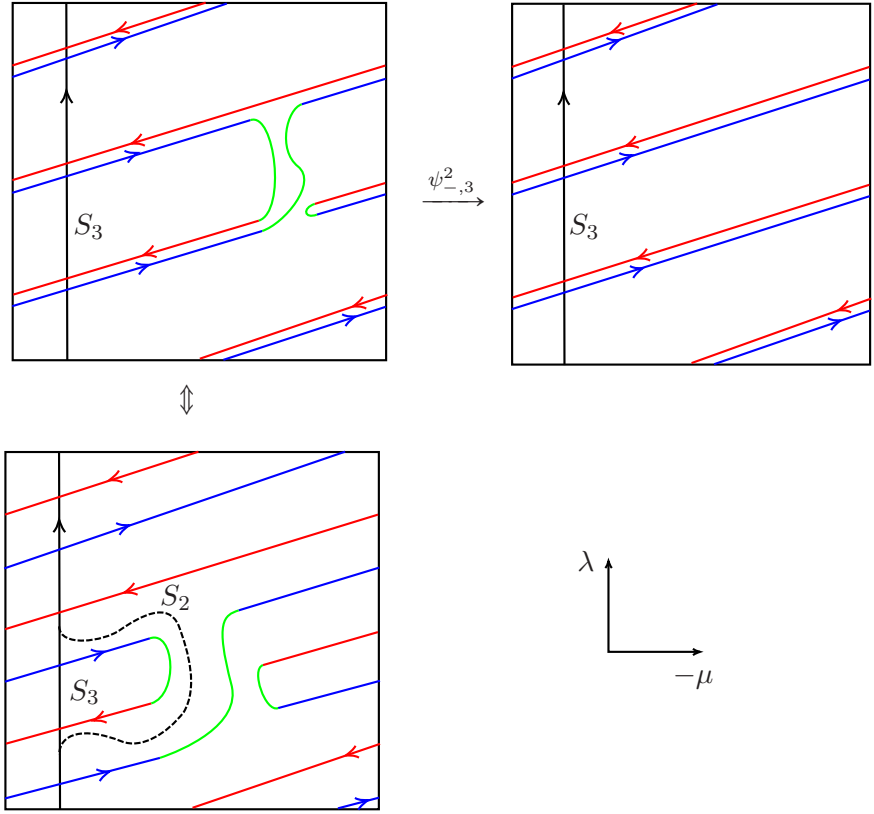}
	\put(66,28){$\lambda$}
	\put(77,15){$-\mu$}
	\put(7.5,13){$S_3$}
	\put(8,67){$S_3$}
	\put(65,67){$S_3$}
	\put(18,24){$S_{2}$}
	\put(48,70){$\xra{\psi_{-,3}^2}$}
	\put(20,47){$\Updownarrow$}
\end{overpic}
\vspace{0.05in}
\caption{The solid vertical arc represents the surface $S_3=S_{2}^-$ and the dashed arc represents $S_{2}$.}\label{fig_by_pass_4}
\end{figure}

\bpf
We only prove the proposition for $\phi_{-,n+1}^n$ with $n$ even. Other cases are similar. In Figure \ref{fig_by_pass_4}, it is clear that the surface $S_{n+1}\subset(Y(K),\Ga_{n})$ can also be obtained from the surface $S_{n}$ by a negative stabilization:
$$S_{n+1}=S^{-}_n.$$

Thus, for any $i\in\intg$, we have
\beq
\shm(-Y(K),-\Ga_n,S^-_n,i)=\shm(-Y(K),-\Ga_n,S_{n+1},i).
\eeq

For $S_n^{-}=S_{n+1}\subset (Y(K),\Ga_{n})$, we can choose some auxiliary data to construct a marked closure
$$\mathcal{D}_{n}^-=(Y_{n}^-,R,r_{n},m_{n},\eta),$$
so that $S_n^-$ extends to a closed surface $\bar{S}_{n}^-\subset Y_{n}^-$ and it induces a grading on $\shm(-Y(K),-\Ga_n)$ that is exactly the one associated to $S_n^-$. (See Definition \ref{defn_grading}.) 

We can obtain $(Y(K),\Ga_{n+1})$ by attaching a by-pass disjoint from $S_{n+1}=S_n^-$. From Baldwin and Sivek \cite{baldwin2016contact}, we know the map $\phi_{-,n+1}^n$ associated to the by-pass can be described as follows: There is a curve $\be\subset {\rm}(m_n(Y(K)))\subset Y_{n}^-$ so that a $0$-framed Dehn surgery on $\be$, with respect to the $\partial{Y(K)}$ framing, will result in a $3$-manifold $Y_{n+1}$. Since $\be$ is disjoint from ${\rm im}(r_n)$, the data $R$, $r_n$ and $\eta$ survive and we get a marked closure
$$\mathcal{D}_{n+1}=(Y_{n+1},R,r_{n+1},m_{n+1},\eta)$$
which is a marked closure of $(Y(K),\Ga_{n+1})$. The surgery description gives rise to a cobordism $W$ from $Y_{n}^-$ to $Y_{n+1}$ and the cobordism map associated to this cobordism induces the by-pass attaching map $\phi_{-,n+1}^n$.

It is a key observation that the surface $S_n^-=S_{n+1}$ is disjoint from the region we attach the by-pass and, hence, is disjoint from the curve $\be$ along which we perform the Dehn surgery. As a result, the surface $\bar{S}_n^-$ remains as a closed surface $\bar{S}_{n+1}\subset Y_{n+1}$ and induces a grading on $\shm(Y(K),\Ga_{n+1})$. It is clear that the grading induced by $\bar{S}_{n+1}$ is nothing but the one associated to the surface $S_{n+1}\subset (Y(K),\Ga_{n+1})$ as in Definition \ref{defn_grading}.

There is a product cobordism $[0,1]\times \bar{S}^-_{n}\subset W$, from $\bar{S}^-_{n}\subset Y_{n}^-$ to $\bar{S}_{n+1}\subset Y_{n+1}$, and, thus, we conclude that
$$\phi_{-,n+1}^n(\shm(Y(K),\Ga_n,S_{n}^-,i))\subset \shm(Y(K),\Ga_{n+1},S_{n+1},i).$$
This concludes the proof of Proposition \ref{prop_general_degree_shifting_for_knot_complement}. 
\epf

The following Figures \ref{fig_even_to_odd} and \ref{fig_odd_to_even} might be helpful for figuring out how do $\psi_{+,n+1}^n$ and $\psi_{-,n+1}^n$ change the gradings. In the figures, $k'=k+g(S)$.

\begin{figure}[h]
\centering
	\begin{picture}(300,180)
	\put(-5,150){$k'$}
	\put(-5,125){$k'-1$}
	\put(-5,100){$k'-2$}
	\put(-5,75){\vdots}
	\put(-5,50){$2-k'$}
	\put(-5,25){$1-k'$}
	\put(-5,0){$-k'$}
	\put(30,175){$D^{-}_{2k}$}
	\put(70,175){$D_{2k+1}$}
	\put(30,150){$\bigcirc$}
	\put(70,150){$\bigcirc$}
	\put(44,153){\vector(1,0){24}}
	\put(30,125){$\bigcirc$}
	\put(70,125){$\bigcirc$}
	\put(44,128){\vector(1,0){24}}
	\put(30,100){$\bigcirc$}
	\put(70,100){$\bigcirc$}
	\put(44,103){\vector(1,0){24}}
	\put(30,50){$\bigcirc$}
	\put(70,50){$\bigcirc$}
	\put(44,53){\vector(1,0){24}}
	\put(30,25){$\bigcirc$}
	\put(70,25){$\bigcirc$}
	\put(44,28){\vector(1,0){24}}
	\put(70,0){$\bigcirc$}

	\put(150,175){$D^{-}_{2k}$}
	\put(190,175){$D^+_{2k}$}
	\put(230,175){$D_{2k+1}$}
	\put(190,125){$\bigcirc$}
	\put(230,125){$\bigcirc$}
	\put(204,128){\vector(1,0){24}}
	\put(190,100){$\bigcirc$}
	\put(230,100){$\bigcirc$}
	\put(204,103){\vector(1,0){24}}
	\put(190,50){$\bigcirc$}
	\put(230,50){$\bigcirc$}
	\put(204,53){\vector(1,0){24}}
	\put(190,25){$\bigcirc$}
	\put(230,25){$\bigcirc$}
	\put(204,28){\vector(1,0){24}}
	\put(190,0){$\bigcirc$}
	\put(230,0){$\bigcirc$}
	\put(204,3){\vector(1,0){24}}
	\put(150,150){$\bigcirc$}
	\dashline{3}(163,151)(189,129)
	\put(230,150){$\bigcirc$}

	\put(150,125){$\bigcirc$}
	\dashline{3}(163,126)(189,104)
	\put(150,100){$\bigcirc$}
	\dashline{3}(163,101)(189,79)
	\dashline{3}(163,76)(189,54)
	\put(150,50){$\bigcirc$}
	\dashline{3}(163,51)(189,29)
	\put(150,25){$\bigcirc$}
	\dashline{3}(163,26)(189,4)
	
	\put(34,75){\vdots}
	\put(74,75){\vdots}
	\put(154,75){\vdots}
	\put(194,75){\vdots}
	\put(234,75){\vdots}

\end{picture}
\caption{The maps $\phi_{\pm}$ from $\shm(-Y(K),-\Ga_{2k})$ to $\shm(-Y(K),-\Ga_{2k+1})$. The map $\phi_{-,2k+1}^{2k}$ is depicted on the left and $\phi_{+,2k+1}^{2k}$ on the right. They are represented by the solid arrows. The circles $\bigcirc$ denote the graded homologies. The dashed lines represent the grading shifting when using different surfaces to construct the grading.}\label{fig_even_to_odd}
\end{figure}
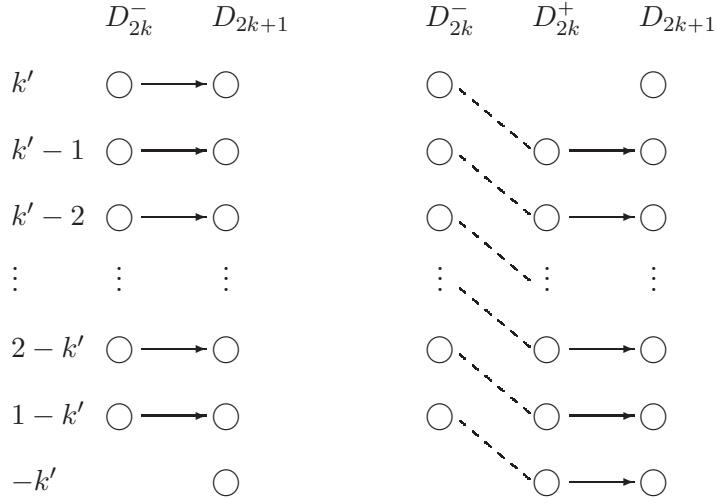

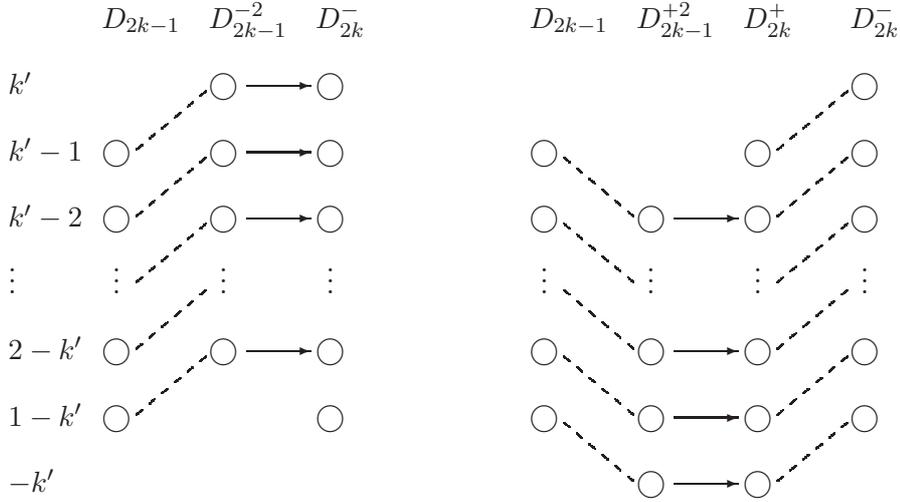
\begin{figure}[h]
\centering
	\begin{picture}(350,180)
	
	\put(-5,150){$k'$}
	\put(-5,125){$k'-1$}
	\put(-5,100){$k'-2$}
	\put(-5,75){\vdots}
	\put(-5,50){$2-k'$}
	\put(-5,25){$1-k'$}
	\put(-5,0){$-k'$}
	
	\put(30,175){$D_{2k-1}$}
	\put(70,175){$D^{-2}_{2k-1}$}
	\put(110,175){$D^-_{2k}$}
	
	\put(70,150){$\bigcirc$}
	\put(110,150){$\bigcirc$}
	\put(84,153){\vector(1,0){24}}
	\put(70,125){$\bigcirc$}
	\put(110,125){$\bigcirc$}
	\put(84,128){\vector(1,0){24}}
	\put(70,100){$\bigcirc$}
	\put(110,100){$\bigcirc$}
	\put(84,103){\vector(1,0){24}}
	\put(70,50){$\bigcirc$}
	\put(110,50){$\bigcirc$}
	\put(84,53){\vector(1,0){24}}
	\put(110,25){$\bigcirc$}

	\put(30,125){$\bigcirc$}
	\dashline{3}(43,129)(69,151)
	\put(30,100){$\bigcirc$}
	\dashline{3}(43,104)(69,126)
	\dashline{3}(43,79)(69,101)
	\put(30,50){$\bigcirc$}
	\dashline{3}(43,54)(69,76)
	\put(30,25){$\bigcirc$}
	\dashline{3}(43,29)(69,51)
	
	\put(190,175){$D_{2k-1}$}
	\put(230,175){$D^{+2}_{2k-1}$}
	\put(270,175){$D^+_{2k}$}
	\put(310,175){$D^-_{2k}$}
	
	\put(270,125){$\bigcirc$}
	\put(230,100){$\bigcirc$}
	\put(270,100){$\bigcirc$}
	\put(244,103){\vector(1,0){24}}
	\put(230,50){$\bigcirc$}
	\put(270,50){$\bigcirc$}
	\put(244,53){\vector(1,0){24}}
	\put(230,25){$\bigcirc$}
	\put(270,25){$\bigcirc$}
	\put(244,28){\vector(1,0){24}}
	\put(230,0){$\bigcirc$}
	\put(270,0){$\bigcirc$}
	\put(244,3){\vector(1,0){24}}

	\put(190,125){$\bigcirc$}
	\dashline{3}(203,126)(229,104)
	\put(190,100){$\bigcirc$}
	\dashline{3}(203,101)(229,79)
	\dashline{3}(203,76)(229,54)
	\put(190,50){$\bigcirc$}
	\dashline{3}(203,51)(229,29)
	\put(190,25){$\bigcirc$}
	\dashline{3}(203,26)(229,4)
	
	\dashline{3}(283,129)(309,151)
	\put(310,150){$\bigcirc$}
	\dashline{3}(283,104)(309,126)
	\put(310,125){$\bigcirc$}
	\dashline{3}(283,79)(309,101)
	\put(310,100){$\bigcirc$}
	\dashline{3}(283,54)(309,76)
	\dashline{3}(283,29)(309,51)
	\put(310,50){$\bigcirc$}
	\dashline{3}(283,4)(309,26)
	\put(310,25){$\bigcirc$}
	
	\put(34,75){\vdots}
	\put(74,75){\vdots}
	\put(114,75){\vdots}
	\put(194,75){\vdots}
	\put(234,75){\vdots}
	\put(274,75){\vdots}
	\put(314,75){\vdots}
	
\end{picture}
\caption{The maps $\phi_{\pm}$ from $\shm(-Y(K),-\Ga_{2k-1})$ to $\shm(-Y(K),-\Ga_{2k})$.}\label{fig_odd_to_even}
\end{figure}

Now, we perform a grading shifting as follows:
$$\shm(-Y(K),-\Ga_n,S_n^{\tau(n)},i)[\sigma(n)]=\shm(-Y(K),-\Ga_n,S_n^{\tau(n)},i+\sigma(n)).$$
Here, $\tau(n)=-1$ if $n$ is even and $\tau(n)=0$ if $n$ is odd, and
$$\sigma(n)=\frac{n-1+\tau(n)}{2}.$$
We will simply write
$$\shm(-Y(K),-\Ga_n,S^{\tau}_n)[\sigma],$$
and the direct system becomes
$$...\ra\shm(-Y(K),-\Ga_n,S^{\tau}_n)[\sigma]\xra{\phi_{-,n+1}^{n}}\shm(-Y(K),-\Ga_{n+1},S^{\tau}_{n+1})[\sigma]\ra...$$
It is straightforward to check that, after the shifting, $\phi_{-,n+1}^n$ is grading preserving and $\phi_{+,n+1}^n$ shifts the grading down by 1. Thus, we conclude the following.
\bprop\label{prop_alexander_grading_on_khm}
If $S$ is a Seifert surface of $K\subset Y$, then $S$ induces a grading on $\khm(-Y,K,p)$, which we write as
$$\khm(-Y,K,p,S,i).$$
Under this grading, the map $U$ is of degree 1. 
\eprop

\bdefn
Suppose $K\subset Y$ is an oriented knot and $S$ is a Seifert surface of $K$. We can define the {\it tau invariant} $\tau(Y,K,S)$ of $K\subset Y$ with respect to $S$ as follows:
$$\tau(Y,K,S)=-\max\{i|\exists x\in\khm(Y,K,p,S,i),~U^jx\neq0~{\rm for}~{\rm any}~j\geq 0.\}$$
Here the base point can be fixed arbitrarily.
\edefn

\begin{quest}
What properties does $\tau(Y,K,S)$ have?
\end{quest}

\subsection{Basic properties}\label{subsec_properties}
\bprop\label{prop_khm_of_unknot}
Suppose $Y$ is a closed oriented $3$-manifold and $K\subset Y$ is a knot so that there exists an embedded disk $S=D^2$ with $\partial{S}=K$. Then
$$\khm(-Y,K,p)\cong \shm(-Y(1),-\delta)\otimes_{\mathcal{R}}\mathcal{R}[U].$$

Here, $p\in K$ is any choice of the base point. $(Y(1),\delta)$ is the balanced sutured manifold obtained from $Y$ by removing a $3$-ball and picking one simple closed curve on the spherical boundary as the suture.
\eprop

\bpf
First assume that $Y=S^3$, then $(Y(1),\delta)$ is a product sutured manifold and $(Y(K),\Ga_n)=(V,\ga^2_{(1,-n)})$, where $(V,\ga_{(1,-n)}^2)$ is the balanced sutured manifold as defined in Subsection \ref{subsec_sutured_solid_torus}. From Proposition \ref{prop_solid_torus_with_two_sutures}, we know that
$$\shm(-V,-\ga^2_{(1,-n)})\cong \mathcal{R}^n.$$
Suppose $S_n$ is a Seifert surface of $K$ that intersects $\Ga_n=\ga^2_{(1,-n)}$ at $2n$ points, then the argument in the proof of Proposition \ref{prop_solid_torus_with_two_sutures} can be applied to calculate the graded homology, and we conclude that: (Note $S_n$ are disks when $K$ is the unknot.)
$$\shm(-V,-\ga^2_{(1,-n)},S_n^{\tau},i)[\sigma]\cong \mathcal{R}$$
for all $i$ such that $1-n\leq i\leq 0$. Moreover, the map
$$\psi_{+,n+1}^n:\shm(-V,-\ga^2_{(1,-n)},S_n^{\tau})[\sigma]\ra\shm(-V,-\ga^2_{(1,-n-1)},S_{n+1}^{\tau})[\sigma]$$
is of degree $-1$ and is an isomorphism for all $i$ such that $1-n\leq i\leq 0$. Thus, we conclude that
$$\khm(-S^3,K,p)\cong \mathcal{R}[U].$$

When $Y$ is an arbitrary $3$-manifold, we know that
$$(Y(K),\Ga_n)=((Y(1),\delta)\sqcup (S^3(K),\ga^2_{(1,-n)}))\cup h,$$
where $h$ is a contact $1$-handle, as introduced in Baldwin and Sivek \cite{baldwin2014invariants}, which connects the two disjoint balanced sutured manifolds $((Y(1),\delta)$ and $(-S^3(K),-\ga^2_{(1,-n)})$. Thus, we know that
$$\shm(-Y(K),-\Ga_n)\cong \shm(-Y(1),-\delta)\otimes (-S^3(K),-\ga^2_{(1,-n)}).$$

Moreover, the the above isomorphism intertwines with the maps $\psi_{\pm,n+1}^n$ on $\shm(-Y(K),-\Ga_n)$ and the maps $id\otimes \psi_{\pm,n+1}^n$ on $\shm(-Y(1),-\delta)\otimes (-S^3(K),-\ga^2_{(1,-n)})$, since the corresponding contact handle attachments are clearly disjoint from each other. Thus, we conclude that
$$\khm(-Y,K,p)\cong \shm(-Y(1),-\delta)\otimes\mathcal{R}[U].$$
\epf

\bprop\label{prop_the_direct_system_stabilizes}
Suppose $K\subset Y$ is a null-homologous knot and $S$ is a minimal genus Seifert surface of $L$. Then, the direct system stabilizes: For any $i\in\intg$, if $n>g(S)-i$, then we have an isomorphism
$$\phi_{-,n+1}^n:\shm(-Y(K),-\Ga_n,S^{\tau}_n,i)[\sigma]{\cong}\shm(-Y(K),-\Ga_{n+1},S^{\tau}_{n+1},i)[\sigma].$$
\eprop
\bpf
We have the following exact triangle:
\begin{equation*}
\xymatrix{
\underline{\rm SHM}(-Y(K),-\Ga_{n+1})\ar[rr]^{\psi_{-,\infty}^{n+1}}&&\shm(-Y(K),-\Ga_{\infty})\ar[dl]^{\psi_{-,n}^{\infty}}
\\
&\shm(-Y(K),-\Ga_{n})\ar[lu]^{\psi_{-,n+1}^{n}}&
}
\end{equation*}

We prove the proposition under the assumption that $n=2k$ is even. The other case, when $n$ is odd, is similar. When $n$ is even, we know from Proposition \ref{prop_general_degree_shifting_for_knot_complement} that
$$\phi_{-,n+1}^n(\shm(-Y(K),-\Ga_n,S_{n}^-,j))\subset \shm(-Y(K),-\Ga_{n+1},S_{n+1},j).$$
By a similar argument, we have
 $$\phi_{-,\infty}^{n+1}(\shm(-Y(K),-\Ga_{n+1},S_{n+1},j))\subset \shm(-Y(K),-\Ga_{\infty},S_{\infty}^{+n},j)$$
where $S_{\infty}$ is a Seifert surface of $K$ that intersects the suture $\Ga_{\infty}$ twice. Proposition  \ref{prop_degree_shifting_formula_for_toroidal_boundary_with_two_sutures} then implies that (recall $n=2k$)
$$\shm(-Y(K),-\Ga_{\infty},S_{\infty}^{+2k},j)=\shm(-Y(K),-\Ga_{\infty},S_{\infty},j+k).$$
However, the adjunction inequality in Lemma \ref{lem_adjunction_inequality} implies that if $j+k>g(S)$, then 
$$\shm(-Y(K),-\Ga_{\infty},S_{\infty},j+k)=0.$$

Thus, for $j\in\intg$ so that $j+k>g(S)$, we have
$$\phi_{-,n+1}^n:\shm(-Y(K),-\Ga_n,S_{n}^-,j)\ra \shm(Y(K),\Ga_{n+1},S_{n+1},j)$$
is an isomorphism. From the way we perform the grading shifting in Proposition \ref{prop_alexander_grading_on_khm}, we know that, for any $j\in\intg$,
$$\shm(-Y(K),-\Ga_n,S^{\tau}_n,j)[\sigma]=\shm(-Y(K),-\Ga_n,S^{\tau}_n,j+k).$$
Thus, for the fixed grading $i\in\intg$ as in the hypothesis of the proposition, when $n>g(S)-i$, we have $(i+k)+k>g(S)$ and this implies that the map
$$\phi_{-,n+1}^n|_{\shm(-Y(K),-\Ga_n,S_{n}^{\tau},i)[\sigma]}=\phi_{-,n+1}^n|_{\shm(-Y(K),-\Ga_n,S_{n}^{\tau},i+k)}$$
is an isomorphism.
\epf

\bcor\label{cor_u_is_isomorphism}
Under the above conditions, there exists an integer $N_0$ so that, for any $i<N_0$, the $U$ map induces an isomorphism:
$$\khm(-Y,K,p,S,i)\cong \khm(-Y,K,p,S,i-1),$$
\ecor

\bpf
The proof is exactly the same as the above proposition.
\epf

\bcor\label{cor_vanishing_of_khm}
For a knot $K\subset Y$, a Seifert surface $S$ of $K$, and a fixed point $p\in K$, we have
$$\khm(-Y,K,p,S,i)=0$$
for $i>g$, and
$$\khm(-Y,K,p,S,g)\cong \underline{\rm KHM}(-Y,K,p,S,g).$$
Here, $g$ is the genus of the Seifert surface $S$, and $\underline{\rm KHM}(-Y,K,p,S,g)$ is defined in Definition \ref{defn_original_knot_Floer_homology}.
\ecor

\bpf
The first statement that
$$\khm(-Y,K,p,S,i)=0$$
for $i>g$ follows from the adjunction inequality in Lemma \ref{lem_adjunction_inequality}.  

For the second part of the statement, we prove the case where $n=2k+1$ is odd and the other case is exactly the same. Suppose $(M',\ga')$ is obtained from $(Y(K),\Ga_n)$ by a sutured manifold decomposition of $S_{n}\subset Y(K)$. It is straight forward to check that if we decompose $(Y(K),\Ga_{\infty})$ along $S_{\infty}$, then we will get exactly the same balanced sutured manifold  $(M',\ga')$. Hence, from Lemma \ref{lem_decomposition_give_top_grading} in \cite{kronheimer2010knots}, we know that
$$\shm(-Y(K),-\Ga_{n},S_{n+1},g(S)+k+1)\cong\shm(M',\ga')\cong \underline{\rm KHM}(-Y,K,p,S_{\infty},g(S)).$$
Then, the corollary follows from Proposition \ref{prop_the_direct_system_stabilizes} and the grading shifting we performed in Proposition \ref{prop_alexander_grading_on_khm}.
\epf

\begin{figure}[h]
\centering
\begin{overpic}[width=2in]{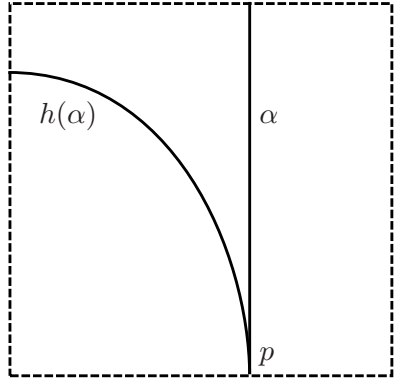}
	\put(8,67){$h(\al)$}
	\put(65,67){$\al$}
	\put(65,5){$p$}
\end{overpic}
\vspace{0.05in}
\caption{Not right-veering}\label{fig_not_right_veering}
\end{figure}

Suppose $K\subset Y$ is a fibred knot with fibre $S$ of genus $g$. Suppose $(S,h)$ is an open book corresponding to the fibration of $K\subset Y$. It supports a contact structure $\xi$ on $Y$. We call $h$ {\it not right-veering} if there is an arc $\al\subset S$ and one end point $p\in\partial{\al}$ so that near $p\subset S$, $h(\al)$ is to the left of $\al$. See Figure \ref{fig_not_right_veering}. See Baldwin and Sivek \cite{baldwin2018khovanov} for more details.

\bcor
Under the above setting, if $h$ is not right-veering, we have
$$\khm(-Y,K,p,S,g)\cong\mathcal{R},$$
and the generator is in the kernel of the $U$ map.
\ecor

\bpf
This result is the main result in Baldwin and Sivek \cite{baldwin2018khovanov}. The only difference is that we translate it into our language involving $\khm$.
\epf

\bprop\label{prop_exact_triangle_for_khm}
We have an exact triangle:
\begin{equation*}
\xymatrix{
\widehat{\rm KHM}(-Y,K,p)\ar[rr]^{U}&&\widehat{\rm KHM}(-Y,K,p)\ar[dl]^{\psi}\\
&\underline{\rm KHM}(-Y,K,p)\ar[lu]^{\psi'}&
}	
\end{equation*}
\eprop

\bpf
We will use the by-pass exact triangle
\begin{equation}\label{eq_positive_by_pass_exact_triangle}	
\xymatrix{
\underline{\rm SHM}(-Y(K),-\Ga_{n+1})\ar[rr]^{\psi_{+,\infty}^{n+1}}&&\shm(-Y(K),-\Ga_{\infty})\ar[dl]^{\psi_{+,n}^{\infty}}
\\
&\shm(-Y(K),-\Ga_{n})\ar[lu]^{\psi_{+,n}^{n+1}}&
}
\end{equation}
The maps $\{\phi_{+,n+1}^n\}_{n\in\intg_+}$ induce the $U$ map. By a similar argument, the maps $\{\phi_{+,\infty}^{n+1}\}\}_{n\in\intg_+}$ and $\{\phi_{+,n}^{\infty}\}_{n\in\intg_+}$ induce the maps $\psi$ and $\psi'$ in the statement of the proposition. Then, it is formal to check that the by-pass exact triangles (\ref{eq_positive_by_pass_exact_triangle}) for all $n\in\intg_{+}$ induce the desired one as stated in the proposition.
\epf

\subsection{Knots representing torsion classes}
In this subsection, we extend the definition of $\khm$ to the case where $K$ is not necessarily null-homologous, but represents a torsion class in $H_1(Y)$. Suppose $Y$ is a closed connected oriented $3$-manifold. Suppose further that $K\subset Y$ is an oriented knot that represents a torsion class in $H_1(Y)$. It is a basic fact that the map
$$i_*:H_1(\partial Y(K);\mathbb{Q})\ra H_1(Y(K);\mathbb{Q}),$$
which is induced by the inclusion map $i: \partial Y(K)\ra Y(K)$, has a kernel of dimension one. Thus, we can find a curve $\al\subset \partial Y(K)$ so that $\al$ bounds a properly embedded surface $S\subset Y(K)$. We always give $S$ an orientation so that $\partial{S}=\al$ is oriented in a coherent way as $K$. This surface is usually called a Rational Seifert surface of $K$. For more details, readers are referred to Ni and Vafaee \cite{ni2019null}. We still look at the knot complement $Y(K)$. On $\partial{Y(K)}\cong T^2$, there is a preferred class $\mu$ which is the meridian of $K$. There is no preferred longitude class, but we can pick any oriented non-separating simple closed curve $\lambda$ on $\partial{Y(K)}$ so that $[\mu]$ and $[\lambda]$ is an oriented basis of $H_1(\partial{Y(K)})$. Then, on $Y(K)$, we can still define the sutures $\Ga_n$ and $\Ga_{\infty}$, and there are by-pass exact triangles as in (\ref{eq_by_pass_on_knot_complement}). Note the same formula as in (\ref{eq_by_pass_on_knot_complement}) holds with our new definitions of $\Ga_n$ and $\Ga_{\infty}$. Furthermore, Corollary \ref{cor_commutative_diagram} continues to hold for exactly the same reason, so we can make the following definition.

\bdefn\label{defn_khm_for_torsion_knot}
Suppose $K\subset Y$ is a knot representing a torsion class in $H_1(Y)$ and $p\in K$ is a base point. Then, define the {\it minus version of monopole knot Floer homology}, which is denoted by $\khm(-Y,K,p)$, to be the direct limit of the direct system
$$...\ra\shm(-Y(K),\Ga_n)\xra{\psi^n_{-,n+1}}\shm(-Y(K),\Ga_{n+1})\ra...$$
Here, the maps $\psi^n_{-,n+1}$ are defined in the exact triangle (\ref{eq_by_pass_on_knot_complement}).
By Corollary \ref{cor_commutative_diagram}, the maps $\{\psi_{+,n+1}^n\}_{n\in\intg_+}$ induce a map on $\khm$, which we call $U$:
$$U:\khm(-Y,K,p)\ra\khm(-Y,K,p).$$
\edefn

It is clear that Definition \ref{defn_khm_for_torsion_knot} is independent of the choice of the longitude $\lambda$ on $\partial{Y(K)}$. Next, we want to use the rational Seifert surface $S$ of $K\subset Y$ to construct a grading on $\khm(-Y,K,p)$. As in Proposition \ref{prop_alexander_grading_on_khm}, we need to perform a grading shifting. Instead of directly writing down the value of the shift, we define the shift in an indirect way. Suppose, for any $n\in\intg_+$, $S_n$ is a rational Seifert surface of $K$, which has the minimal possible intersection with the suture $\Ga_n$. Suppose $S_n^{\tau}$ is exactly the surface $S_n$ if $|S_n\cap \Ga_n|$ is of the form $4k+2$, and $S_n^{\tau}$ is obtained from $S_n$ by performing a negative stabilization if else. We define a grading shifting, $\shm(-Y(K),-\Ga_n,S_n^{\tau})[\sigma]$, of $\shm(-Y(K),-\Ga_n,S_n^{\tau})$, so that
$$\shm(-Y(K),-\Ga_n,S_n^{\tau},i)[\sigma]=\shm(-Y(K),-\Ga_n,S_n^{\tau},i+\sigma(n)).$$
Here, the value $\sigma(n)\in\intg$ is determined by the following property: The top non-vanishing grading of $\shm(-Y(K),-\Ga_n,S_n^{\tau})[\sigma]$ equals g(S), the genus of $S$.

\brem
Note the grading shifting we performed in Proposition \ref{prop_alexander_grading_on_khm} can also be described in the above way.
\erem

\bprop\label{prop_alexander_grading_on_khm_torsion_knot}
If $S$ is a ration Seifert surface of $K\subset Y$, then $S$ induces a $\intg$-grading on $\khm(-Y,K,p)$, which we write as
$$\khm(-Y,K,p,S,i).$$
Under this grading, the map $U$ is of degree $l$, where $l$ is an integer depending on the knot $K\subset Y$. 
\eprop

As we did in Subsection \ref{subsec_properties}, we can prove that the direct system in Definition \ref{defn_khm_for_torsion_knot} stabilizes.

\bprop\label{prop_the_direct_system_stabilizes_torsion_knot}
Suppose $K\subset Y$ is a knot representing a torsion class in $H_1(Y)$, and $S$ is a rational Seifert surface of $L$. Then, the direct system stabilizes: For any $i\in\intg$, there exists $N$ so that if $n>N$, then we have an isomorphism
$$\phi_{-,n+1}^n:\shm(-Y(K),-\Ga_n,S^{\tau}_n,i)[\sigma]{\cong}\shm(-Y(K),-\Ga_{n+1},S^{\tau}_{n+1},i)[\sigma].$$
\eprop

The most common cases we might encounter a knot which represents a torsion first homology class is when performing Dehn surgeries. Suppose $K\subset Y$ is a null-homologous knot, and $S$ is a Seifert surface of $K$. Let $Y(K)$ be the knot complement. Let $\lambda$ and $\mu$ represent the longitude and meridian on $\partial{Y(K)}$, respectively, according to the framing induced $S$. We can perform a Dehn surgery along the knot $K$ and obtain a surgery manifold
$$Y_{\phi}=Y(K)\mathop{\cup}_{\phi}S^1\times D^2.$$
Suppose $\mu_{\phi}=\phi(\{1\}\times\partial{D}^2)=q_0\lambda-p_0\mu$ and $\lambda_{\phi}=\phi(S^1\times\{1\})=r_0\lambda-s_0\mu$. This results in a surgery of slope $-\frac{p_0}{q_0}$. Now $\lambda_{\phi}$ and $\mu_{\phi}$ form another framing on $\partial{Y(K)}$, so that $\mu_{\phi}$ is the meridian of the knot $K_{\phi}=S^1\times\{0\}\subset Y_{\phi}$. Note $Y(K)$ is also a knot complement of $K_{\phi}\subset Y_{\phi}$ and $K_{\phi}$ is a knot inside $Y_{\phi}$ which represents a torsion class in $H_1(Y_{\phi})$. Hence, we can use the new framing to construct a minus version of monopole knot Floer homology $\khm(-Y_{\phi},K_{\phi})$ of $(Y_{\phi},K_{\phi})$. Here, we omit the choice of base points, since the discussion will be carried out on a fixed knot complement. We have the following property.

\bprop\label{prop_surgery_formula_for_khm}
For any fixed $i_0\in\intg$, there exists $N$ so that for any surgery slope $-\frac{p_0}{q_0}<-N$, we have
$$\khm(-Y,K,S,i_0)\cong \khm(-Y_{\phi},K_{\phi},S,i_0).$$
\eprop
\bpf
We use the framing $(\lambda,\mu)$ intricately and write both the curve $q\lambda-p\mu$ or the slope $-\frac{p}{q}$ as $(q,-p)$. We use $\ga_{(q\lambda-p\mu)}$ or $\ga_{(q,-p)}$ to denote the suture consisting of two curves of slope $(q,-p)$. Note, $\ga_{(1,-n)}=\Ga_n$, for the notation $\Ga_n$ as used in Subsection \ref{subsec_construction}.

From the stabilization properties in Propositions \ref{prop_the_direct_system_stabilizes} and \ref{prop_the_direct_system_stabilizes_torsion_knot}, we know that there exists $N_1>g(S)-i_0$ such that for any $n>N_1$, we have
\begin{equation}\label{eq_identify_khm_with_shm_original_knot}
\khm(-Y,K,S,i_0)\cong \shm(-Y(K),-\ga_{(1,-n)},S^{\tau},i_0)[\sigma],
\end{equation}
and 
\begin{equation}\label{eq_identify_khm_with_shm_surgery_knot}
\khm(-Y_{\phi},K_{\phi},S,i_0)\cong \shm(-Y(K), -\ga_{(\lambda_{\phi}-n\mu_{\phi})},S^{\tau},i_0)[\sigma].
\end{equation}

Hence to prove the theorem, it is suffice to prove that for large enough $n$ and large enough surgery slope, we have
\begin{equation}\label{eq_target_grading_shifted}
\shm(-Y(K),-\ga_{(1,-n)},S^{\tau},i_0)[\sigma]\cong \shm(-Y(K), -\ga_{(\lambda_{\phi}-n\mu_{\phi})},S^{\tau},i_0)[\sigma].
\end{equation}

Fix an $n_2>N_2$, and write $\lambda_{\phi}-n_2\mu_{\phi}=q\lambda-p\mu.$ From the proof of Proposition \ref{prop_solid_torus_with_two_sutures}, we can construct two sequences of slopes $\{(q'_{j},-p'_{j})\}$ and $\{(q_j'',-p_j'')\}$ inductively as follows: Let $(q_0',-p_0')=(q,-p)$, and, for any $j\geq 1$, suppose we have the continued fraction of $(q_{j-1}',-p_{j-1}')$ to be
$$(q_{j-1}',-p_{j-1}')=[r_1,...,r_{k-j},r_{k-j+1}],$$
then define
$$(q_{j}'',-p_{j}'')=[r_1,...,r_{k-j},r_{k-j+1}+1],~(q_{j}',-p_{j}')=[r_1,...,r_{k-j}].$$
Note we identify $[r_1,...,r_{l},-1]$ as $[r_1,...,r_{l-1},r_{l}+1].$
We end the sequence when
\begin{equation}\label{eq_first_term_of_continued_fracion}
(q_{k-1}',-p_{k-1}')=[r_1]=(1,r_1).
\end{equation}
Here $r_1\leq -2$ is the first term in the continued fraction of $(q,-p)=(\lambda_{\phi}-n_2\mu_{\phi})$.

\brem
Note $(q_0,-p_0)$ is the slope of the surgery that gives rise to $(Y_{\phi},K_{\phi})$, while $(q_0',-p_0')=(q,-p)=\lambda_{\phi}-n_2\mu_{\phi}$. Also we can pick $n_2$ as large as we want.
\erem

To proceed, we only carry out the proof in the case where $n$ is odd, and for any $j$, $p_{j}'$ is odd and $p_j''$ is even. Other cases are similar. Under this assumption, we can un-package the grading shifting we performed in Propositions \ref{prop_alexander_grading_on_khm} and \ref{prop_alexander_grading_on_khm_torsion_knot}, and to prove (\ref{eq_target_grading_shifted}) is equivalent to proving (we omit the surface $S$ from the notation):
\begin{equation}\label{eq_target_grading}
\shm(-Y(K),-\ga_{(1,-n)},i_{0}+\frac{n-1}{2})\cong \shm(-Y(K), -\ga_{(q_0',-p_0')},i_{0}+\frac{p_0'-1}{2}).
\end{equation}
For $l=0,...,k-1$, write 
$$i_{l}'=i_0+\frac{p_l'-1}{2}$$

{\bf Claim 1.} There exists an $N>0$ so that if the surgery slope $-\frac{p_0}{q_0}<-N$, then, for any $l\in\{0,...,k-2\}$, there is an isomorphism:
$$\shm(-Y(K), -\ga_{(q_l',-p_l')}, S, i_l')\cong \shm(-Y(K), -\ga_{(q_{l+1}',-p_{l+1}')},S,i_{l+1}').$$

{\bf Claim 2.} There exists an $N>0$ so that if the surgery slope $-\frac{p_0}{q_0}<-N$, then we have $r_1>g(S)-i_0$.

Assuming Claim 1 and 2, we now prove the proposition. By Claim 1, Claim 2, and Proposition \ref{prop_the_direct_system_stabilizes}, we have (note we have assumed that $r_1=-p_{k-1}'$ is odd)
\beq
\shm(-Y(K),-\ga_{(1,-n)},S^{\tau},i_0)[\sigma]&\cong\shm(-Y(K),-\ga_{(1,r_1)},S^{\tau},i_0)[\sigma]\\
&=\shm(-Y(K), -\ga_{(1,-r_1)}, S, i_{0}+\frac{-r_1-1}{2})\\
&=\shm(-Y(K), -\ga_{(q_{k-1}',-p_{k-1}')}, S, i_{k-1}')\\
&\cong \shm(-Y(K), -\ga_{(q_0',-p_0')}, S, i_0')\\
&= \shm(-Y(K), -\ga_{(q_0',-p_0')},i_{0}+\frac{p_0'-1}{2}).
\eeq
Thus (\ref{eq_target_grading}) is proved, and Proposition \ref{prop_surgery_formula_for_khm} follows.

To prove Claim 2, by definition, we have
\begin{equation}\label{eq_r_1}
r_1=-(\lfloor\frac{p}{q}\rfloor+1)~{\rm and}~\frac{p}{q}=\frac{s_0+n_2p_0}{r_0+n_2q_0}.
\end{equation}
If we choose large enough $n_2$ (we can freely make $n_2$ larger), then we know that
\begin{equation}\label{eq_p_over_q}
\lfloor\frac{p}{q}\rfloor\geq \lfloor\frac{p_0}{q_0}\rfloor-1.
\end{equation}
Hence, for any surgery slope $-\frac{p_0}{q_0}<N=-(g(S)-i_0)$, Claim 2 holds.

It remains to prove Claim 1. As in Subsection \ref{subsec_sutured_solid_torus}, the sutures of slopes $(q_l',-p_l')$ and $(q_l'',-p_l'')$ fit into a by-pass exact triangle:
\begin{equation}\label{eq_exact_triangle_on_knot_complement}
	\xymatrix{
	&\shm(-Y(K),-\ga_{(q_{l-1}',p_{l-1}')})\ar[rd]^{\psi_{l,2}}&\\
	\shm(-Y(K),-\ga_{(q_{l}'',p_{l}'')})\ar[ru]^{\psi_{l,1}}&&\shm(-Y(K),-\ga_{(q_{l}',p_{l}')})\ar[ll]^{\psi_{l,0}}
	}
\end{equation}
If $Y=S^3$ and $K$ is the unknot, then $\psi_{j,k}=\psi_{-,k}$ for $k=0,1,2$ in the previous exact triangle (\ref{eq_exact_triangle_for_two_sutures}). As in Subsection \ref{subsec_sutured_solid_torus}, for all $l\in\{1,...,k-1\}$ and $j\in\intg$, we have
$$\psi_{l,0}:(\shm(-Y(K),-\ga_{(q_{l}',-p_{l}')},S^{-p_l''},j)\ra \shm(-Y(K),-\ga_{(q_{l}'',-p_{l}'')},S^{+p_l'},j),$$
$$\psi_{l,1}:\shm(-Y(K),-\ga_{(q_{l}'',-p_{l}'')},S^{+p_l'},j)\ra\shm(-Y(K),-\ga_{(q_{l-1}',-p_{l-1}')},S,j),$$
$$\psi_{l,2}:\shm(-Y(K),-\ga_{(q_{l-1}',-p_{l-1}')},S,j)\ra(\shm(-Y(K),-\ga_{(q_{l}',-p_{l}')},S^{-p_l''},j).$$
Note, in above formulae, we have assume that $p_{l-1}''$ is odd for all $l$. 
From them, Claim 1 is equivalent to the fact that $\psi_{l_2}$ is an isomorphism at the grading 
$$j=i_{l-1}'=i_0+\frac{p_{l-1}'-1}{2},$$
which is further equivalent to that
\begin{equation}\label{eq_target_vanishing_grading_shifted}
\shm(-Y(K),-\ga_{(q_{l}'',-p_{l}'')},S^{+p_l'},i_{l-1}')=0.
\end{equation}

Note, by assumption, $p_l''=p_{l-1}'-p_l'$ is even. From the grading shifting property, Proposition \ref{prop_degree_shifting_formula_for_toroidal_boundary_with_two_sutures}, we know that (\ref{eq_target_vanishing_grading_shifted}) is equivalent to
\begin{equation}\label{eq_target_vanishing_grading}
	\shm(-Y(K),-\ga_{(q_{l}'',-p_{l}'')},S^{+},i_{l-1}'+\frac{p_l'-1}{2})=0.
\end{equation}

Note we have $|\partial S^+\cap\ga_{(q_{l}'',-p_{l}'')}|=2p_l''+2$. From (the vanishing statement of) Lemma \ref{lem_decomposition_give_top_grading}, we know that (\ref{eq_target_vanishing_grading}) happens if
\begin{equation}\label{eq_target_inequality_on_grading}
i_{l-1}'+\frac{p_l'-1}{2}>g(S)+\frac{p_l''}{2}.
\end{equation}
Recall that
$$i_{l-1}'=i_0+\frac{p_{l-1}'-1}{2},$$
so we know that
\beq
&i_{l-1}'+\frac{p_l'-1}{2}>g(S)+\frac{p_l''}{2}\\
\Leftrightarrow&i_0+\frac{p_{l-1}'-1}{2}+\frac{p_l'-1}{2}>g(S)+\frac{p_l''}{2}\\
\Leftrightarrow&p_l'>g(S)-i_0+1.
\eeq
Since, by (\ref{eq_r_1}) and (\ref{eq_p_over_q}), we have
$$p_l'\geq p_{k-1}'=-r_1\geq \lfloor\frac{p_0}{q_0}\rfloor>N.$$
Thus, if we pick $N=-(g(S)-i_0)$, then (\ref{eq_target_inequality_on_grading}) holds and Claim 1 follows. This concludes the proof of Proposition \ref{prop_surgery_formula_for_khm}.
\epf

\section{Instantons and knot Floer homology}
\subsection{Instanton Floer homology and generalized eigenspace decompositions}
Suppose $Y$ is a closed connected oriented $3$-manifold, and $\omega$ is a fixed Hermitian line bundle whose first Chern class $c_1(\omega)$ has an odd pairing with the fundamental class of some surface. Suppose further that $E$ is an $U(2)$-bundle whose determinant line bundle $\Lambda^2E$ is isomorphic to $\omega$. Let $\mathfrak{g}_E$ be the bundle of traceless skew-Hermitian endomorphisms of $E$, and let $\mathcal{A}_{E}$ be the space of $SO(3)$-connections on $\mathfrak{g}_E$. Let $\mathcal{G}_E$ be the group of determinant-one gauge transformations and let $\mathcal{B}_E=\mathcal{A}_E\slash \mathcal{G}_E$. Then, we can use the Chern-Simons functional to construct a well defined $SO(3)$ instanton Floer homology over $\mathbb{C}$, which we denote by $I^{\omega}(Y)$.

If $x\in Y$ is a point, then there is an action $\mu(x)$ on $I^{\omega}(Y)$. The action $\mu(x)$ has eigenvalues $2$ and $-2$. By slightly abusing the notations, from now on we use $I^{\omega}(Y)$ to denote only the generalized eigenspace of $\mu(x)$ corresponding to eigenvalue 2.

Suppose $\Sigma\subset Y$ is a closed oriented embedded surface inside $Y$. Then, there is also an action $\mu(\Sigma)$ on $I^{\omega}(Y)$. We have the following result about the eigenvalues:

\bprop[Kronheimer, Mrowka, \cite{kronheimer2010knots}]\label{prop_possible_generalized_eigenvalues}
If $c_1(\omega)$ and $\Sigma$ has an odd pairing, then the eigenvalues of the action $\mu(K)$ on $I^{\omega}(Y)$ belongs to the set of even integers ranged from $2-2g(\Sigma)$ to $2g(\Sigma)-2$.
\eprop

If $\Sigma$ and $\Sigma'$ are two such embedded surfaces, then the action $\mu(\Sigma)$ and $\mu(\Sigma')$ commute. Then, we can look at the simultaneous generalized eigenspace. Similar to Corollary 7.6 in Kronheimer and Mrowka \cite{kronheimer2010knots}, we can make the following definition.

\bdefn\label{defn_generalized_eigenspace_decomposition}
Suppose we have a linear function $\lambda: H_2(Y;\intg)\ra 2\intg$, then we can define
$$I^{\omega}(Y)_{\lambda}=\bigcap_{\sigma\in H_2(Y;Z)}\bigcup_{N\geq0}ker(\mu(\sigma)-\lambda(\sigma))^N.$$
Such a function $\lambda$ is a called an {\it eigenvalue function}.
\edefn

If the embedded surface $\Sigma$ represents a zero class in $H_2(Y;\mathbb{Q})$, then the action $\mu(\Sigma)$ is trivial. This means that if $I^{\omega}(Y)_{\lambda}\neq 0$ then we can lift $\lambda$ to a linear map (which we will use the same notation to denote)
$$\lambda: H_2(Y;\mathbb{Q})\ra\mathbb{Q}.$$
Thus, from now on, we regard $\lambda$ as an element in $H^2(Y;\mathbb{Q})$. We then have a decomposition
$$I^{\omega}(Y)=\bigoplus_{\lambda\in H^2(Y;\mathbb{Q})}I^{\omega}(Y)_{\lambda}.$$

Suppose $R\subset Y$ is a closed oriented embedded surface inside $Y$, then as we did in Definition \ref{defn_set_of_spin_c_structures}, we can define the following.
\bdefn
Suppose the pair $(Y,R)$ is as above. Then, we can define the set
$$\mathfrak{H}^*(Y|R)=\{\lambda\in H^2(Y;\mathbb{Q})|\lambda([R])=2g(R)-2,~I^{\omega}(Y)_{\lambda}\neq0\},$$
The elements $\lambda\in \mathfrak{H}^*(Y|R)$ are called {\it supporting eigenspace functions.}
\edefn

We have the following lemma which is the instanton correspondence to Lemma \ref{lem_spin_c_structure_shall_extend} for monopole theory.
\blem\label{lem_eigenvalue_funciton_shall_extend}
Suppose $(W,\nu)$ is a cobordism between $(Y,\omega)$ and $(Y',\omega')$. Suppose further that $\lambda\in H^2(Y;\mathbb{Q})$ and $\lambda'\in H^2(Y';\mathbb{Q})$ are two eigenvalue functions. Let $i:Y\ra W$ and $i':Y'\ra W$ are the inclusion map. If
$$I(W,\nu)(I^{\omega}(Y)_{\lambda})\cap I^{\omega'}(Y')_{\lambda'}\neq\{0\},$$
then there must be an element $\tau\in H^2(W;\mathbb{Q})$ so that $i^*(\tau)=\lambda$ and $(i')^*(\tau)=\lambda'$. 
\elem

\bpf
For a second homology class $\sigma$ and a rational number $r\in\mathbb{Q}$ we can define
$$I^{\omega}(Y,\sigma,r)=\bigcup_{N\geq0}ker(\mu(\sigma)-r)^N.$$
By definition, we know that
$$I^{\omega}(Y)_{\lambda}=\bigcap_{\sigma\in H_2(Y;\mathbb{Q})}I^{\omega}(Y,\sigma,\lambda(\sigma)).$$
Similarly, we can define $I^{\omega'}(Y',\sigma',r')$.

Note we can regard an element $\tau\in H^2(W;\mathbb{Q})$ as a map
$$\tau:H_2(W;\mathbb{Q})\ra \mathbb{Q}.$$
Suppose there are no such $\tau$ as in the statement of the lemma, then there is a class $\sigma_0\in H_2(Y;\mathbb{Q})$ and a class $\sigma_0'\in H_2(Y';\mathbb{Q})$ so that
$$i_*(\sigma_0)=i'_*(\sigma'_0)\in H_2(W),$$
while
$$\lambda(\sigma_0)\neq\lambda'(\sigma'_0).$$
 
Thus, we know that
$$I(W,\nu)(I^{\omega}(Y)_{\lambda})\subset I(W,\nu)(I^{\omega}(Y,\sigma_0,\lambda(\sigma_0)))\subset I^{\omega'}(Y',\sigma'_0,\lambda(\sigma_0)).$$
The last inclusion follows from Lemma 2.6 in \cite{baldwin2018khovanov}. However, $\lambda(\sigma)\neq\lambda'(\sigma')$ so
$$I^{\omega'}(Y',\sigma'_0,\lambda(\sigma_0))\cap I^{\omega'}(Y',\sigma'_0,\lambda'(\sigma'_0))=\{0\}.$$
Hence, we conclude
$$I(W,\nu)(I^{\omega}(Y)_{\lambda})\cap I^{\omega'}(Y')_{\lambda'}=\{0\},$$
which is a contradiction. Thus Lemma \ref{lem_eigenvalue_funciton_shall_extend} follows.
\epf

\subsection{Sutured instanton Floer homology}
Suppose $(M,\ga)$ is a balanced sutured manifold, then, as we did for monopole theory, we can construct a closure of $(M,\ga)$ and apply the construction of instanton Floer homology in the previous subsection. Pick a connected auxiliary surface $T$ of large enough genus, then we can get a pre-closure
$$\widetilde{M}=M\cup T\times [-1,1],~{\rm with}~\partial {\widetilde{M}}=R_+\sqcup R_-.$$
For the construction in instanton theory, we also need to pick a point $p\in T$ so that there are corresponding points $p_{\pm}\in R_{\pm}$. When choosing the gluing diffeomorphism $h:R_+\ra R_-$, we also require that $h(p_+)=p_-$. Thus, we know that, inside the closure $(Y,R)$, there is a closed curve $p\times S^1\subset Y$. Let $\omega$ be a complex line bundle over $Y$ whose first Chern class is dual to the curve $p\times S^1$. Then, we can make the following definition.

\bdefn[Kronheimer, Mrowka \cite{kronheimer2010knots}]
Define the {\it sutured instanton Floer homology} of $(M,\ga)$ to be
$$SHI(M,\ga)=I^{\omega}(Y|R)=\bigoplus_{\lambda\in\mathfrak{H}^*(Y|R)}I^{\omega}(Y)_{\lambda}.$$
\edefn

Baldwin and Sivek \cite{baldwin2015naturality} also made refinements of closures and constructed canonical maps for the sutured instanton Floer homology.

\bdefn
A marked odd closure $\mathcal{D}=(Y,R,r,m,\eta,\al)$ of $(M,\ga)$ is a tuple so that $(Y,R,r,m,\eta)$ is a marked closure of $(M,\ga)$ as in definition \ref{defn_marked_closure}, the simple closed curve $\al$ is disjoint from ${\rm im}(m)$, and $\al\cap r(R\times[-1,1])$ is of the form $r(p\times[-1,1])$.

We can pick a complex line bundle $\omega$ whose first Chern class is dual to $\al\sqcup \eta$. Then we can define
$$SHI(\mathcal{D})=I^{\omega}(Y|r(R\times\{0\})).$$
\edefn 

\bthm[Baldwin, Sivek \cite{baldwin2015naturality}]\label{thm_instanton_canonical_maps}
Suppose $(M,\ga)$ is a balanced sutured manifold, and $\mathcal{D}$ and $\mathcal{D}'$ are two marked odd closures of $(M,\ga)$. Then, there is a canonical map 
$$\Phi_{\mathcal{D},\mathcal{D}'}:SHI(\mathcal{D})\ra SHM(\mathcal{D}'),$$
which is an isomorphism well defined up to multiplication by a non-zero element in $\mathbb{C}$. Furthermore, the canonical map satisfies the same functoriality properties as those of the canonical maps for sutured monopole Floer homology in Theorem \ref{thm_canonical_maps}.
\ethm

Hence, we have a well defined projective transitive system
$$\shi(M,\ga)$$
associated to $(M,\ga)$. For a knot, there is a similar discussion as in Subsection \ref{subsec_naturality} and we have a well defined projective transitive system
$$\khi(Y,K,p)$$
associates to a triple $(Y,K,p)$ for a knot $K\subset Y$ and a base point $p\in K$.

There are similar results for the contact gluing maps and by-pass exact triangles.

\bthm[Li \cite{li2018gluing}]\label{thm_instanton_gluing_maps}
There is a gluing map for sutured instanton Floer homology, satisfying the same properties as in Theorem \ref{thm_gluing_map}.
\ethm

\bthm[Baldwin and Sivek \cite{baldwin2018khovanov}]\label{thm_instanton_bypass_exact_triangle}
Suppose $(M,\ga_1)$, $(M,\ga_2)$ and $(M,\ga_3)$ are three balanced sutured manifolds which are related in the same way as in theorem \ref{thm_by_pass_attachment}. Then there is still a by-pass exact triangle
\begin{equation*}
\xymatrix{
\shi(-M,-\ga_1)\ar[rr]^{\psi_{12}}&&\shi(-M,-\ga_2)\ar[dl]^{\psi_{23}}\\
&\shi(-M,-\ga_3)\ar[ul]^{\psi_{31}}&
}
\end{equation*}
where the maps $\psi_{ij}$ comes from the gluing maps in sutured instanton Floer homology, just as the monopole case in Subsection \ref{subsec_contact_elements_and_contact_structures}.
\ethm

\subsection{Statement of results}
With Lemma \ref{lem_eigenvalue_funciton_shall_extend} and Theorem \ref{thm_instanton_bypass_exact_triangle} in place of Lemma \ref{lem_spin_c_structure_shall_extend} and Theorem \ref{thm_by_pass_attachment}, we can recover all results we obtained in this paper for sutured monopole Floer homology. We present those results without further proofs.

\bprop
Suppose $(M,\ga)$ is a balanced sutured manifold and $\mathcal{D}$, and $\mathcal{D}'$ are two marked odd closures of the same genus. Then, the canonical map $\Phi_{\mathcal{D},\mathcal{D}'}$ in sutured instanton Floer theory can be interpreted in terms of the Floer excision cobordism, in the same way as in Proposition \ref{prop_new_definition_for_canonical_map}.
\eprop

\bthm
Suppose $(M,\ga)$ is a balanced sutured manifold, and $S$ is a properly embedded surface inside $M$ so that $\partial{S}$ is connected and $|\partial{S}\cap \ga|=2n$ with $n$ odd. Then, $S$ induces a grading on $\shi(M,\ga)$ which we denote by
$$\shi(M,\ga,S,i).$$
\ethm

\bprop
Suppose $(M,\ga)$ is a balanced sutured manifold so that $M$ is the complement of a null-homologous knot $K\subset X$ and $\ga$ has two components. Suppose further that $S$ is a Seifert surface of $K$, viewed as a properly embedded surface in $M$, so that $|\partial{S}\cap \ga|=2n$. Then, for any $p,l,k\in \intg$ such that $n+p$ is odd, we have
$$\shi(-M,-\ga,S^{p},l)=\shi(-M,-\ga,S^{p+2k},l-k).$$
\eprop

\bprop
Suppose $V$ is a solid torus and $\ga$ is a suture on $\partial{V}$ with $2n$ components and slope $\frac{p}{q}$, then
$$\shi(-V,-\ga)\cong\mathbb{C}^{(2^{n-1}\cdot |p|)}.$$
\eprop

\bthm
Suppose $K$ is a null-homologous knot inside an closed connected oriented $3$-manifold $Y$ and $p\in K$ is a base point. Then, there is a projective $\mathbb{C}$-vector space $\khi(Y,K,p)$, whose elements are well defined up to multiplication by a non-zero element in $\mathbb{C}$, associated to the triple $(Y,K,p)$. Also, there is a homomorphism
$$U:\khi(Y,K,p)\ra\khi(Y,K,p).$$
If $S$ is a Seifert surface of $K$, then $S$ induces a $\intg$ grading on $\khi(Y,K,p)$ so that $U$ is of degree $-1$. Furthermore, analogous results to Proposition \ref{prop_khm_of_unknot}, Proposition \ref{prop_the_direct_system_stabilizes}, Corollary \ref{cor_u_is_isomorphism}, Corollary \ref{cor_vanishing_of_khm}, Proposition \ref{prop_exact_triangle_for_khm}, and Proposition \ref{prop_surgery_formula_for_khm} all hold for $\khi(Y,K,p)$. 
\ethm

\bibliography{Index}
\end{document}